\documentclass[12pt,preprint]{elsarticle}

\usepackage{amssymb}

\usepackage[utf8]{inputenc}
\usepackage{graphicx}
\usepackage{dcolumn}
\usepackage{dsfont}
\usepackage{amsthm}
\usepackage{epstopdf}
\usepackage{mathtools}
\usepackage{textcomp}
\usepackage{colortbl}
\usepackage{float}
\usepackage{bm}
\usepackage[format=plain,font={small},justification=centerlast]{caption}
\usepackage{subcaption}
\usepackage{color}
\usepackage{tabularx,ragged2e,booktabs}
\usepackage{MnSymbol}
\usepackage{wasysym}
\usepackage{anysize}
\usepackage{epstopdf}
\usepackage{verbatim}
\usepackage{tikz}
\usepackage{pgfplots}
\usepackage{mathtools}
\usepackage{natbib}
\usepackage{xcolor}
\DeclareGraphicsExtensions{.pdf,.png,.jpg,.pdf}



\biboptions{comma,square}

\def\im{\mathop{\rm \od{\iota}}\nolimits}
\newcommand{\BR}{{\mathbb{R}}}
\newcommand{\p}{\partial}

\newcommand{\ts}[1]{\textstyle #1}
\newcommand{\bn}[1]{\mbox{\boldmath $#1$}}
\newcommand{\bc}{\begin{center}}
	\newcommand{\ec}{\end{center}}
\newcommand{\be}{\begin{equation}}
\newcommand{\ee}{\end{equation}}
\newcommand{\bea}{\begin{eqnarray}}
\newcommand{\eea}{\end{eqnarray}}
\newcommand{\ba}{\begin{array}}
	\newcommand{\ea}{\end{array}}

\newcommand{\JJ}{\mathbf{J}}
\newcommand{\II}{\mathbf{I}}

\newcommand{\N}{\mathds{N}}
\newcommand{\x}{\mathbf{x}}
\newcommand{\n}{\mathbf{n}}
\newcommand{\e}{\mathbf{e}}
\newcommand{\y}{\mathbf{y}}
\newcommand{\ds}{\displaystyle}
\renewcommand{\a}{\mathbf{a}}
\renewcommand{\r}{\mathbf{r}}
\newcommand{\h}{\mathbf{h}}
\renewcommand{\v}{\mathbf{v}}
\renewcommand{\u}{\mathbf{u}}
\newcommand{\bO}{\mathcal{O}}
\newcommand{\Pcal}{\mathcal{P}}

\newcommand{\Z}{\mathds{Z}}
\newcommand{\Q}{\mathds{Q}}
\newcommand{\R}{\mathds{R}}
\newcommand{\C}{\mathds{C}}
\renewcommand{\P}{\mathds{P}}
\newcommand{\phiep}{\varphi_{\epsilon}}
\newcommand{\mt}{\mathcal{T}}
\newcommand{\mb}{\mathcal{B}}
\newcommand{\ml}{\mathcal{L}}
\newcommand{\ve}{\varepsilon}
\newcommand{\om}{\Omega}
\newcommand{\td}{\tilde}
\newcommand{\al}{\alpha}
\newcommand{\vp}{\varphi}
\newcommand{\rd}{{\mathbb{R}^{d}}}
\newcommand{\de}{\delta}

\newcommand{\ri}{\mathrm{i}}

\newtheorem{rmrk}{Remark}

\pgfplotsset{compat=1.14}
\journal{Elsevier}

\begin{document}
	
	\begin{frontmatter}
		
		
		\title{Radial basis function-generated finite differences with Bessel weights for the 2$D$ Helmholtz equation}


		
		\author{Mauricio A. Londo\~{n}o}
		\ead{alejandro.londono@udea.edu.co}
		\author{Hebert Montegranario}
		\ead{hebert.montegranario@udea.edu.co}
		
		\address{Instituto de Matemáticas \\
			Universidad de Antioquia\\
			Calle 67 53-108, Medell\'in, Colombia}
		
		\begin{abstract}
			In this paper we obtain approximated numerical solutions for the 2D Helmholtz equation using a Radial Basis Function-generated Finite Difference (RBF-FD) scheme, where  weights are calculated by applying an oscillatory radial basis function given in terms of Bessel functions of the first kind. The problem of obtaining  weights by local interpolation is ill-conditioned;  we overcome this difficulty by means of regularization of the interpolation matrix by perturbing its diagonal. The condition number of this perturbed matrix is controlled according to a prescribed value of a regularization parameter. Different numerical tests are  performed  in order to study  convergence and  algorithmic complexity. As a result, we verify that dispersion and  pollution effects are mitigated. 
		\end{abstract}
		
		
	\end{frontmatter}
	
	
	\section{Introduction}
    
    The Helmholtz equation is an elliptic Partial Differential Equation (PDE) which represents time-independent solutions of the wave equation. This equation models a wide variety of physical phenomena. These include among others, acoustic wave scattering, time harmonic acoustic, electromagnetic fields, water wave propagation, membrane vibration and radar scattering. One of the objectives of  numerical solutions of Helmholtz equation is to build a solver dealing with (i) a wide range of wave numbers and (ii) decrease the  accumulation of spurious dispersion in computation due to the pollution effect. Given that an increase in wave number requires an appropriate
increase in the mesh resolution for maintaining the level of accuracy, most numerical methods face difficulties for tackling the pollution effect.

In this work we consider to find numerical solutions for the 2$D$ Helmholtz equation given by
    
    \begin{equation}\label{eq:helmholtz1}
    \left\{
    \begin{array}{rcll}
    -\Delta u(\x)-\omega^2c(\x)^{-2}u(\x)&=&f(\x), &\mbox{ in } \Omega\\
    b\frac{\partial}{\partial\n}u(\x)+\ri\omega c(\x)^{-1}\mathcal{B}u(\x)&=&g(\x), & \mbox{ on } \Gamma=\partial\Omega
    \end{array}
    \right.
    \end{equation}
where $\omega$ is the angular frequency, $c(\x)>0$ is the sound speed of the continuous media, $f(\x)$ is the source term, $\n$ is unitary normal vector to the boundary $\Gamma$, $b$ takes values zero or one, $\mathcal{B}$ is a certain linear operator and $g(\x)$ is certain exact data on $\Gamma$. Here $\ri=\sqrt{-1}$.

We apply a basic idea inspired in a combination of Trefftz method and Radial basis functions (RBF). Trefftz method consider  linear combinations of solutions of the equation itself to solve a PDE. In this case the solutions we apply are given by oscillatory radial basis functions in terms of Bessel functions of the first kind. The weights of the linear combination are found by considering local stencils in the way of finite difference method(FD). This joint formulation of RBF and FD is well known as RBF-FD method. RBF-FD methods have been widely applied in the solution of partial differential equations. For a better idea the reader may consult \cite{fassbook,fornb1} and references therein; in particular, Fornberg \cite{FORNBERG2006_oscillatory_RBF} studied the family of oscillatory Bessel RBF's  we use here.

When solving a differential equation is very important to apply methods adaptable to the geometry or node distribution in the solution domain.  RBFs are an appropriate meshless tool, given that they only depend on the distance between points, do not require any prescribed structure on them. This method consider solutions in the form
$S(\x)=\sum_{j=1}^n\varphi(\|\x-\x_j\|) + p(\x)$; where $\varphi(r)$ is a radial function such as the Gaussian family of radial basis functions (GRBF) $\varphi(r)=e^{-(\ve r)^2}$, $\|\cdot\|$ is the Euclidean norm and $X=\{\x_1,\x_2,\dots,\x_n\}$ is a set of scattered points in the domain $\Omega$. Depending of the chosen RBF the low degree polynomial term $p(\x)$ could be or not included.

From the first publications of Kansa \cite{kansa1,kansa2} until now, there has been an increasing interest and success in these methods with a wide range of applications  \cite{fassbook,yangu1,yangu2}. In particular, several RBFs methods have been developed and applied to obtain numerical solutions of PDEs. There exist  a large number of these functions, in fact a theorem  from Bochner (1932) \cite{iske1} shows that (under certain conditions)  if  the Fourier transform
\(\widehat \Phi\) of $\Phi$ is positive on $\rd$ with $\widehat{\Phi}>0$ then $\Phi$ is positive definite in $\rd$. Nevertheless, only  some of them are usually chosen, depending on every particular application.

Usually, RBFs as multiquadrics $\phi(r)=\sqrt{1+(\ve r)^2}$  contain a shape parameter $\ve$ which decides the flatness of the function and by consequence the  condition of the interpolation matrix. As $\ve \rightarrow  0$  the shape of $\phi (r)$ goes from very peaked($\ve$ large)to nearly flat ($\ve$ small), in this last case the interpolation has shown to be remarkably accurate \cite{fornb1}. Until recently, the literature only have shown non-oscillatory RBFs, nevertheless   the family of Bessel oscillatory RBFs  \eqref{eq:bessel_RBF} applied here provides existence and uniqueness of the interpolation problem, they do not diverge in the limit of flat basis functions for any node geometry and have exact polynomial reproduction of arbitrary order \cite{flyer1}. In this paper the role of the shape parameter is taken by the wave number $k$. The main argument for applying these functions in Helmholtz problem is that they are in themselves solutions of the equation.

	The rest of the paper is organized as follows. In the next section we define the oscillatory Bessel function we are going to  work with. In section 3 we give the general setting of the RBF-FD method. In Section 4 we describe the method of diagonal increments which deals with the ill-conditioning of the interpolation matrix.  Sections 5,6  show results obtained in testing our methods with some well-known problems and benchmarks of current literature on Helmholtz equation.
	
    \section{Oscillatory RBF}

There exist a wide number of Trefftz methods for the Helmholtz problems that have been surveyed in   \cite{trefftz_paper}. These are schemes of type finite elements where test and trial functions are local solutions of the differential equation to solve. Inspired by Trefftz methods, in this paper we work with a family of oscillatory RBF's whose members are solutions of the homogeneous Helmholtz equation. Besides,  given the oscillatory behavior of solutions of Helmholtz equation, it makes sense to consider  such a family, whose members are given in terms of Bessel functions of the first kind.

The oscillatory RBF class $\varphi_k^{(d)}(r)$ (BRBF) is the family of radial basis functions given by
\begin{equation}\label{eq:bessel_RBF}
\varphi_k^{(d)}(r)=\frac{J_{d/2-1}(kr)}{(kr)^{d/2-1}}, \ \ \ d=1,2,\ldots,
\end{equation}
 which  are detailed studied in \cite{FORNBERG2006_oscillatory_RBF}. Here $J_{\alpha}(r)$ is denoting the Bessel function of the first kind and order $\alpha$. Two remarkable properties of these functions  are:
\begin{itemize}
	\item the non-singularity of the interpolation matrix for arbitrarily scattered data in up to $d$ dimensions, when $d>1$,

	\item and that the Laplace eigenvalue problem $\Delta \varphi + k^2\varphi =0$ has as bounded solutions, at the origin, the functions given in \eqref{eq:bessel_RBF}, thus any interpolant of the form
	\begin{equation}\label{eq:oscilla_interpolant}
	S(\x)=\sum_{j=1}^n\alpha_j\varphi^{(d)}_k(\|\x-\x_j\|)
	\end{equation}
	will also satisfy  $\Delta S+k^2 S=0$.  

\end{itemize}

In the case $d=2$ the current literature shows very few applications of the oscillatory RBF \eqref{eq:bessel_RBF}; this has been because the function \eqref{eq:oscilla_interpolant} implies that $\Delta S=-k^2S$ so, by the weak maximum principle, the function $S$ cannot have local maximum at points where for some neighborhood it is negative; a fact that put restrictions to be used for general 2D interpolation. But in this work such a feature becomes a strength, since we are just interpolating solutions of Helmholtz problems, which locally can be seen as plane waves satisfying the homogeneous Helmholtz equation. In early works, as in \cite{LIN2012_oscillatory_radial_Basis_function_helmholtz}, oscillatory RBF's based on Bessel functions have been employed to solve the 2D Helmholtz equation with constant wavenumber within the approach of global collocation method and using the RBF  
\begin{equation}\label{eq:RBFlin2012}
\phi_{C,k}(r)=J_0(k\sqrt{r^2+C^2}),
\end{equation}
which has two shape parameters with $k$ corresponding to the wavenumber and $C$  is  empirically chosen. The ill-conditioning of the interpolation matrix that arises from \eqref{eq:RBFlin2012} is overcome by way of a regularized singular value decomposition method. 

For our interest, the 2D Helmholtz problem  with large wavenumber, we take the special case $d=2$. So we work with the oscillatory RBF
\begin{equation}\label{eq:bessel_RBF2D}
\phi_{k}(r)=J_0(kr),
\end{equation}
such that in the approach RBF-FD the shape parameter $k$ will be evaluated at the wave number $k(\x)=\omega/c(\x)$ corresponding to the center of the stencil.

Among the strengths of the oscillatory RBF family $\varphi_k^{(d)}(r)$ over other radial functions we must remark that the Gaussian family is contained in the Bessel RBF class in the limiting case
\begin{equation*}
\lim_{\de\to\infty}2^{\de}\de ! \frac{J_{\de}(2\sqrt{\de }r)}{2\sqrt{\de}r}=e^{-r^2}
\end{equation*}
In fact, all other RBFs could suffer divergence when  $\ve\to 0$. It was shown in \cite{gibbs} that such divergence can never occur when using GRBF, independently of the node distribution.\\

It is well known that for assembling the sparse matrix, which discretizes the Helmholtz problem is necessary to solve a small linear equation system at each node.  As it will be seen, interpolation matrices are ill-conditioned and we  deal with this issue by the Method of Diagonal Increments (MDI) \cite{method_of_diagonal_increments}, \cite{SARRA2014_regularized_positive_definite} adding to the diagonal entries a small regularization parameter $\beta>0$, thus we solve, instead of the linear system  $\mathbf{b}=\mathbf{A}\mathbf{y}$,  the equation
\begin{equation}
\mathbf{b}=(\mathbf{A}+\beta\mathbf{I})\widetilde{\mathbf{y}},
\end{equation}
where $\II$ is the identity matrix. An explanation of MDI will be given in the section \ref{sec:MDI}, where it is  shown that the matrix $\widetilde{\mathbf{A}}=\mathbf{A}+\beta\II$ is better conditioned than $\mathbf{A}$ and  $\widetilde{\mathbf{y}}\approx\mathbf{y}$. Now are  described the properties of discretizing Helmholtz problems with Bessel RFB.\\


\section{Discretization by RBF-FD method }

Under the RBF interpolation framework, we want to approximate the solution of boundary value problems in the form\footnote{Helmholtz problems we are dealing with in this paper, can be seen as particular cases of \eqref{eq:boundary_pro}.}
\begin{equation}\label{eq:boundary_pro}
\begin{cases}
\mathcal{L} u(\x)=f(\x), & \mbox{ \ if \ } \x\in\Omega\\
\mathcal{B} u(\x)=g(\x), & \mbox{ \ if \ } \x\in\partial\Omega,
\end{cases}
\end{equation}
where $\ml$ and  $\mb$ are linear partial differential operators whose coefficients have a  good enough regularity, and $\Omega$ is a bounded, open and connected set in $\R^d$. It is assumed that (\ref{eq:boundary_pro}) it is a well-posed problem.

In interpolation with Radial Basis Functions the goal is to reconstruct a real-valued or complex-valued function $u$ defined on a bounded domain $\Omega\subset\R^d$ from the values $u(\x_k)$ of $u$ on a finite set of $N$ scattered nodes $ X = \{\x_1,\x_2, \ldots , \x_N\}\subset \Omega \subset \R^{d}$, where $d$ is a positive integer. A radial basis function with \emph{shape parameter} $\ve$ is defined as a function $\Phi_{\ve}:\R^d\times\R^d\rightarrow\R$  such that $\Phi_{\ve}(\x,\y)=\phi(\ve \|\x-\y\|)$, where $\phi:[0,\infty)\rightarrow\R$ is a single variable function \cite{wendland1,fassbook,schaback}. A sufficiently smooth function $u:\Omega\subset\R^d\rightarrow \R $, with $\Omega$ an open set whose boundary is regular enough, can be approximated by the \emph{interpolant}
\begin{equation} \label{eq:pro_app}	\mathcal{P}_{X,\ve}u(\x)=\sum_{j=1}^N\alpha_j\Phi_{\ve}(\x,\x_j),
\end{equation}
by forcing the condition $\mathcal{P}_{X,\ve}u(\x_k)=u(\x_k)$ for $k=1,\ldots,N$, where weights $\alpha_j$ can be determined by solving the linear system 
\begin{equation} \label{eq:pro_app_forcing}
u(\x_k)=\sum_{j=1}^N\alpha_j\Phi_{\ve}(\x_k,\x_j), \mbox{ with } k=1,\ldots,N.
\end{equation}

Provided that the \emph{interpolation matrix} $\Phi_{X,\ve}=(\Phi_{\ve}(\x_k,\x_j))_{1\leq k,j\leq N}$ is non-singular, we have the solution
\begin{equation}\label{eq:alphas}
\left(
\begin{array}{c}
\alpha_1 \\
\vdots \\
\alpha_N \\
\end{array}
\right)=(\Phi_{X,\ve})^{-1}u|_{X},
\end{equation} 
where $u|_{X}=\left(
\begin{array}{cccc}
u(\x_1) & u(\x_2) & \cdots &  u(\x_N) 
\end{array}
\right)^T$, hence the interpolant $\mathcal{P}_{X,\ve}u$ in \eqref{eq:pro_app} is known.\\

When is necessary to compute solutions of  \eqref{eq:boundary_pro} on large set of nodes, the resultant matrix for collocation method is dense, huge and ill-conditioned, carrying a prohibited computational cost. A variant of collocation method that allows to deal with large domains is the local version \cite{Tolstykh2003,fornb3}.  Making a local interpolation is possible to obtain a sparse matrix which discretizes the linear partial differential operator.  We now give a description of the RBF-FD method.

For an open and bounded set $\Omega\subset\R^d$,  let $X=\{ \x_i  \}_{i=1}^{N}\subset\om\cup\partial\Omega$ be a set of interpolation
nodes. For every $\x_i \in  X$, we create an influence domain $S_i\subset X$, which is formed by the $n_i$ nearest neighbor interpolation nodes, where $n_i$ is a positive integer, for $i=1,2,\ldots,N$. That is, we consider an $n_i$-stencil $S_i =\{\x_j^{i}  \}_{j=1}^{n_i}\subset X$, where $\x_1^i\equiv\x_i$ and  we denote the convex hull of the stencil $S_i$ by $\Omega_i$, i.e.,  $\Omega_i=\mbox{ConvexHull}(S_i)$. Thus we have a collection of subsets $\{S_i\}_{i=1}^{N}$ formed by nodes of $X$. By RBF interpolation, for any $\x\in\Omega$ we choose an $\Omega_i$ such that $\x\in \Omega_i$, so we can approximate $u(\x)$ as
\begin{equation}\label{eq:interpolant}
u(\x)\approx\widetilde{u}(\x)= \Pcal_{S_i,\ve_i}u(\x)=\sum_{j=1}^{n_i}\al_j^i\Phi_{\ve_i}(\x,\x_{j}^{i}).
\end{equation}
Note we have taken the shape parameter $\ve_i$ depending of the location $\x_i$, thus the shape parameter of the RBF is conveniently manipulated, according to local known information related with the PDE. 

By collocating the $n_i$ nodes of the stencil $S_i$, we obtain a small linear system
\begin{equation}\label{eq7}
\Phi_{S_i,\ve_i} \bm{\alpha}^i={ U}_{i}
\end{equation}
where ${U}_{i}=\left(\begin{array}{cccc}
\widetilde{u}(\x_1^i) & \widetilde{u}(\x_2^i)& \cdots & \widetilde{u}(\x_{n_i}^i)
\end{array}\right)^T$, $\Phi_{S_i,\ve_i}=\left(\Phi_{\ve_i}(\x_j^i,\x_k^i) \right)_{1\le j,k\le n_i}$ is the \emph{local interpolation matrix} and  $\bm{\alpha}^i =\left( \al_1^i\ \ \al_2^i\ \ \cdots\ \ \al_{n_i}^i\right)^T$.

The unknown coefficients $\bm{\alpha}^i$ in (\ref{eq7}) can be expressed in terms of the function values at the local interpolation nodes as
\begin{equation}\label{eq8}
\bm{\alpha}^i=\Phi_{S_i,\ve_i}^{-1}U_i.
\end{equation}
The inverse matrix $\Phi_{S_i,\ve_i}^{-1}$ exists provided that $\Phi_{S_i,\ve_i}$ is positive definite, which is true for Bessel RBF (\ref{eq:bessel_RBF}).

Now, with the aim to get a local discretized version for (\ref{eq:boundary_pro}) we consider $\x_i\in\Omega\cap X$ (or $\x_i\in\partial\Omega\cap X$). In both cases it must be applied a linear partial differential operator, either $\ml$ or $\mb$, to equation (\ref{eq:interpolant}). For $\x_i\in\Omega$, with \eqref{eq8}, we have
\begin{align}
\ml\Pcal_{S_i,\ve_i}u(\x_i)&=\sum_{j=1}^{n_i}\al_j^i\ml\Phi_{\ve_i}(\x,\x_j^i)|_{\x=\x_i}\nonumber\\
&=\ml\Phi_{S_i,\ve_i}^1\ \bm{\alpha}^i\nonumber\\
&=\ml\Phi_{S_i,\ve_i}^1\Phi_{S_i,\ve_i}^{-1}U_i\label{eq:weights_rbf}
\end{align}
where $$\ml\Phi_{S_i,\ve_i}^1=\left(\begin{array}{cccc}
\ml\Phi_{\ve_i,}(\x,\x_1^i)|_{\x=\x_i}&\ml\Phi_{\ve_i,}(\x,\x_2^i)|_{\x=\x_i}&\cdots&\ml\Phi_{\ve_i,}(\x,\x_{n_i}^i)|_{\x=\x_i}
\end{array}\right)$$
is a row matrix. Similarly, for  $\x_i\in\partial\om\cap X$  we have
\begin{equation}\label{eq:weights_rbf_bound}
\mb\Pcal_{S_i,\ve_i}u(\x_i)=\mb\Phi_{S_i,\ve_i}^1\Phi_{_{S_i,\ve_i}}^{-1}U_i.
\end{equation}
We denote $W\ml_{S_i,\ve_i}=\ml\Phi_{S_i,\ve_i}^1\Phi_{S_i,\ve_i}^{-1}$ and $W\mb_{S_i,\ve_i}=\mb\Phi_{S_i,\ve_i}^1\Phi_{S_i,\ve_i}^{-1}$. From  (\ref{eq:weights_rbf}) and (\ref{eq:weights_rbf_bound}) it follows a discretized local version of (\ref{eq:boundary_pro}) 
\begin{align*}
W\ml_{S_i,\ve_i}U_i&=f(\x_i),\mbox{ if } \x_i\in\Omega\cap X\\
W\mb_{S_i,\ve_i}U_i&=g(\x_i), \mbox{ if } \x_i\in\partial\Omega\cap X,
\end{align*}
for $1\leq i \leq N$.

The above system of linear equations can be assembled forming a sparse matrix $\mathbf{H}$ of size $N\times N$ where the $i-$th row, associated to $\x_i\in X$, has at most  $n_i$ nonzero entries, and the unknown column matrix is given by $U=\left(\begin{array}{cccc}
\widetilde{u}(\x_1)&\widetilde{u}(\x_2)&\cdots&\widetilde{u}(\x_N) 
\end{array}\right)^T$ (we recall that $\x_i\equiv\x_1^i$). We can then obtain an
approximated solution $\td{u}(\x)$ at all interpolation nodes by solving $\mathbf{H}U=F$. The $\mathbf{H}$ can be thought as a discretized version of the problem \eqref{eq:boundary_pro}.

 \subsection{Bessel RBF-FD}
Suppose that $u$ is a solution of the Helmholtz equation 
$\Delta u(\x)+k(\x)^2u(\x)=0$, for $\x\in\Omega$, and $u(\x)$, with $\x\in\partial\Omega$, it satisfies certain boundary condition. If $X=\{\x_i\}_{i=1}^N\subset\Omega\cup\partial\Omega$ is a set of nodes, for $\x_i\in X\cap\Omega$ we take a stencil $S_i=\{\x_j^i\}_{j=1}^{n_i}\subset X$ based on $\x_i$, with $\x_1^i=\x_i$. For $\x\in\Omega_i=\mbox{ConvexHull}(S_i)$ we define, with $k_i=k(\x_i)$, the interpolant 
\begin{equation}\label{eq:interpolant_bessel}
\widetilde{u}(\x)=\sum_{j=1}^{n_i}\alpha_j^i J_0(k_i\|\x-\x_j^i\|).
\end{equation}
With the local interpolation matrix, $\mathbf{J}_{k_i}=(J_0(k_i\|\x_l^i-\x_j^i\|))_{1\leq l,j\leq n_i}$, which is positive definite \cite{FORNBERG2006_oscillatory_RBF}, and forcing the condition $\widetilde {u}(\x_l^i)=u(\x_l^i)$, then from \eqref{eq:interpolant_bessel} we have the linear equation
 \begin{equation}\label{eq:system_interpolation_ill}
 U_i=\mathbf{J}_{k_i}\bm{\alpha}_i,
 \end{equation}
where $U_i=\left(\begin{array}{cccc}
u(\x_1^i)&u(\x_2^i)&\cdots&u(\x_{n_i}^i)
\end{array}\right)^T$ and $\bm{\alpha}_i=\left(\begin{array}{cccc}
\alpha_1^i&\alpha_2^i&\cdots&\alpha_{n_i}^i
\end{array}\right)^T$.
In view that $\phi_k$, defined in \eqref{eq:bessel_RBF2D}, satisfies the homogeneous Helmholtz equation, then
\begin{eqnarray*}
\Delta\widetilde{u}(\x)|_{\x=\x_i}&=&\sum_{j=1}^{n_i}\alpha_j^i \Delta J_0(k_i\|\x-\x_j^i\|)|_{\x=\x_i}\\
&=&-k_i^2\sum_{j=1}^n\alpha_j^i J_0(k_i\|\x_i-\x_j^i\|)\\
&=&-k_i^2\widetilde{u}(\x_i).
\end{eqnarray*}
Hence the interpolant \eqref{eq:interpolant_bessel} satisfies the homogeneous Helmholtz equation. On the other hand, applying  $\Delta_{S_i,k_i}$  to the solution $u$, we have
\begin{eqnarray*}
\Delta_{S_i,k_i}u(\x_i)&=&\Delta \mathbf{J}_{S_i,k_i}^1\mathbf{J}_{k_i}^{-1}U_i\\
&=&-k^2_i\mathbf{e}_1U_i, (\mbox{ where } \mathbf{e}_1=(1\ \ 0\ \ 0\ \ \cdots\ \ 0))\\
&=&-k_i^2u(\x_i).
\end{eqnarray*}
Note that $\Delta_{S_i,k_i}u(\x_i)-\Delta\widetilde{u}(\x_i)=-k_i^2(u(\x_i)-\widetilde{u}(\x_i))$, thus, for solutions of the homogeneous Helmholtz problem the local truncation error for the Laplace operator has a theoretical error depending of wavenumber at $\x_i$ and of the error of the local interpolant. The error of the approximated  solutions is produced by the interpolant \eqref{eq:interpolant_bessel} and by the ill-conditioning of the matrix $\JJ_{k_i}$, in solving the linear system \eqref{eq:system_interpolation_ill}. Next we will deal with solutions of these systems.

\section{Method of Diagonal Increments (MDI)}\label{sec:MDI}
The interpolation matrix $\JJ_k$ is ill-conditioned, especially for certain node distributions. In literature there are several methods for dealing with the ill-conditioning when the shape parameter is small \cite{stable_general_bessel_rbf}, but in our case we are taking the shape parameter as the wavenumber $k$, which can be large. So we have chosen the MDI. For our case $\JJ_k$ will be considered  ill-conditioned when the condition number in the spectral norm \footnote{The spectral norm matches with the matrix norm induced by Euclidean norm for vectors, i.e., \\
$\| \mathbf{A} \|_2=\sup\{\|\mathbf{A}\bm{\alpha}\|_2 \ :\ \bm{\alpha}\in \R^d \mbox{ with } \|\bm{\alpha}\|_2=1\}$. }, $\kappa(\JJ_k)=\|\JJ_k\|_2\|\JJ_k^{-1}\|_2$, satisfies  $\kappa(\JJ_k)>10^{15}$, which hinders that the solution $\bm{\alpha}$ of $U=\mathbf{J}_k\bm{\alpha}$ be accurately calculated, in double precision, through Cholesky factorization. An alternative to compute $\bm{\alpha}$ with better tolerance respect to  large condition numbers, allowing roughly 2 orders of magnitude more,  i.e.,  up to $\kappa(\JJ_k)\sim 10^{17}$, is the Block-$LDL^T$-decomposition (LDL$^T$)  \cite{LDL_factorization,LDL_factorization_analysis}. When $\JJ_k$ is ill-conditioned we solve the better conditioned problem $U=(\mathbf{J}_k+\beta\mathbf{I})\widetilde{\bm{\alpha}}$ instead, where $\II$ is  the identity matrix and $\beta$  a small positive real number.  Next, we will give some important aspects about the spectrum of $\JJ_k+\beta\II$. The following development is based on the Riley's method \cite{riley_method}.

\begin{rmrk}\label{rmrk:convergence_neumann_series_jk}

	$\JJ_k$ is positive definite \cite{FORNBERG2006_oscillatory_RBF}, thus its spectrum is real and positive. If $\ds\{\lambda_m\}_{m=1}^n$ is its spectrum,  with $\lambda_1\geq\lambda_2\geq\cdots\geq\lambda_n$, then $\{\lambda_m+\beta\}_{m=1}^n$ is the spectrum of $\widetilde{\JJ}_k=\JJ_k+\beta\II$ and $\{\frac{\beta}{\lambda_m+\beta}\}_{m=1}^n$ is the spectrum of $\beta\widetilde{\JJ}_k^{-1}$, hence we have the spectral norms,  $\|\beta\widetilde{\JJ}_k^{-1}\|_2=\frac{\beta}{\lambda_n+\beta}$ and $\|	(\II-\beta\widetilde{\JJ}_k^{-1})^{-1}\|_2=\frac{\lambda_n+\beta}{\lambda_n}$. The  above implies that the Neumann series  $\sum_{m=0}^{\infty}(\beta\widetilde{\JJ}_k^{-1})^m$ converges and the equality 
	\begin{equation}\label{eq:neumann_series}
	(\II-\beta\widetilde{\JJ}_k^{-1})^{-1}=\sum_{m=0}^{\infty}(\beta\widetilde{\JJ}_k^{-1})^m
	\end{equation} holds.

\end{rmrk}

\begin{rmrk}\label{rmrk:better_coditioning_wideJk}

	If $\{\lambda_m\}_{m=1}^n$ is the spectrum  of $\JJ_k$, as in remark \ref{rmrk:convergence_neumann_series_jk},  then the condition number of $\JJ_k$ is given by $\kappa(\JJ_k)=\frac{\lambda_1}{\lambda_n}$ and $\kappa(\widetilde{\JJ}_k)=\frac{\lambda_1+\beta}{\lambda_n+\beta}$, which implies that 
	\begin{equation*}
	\kappa(\widetilde{\JJ}_k)<\kappa({\JJ}_k).
	\end{equation*}
	With this, the matrix $\widetilde{\JJ}_k$ is better conditioned than $\JJ_k$.

\end{rmrk}

\begin{rmrk}\label{rmrk:aprrox_jk_by_neumann_series}

Given that $\widetilde{\JJ}_k=\JJ_k+\beta\II$, then $\JJ_k^{-1}=\widetilde{\JJ}_k^{-1}(\II-\beta\widetilde{\JJ}_k^{-1})^{-1}$.
If $\bm{\alpha}$ is the true solution of the equation $\JJ_k\bm{\alpha}=U$ and $\widetilde{\bm{\alpha}}$ is the solution for the perturbed system $\widetilde{\JJ}_k\widetilde{\bm{\alpha}}=U$, we can compare $\bm{\alpha}$ and $\widetilde{\bm{\alpha}}$. Note the following:
\begin{equation}
\begin{array}{ccc}
{\JJ}_k^{-1}-\widetilde{\JJ}_k^{-1}&=&\widetilde{\JJ}_k^{-1}(\II-\beta\widetilde{\JJ}_k^{-1})^{-1}-\widetilde{\JJ}_k^{-1}\\
&=&\widetilde{\JJ}_k^{-1}\left((\II-\beta\widetilde{\JJ}_k^{-1})^{-1}-\II\right)
\end{array}
\end{equation}
thus, 
\begin{equation}
\bm{\alpha}-\widetilde{\bm{\alpha}}=\frac{1}{\beta}(\beta\widetilde{\JJ}_k^{-1})\left((\II-\beta\widetilde{\JJ}_k^{-1})^{-1}-\II\right)U.
\end{equation}
Since $\|(\II-\beta\widetilde{\JJ}_k^{-1})^{-1}-\II\|_2=\frac{\beta}{\lambda_n}$, finally we have, 
\begin{equation}\label{eq:bound_error_alpha_bar}
\|\bm{\alpha}-\widetilde{\bm{\alpha}}\|_2\leq \frac{1}{\lambda_n}\left(\frac{\beta}{\lambda_n+\beta}\right)\|U\|_2.
\end{equation}
\end{rmrk}

Now, with the purpose of obtaining closer solutions to the true one $\bm{\alpha}$, and improve the error bound \eqref{eq:bound_error_alpha_bar}, we consider the following. By using the Neumann series \eqref{eq:neumann_series} we have
	 \begin{eqnarray*}
	 	\JJ_k^{-1}&=&\widetilde{\JJ}_k^{-1}\sum_{m=0}^{\infty}(\beta\widetilde{\JJ}_k^{-1})^m\\
	 	&=&\frac{1}{\beta}\sum_{m=1}^{\infty}(\beta\widetilde{\JJ}_k^{-1})^m.
	 \end{eqnarray*}
	If  $\JJ_k\bm{\alpha}=U$ and  $\widetilde{\JJ}_k\widetilde{\bm{\alpha}}=U$ then, from the remark \ref{rmrk:aprrox_jk_by_neumann_series}, 
	\begin{equation}\label{eq:alpha_series}
	\bm{\alpha}=\frac{1}{\beta}\sum_{m=1}^{\infty}(\beta\widetilde{\JJ}_k^{-1})^mU
	\end{equation} and
	\begin{equation}\label{eq:alpjaM_sum}
	\bm{\alpha}=\sum_{m=0}^{\infty}(\beta\widetilde{\JJ}_k^{-1})^m\widetilde{\bm{\alpha}}.
	\end{equation}
If we truncate the series in \eqref{eq:alpha_series} up to order $M$, we obtain an approximation of the true   solution $\bm{\alpha}$, we denote it by  
\begin{equation}\label{eq:approx_alpha_barM}
\widetilde{\bm{\alpha}}_M=\frac{1}{\beta}\sum_{m=1}^{M}(\beta\widetilde{\JJ}_k^{-1})^mU.
\end{equation}
From \eqref{eq:alpha_series} and \eqref{eq:approx_alpha_barM}, the error of the approximation $\widetilde{\bm{\alpha}}_M$ can be bounded by using the formula
\begin{eqnarray*}
\bm{\alpha}-\widetilde{\bm{\alpha}}_M&=&\frac{1}{\beta}\sum_{m=M+1}^{\infty}(\beta\widetilde{\JJ}_k^{-1})^mU\\
&=&\frac{1}{\beta}(\II-\beta\widetilde{\JJ}_k^{-1})^{-1}(\beta\widetilde{\JJ}_k^{-1})^{M+1}U.
\end{eqnarray*}
Finally, by taking the Euclidean norm, we have in terms of the spectral norm,
\begin{equation*}
\|\bm{\alpha}-\widetilde{\bm{\alpha}}_M\|_2\leq\frac{1}{\beta}\| (\beta\widetilde{\JJ}_k^{-1})^{M+1}\|_2\|(\II-\beta\widetilde{\JJ}_k^{-1})^{-1} \|_2\|U\|_2.
\end{equation*}
From the remark \ref{rmrk:convergence_neumann_series_jk} we can conclude that
\begin{equation}\label{eq:error_bound_MDI}
\|\bm{\alpha}-\widetilde{\bm{\alpha}}_M\|_2\leq\frac{1}{\lambda_n}\left(\frac{\beta}{\lambda_n+\beta}\right)^{M}\|U\|_2.
\end{equation}

An iterative procedure to compute \eqref{eq:approx_alpha_barM}, with the better conditioned matrix $\beta\widetilde{\JJ}_k$, can be obtained just by noting that, with $\widetilde{\bm{\alpha}}=\widetilde{\JJ}_k^{-1}U$,
  \begin{eqnarray*}\label{eq:approx_alpha_barM_iter}
  \widetilde{\bm{\alpha}}_M&=&\frac{1}{\beta}\sum_{m=1}^{M}(\beta\widetilde{\JJ}_k^{-1})^mU\\
  &=&\sum_{m=1}^{M}(\beta\widetilde{\JJ}_k^{-1})^{m-1}\widetilde{\bm{\alpha}}\\
  &=&\widetilde{\bm{\alpha}}+\beta\widetilde{\JJ}_k^{-1}\left(\widetilde{\bm{\alpha}}+\beta\widetilde{\JJ}_k^{-1}\left(\widetilde{\bm{\alpha}}+\cdots\right)\right).
  \end{eqnarray*}
  Hence we can compute $\widetilde{\bm{\alpha}}_M$  as: 
  \begin{equation}
  	\begin{aligned}
  	\widetilde{\bm{\alpha}}_0=&\ \widetilde{\JJ}_k^{-1}U\\
  	\widetilde{\bm{\alpha}}_{m}=&\ \widetilde{\bm{\alpha}}_0+\beta\widetilde{\JJ}^{-1}_k\widetilde{\bm{\alpha}}_{m-1}, \ \ \ \mbox{\ \  for } m=1,2,\ldots,M.
  	\end{aligned}\label{eq:ITMDI}
  \end{equation}
We call this algorithm the Iterative MDI (ITMDI) and its origin goes back to \cite{riley_method}. 

Since  $0<\frac{\beta}{\lambda_n+\beta}<1$, it is important to note that the parameter $\beta$ needs to be selected to be large enough to improve conditioning but small enough so that the convergence of the method is faster \cite{SARRA2014_regularized_positive_definite, riley_method}. An undesirable situation is when $\lambda_n$ is near to the machine epsilon,\footnote{In double precision the machine epsilon is approximately 2.22e-16.} in theses cases the ratio $\frac{\beta}{\lambda_n+\beta}$ is very close to 1 and in the convergence may be too slow. However, we have seen that with a small $M$ (for example $M=15$) is enough to improve the error \eqref{eq:bound_error_alpha_bar}. All codes were typed in Matlab R2016a and run in a laptop \textcolor{black}{ with Core i7 processor at 2.8 Ghz with 12\,GiB  of RAM.}

\subsection{Local truncation error and condition number}
\begin{figure}[ht]
	\begin{center}	
		\begin{tabular}{|c|c|c|c|c|c|}
			\hline
			\hline
			\includegraphics[scale=0.25]{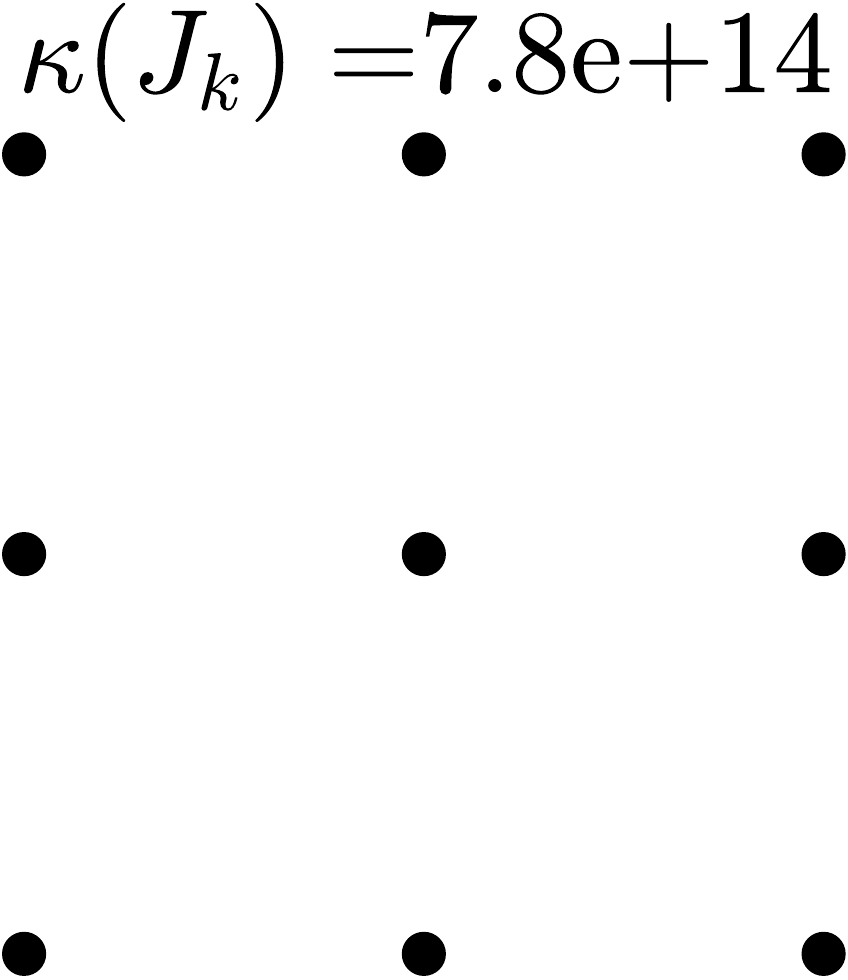}&\includegraphics[scale=0.25]{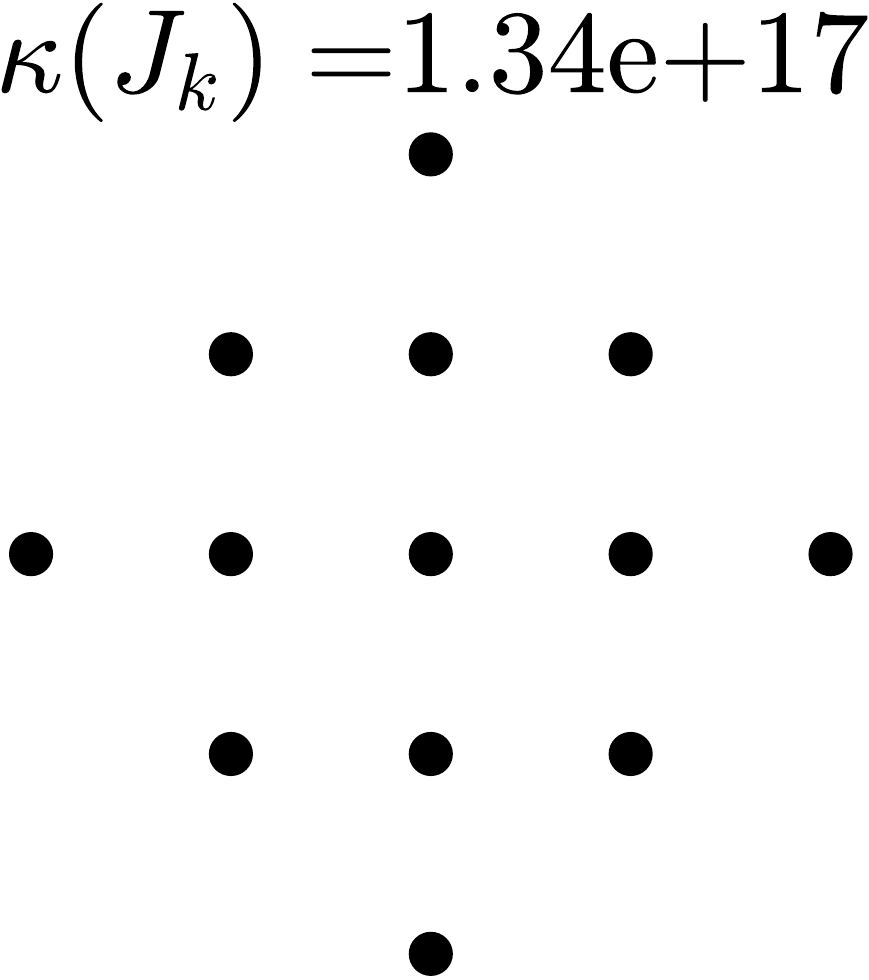} &\includegraphics[scale=0.25]{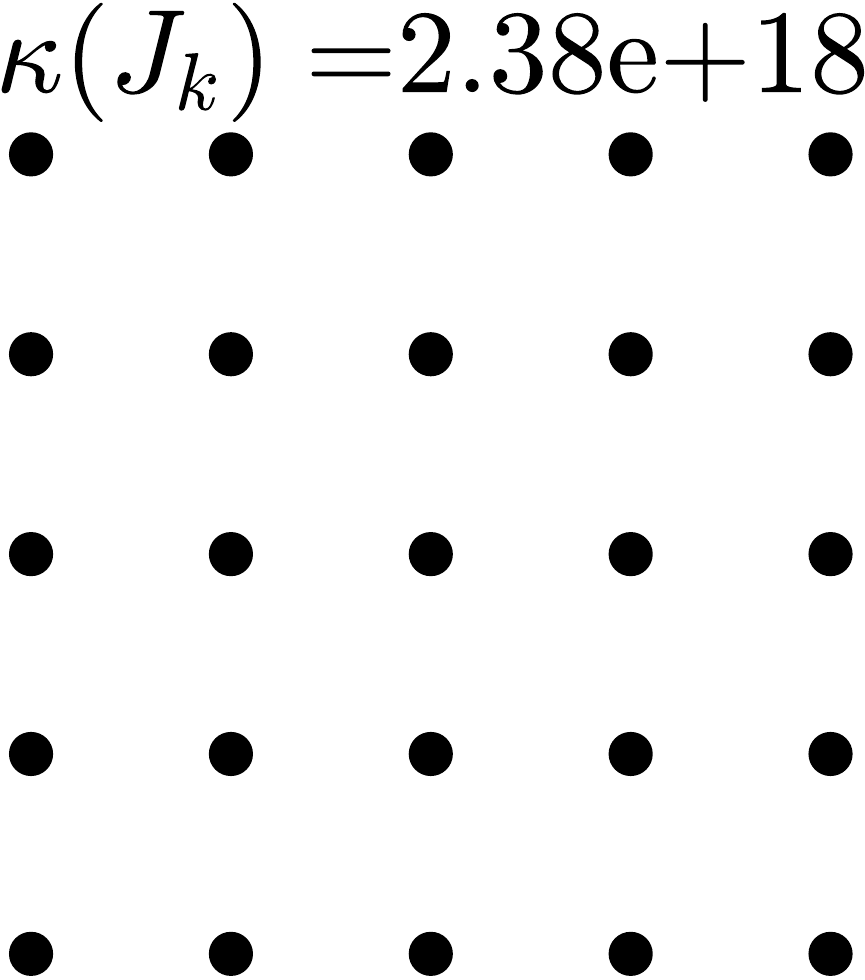} &\includegraphics[scale=0.25]{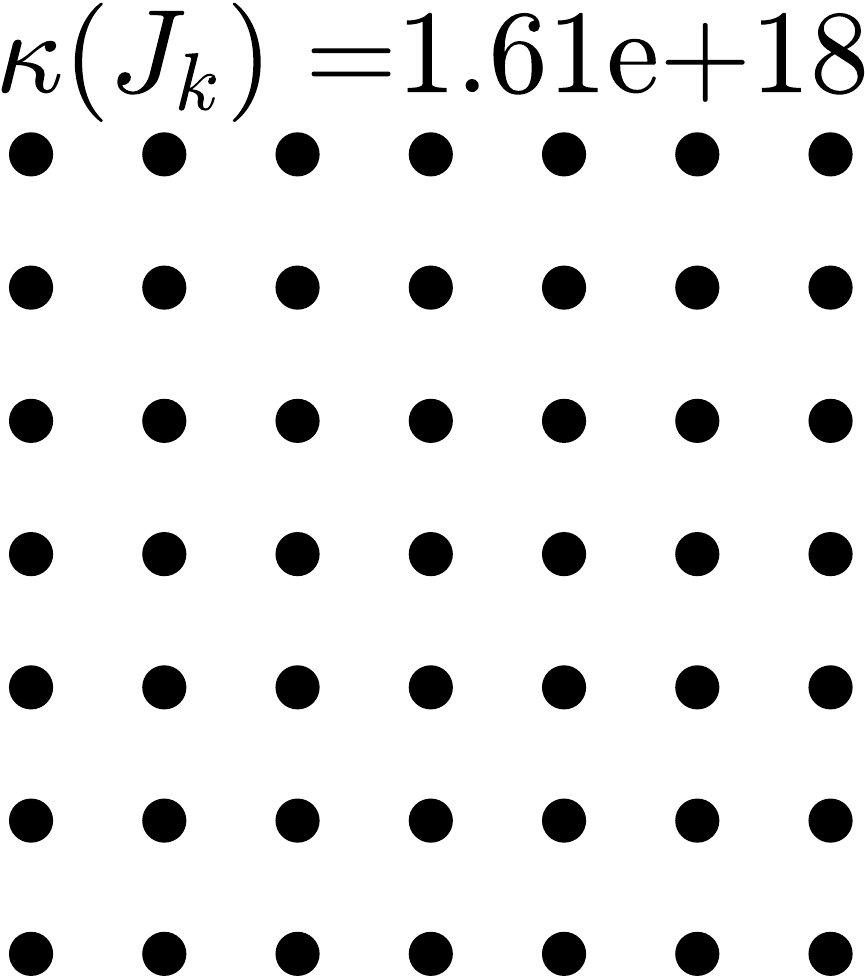} &\includegraphics[scale=0.25]{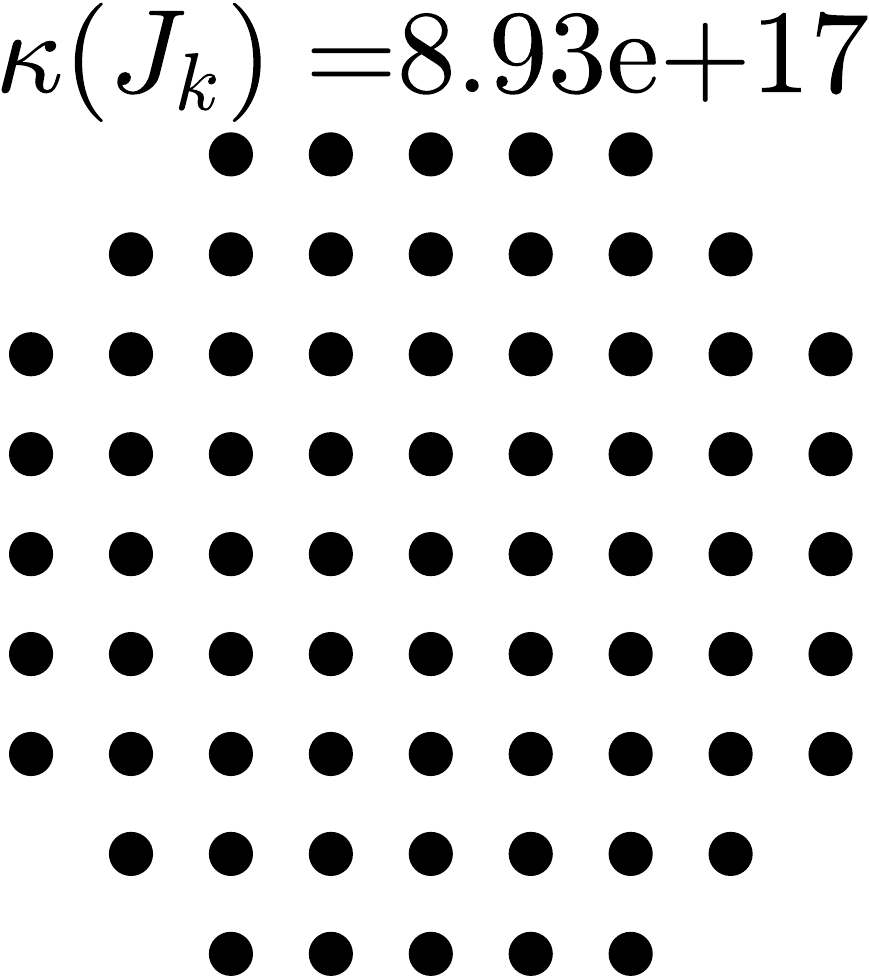} &\includegraphics[scale=0.25]{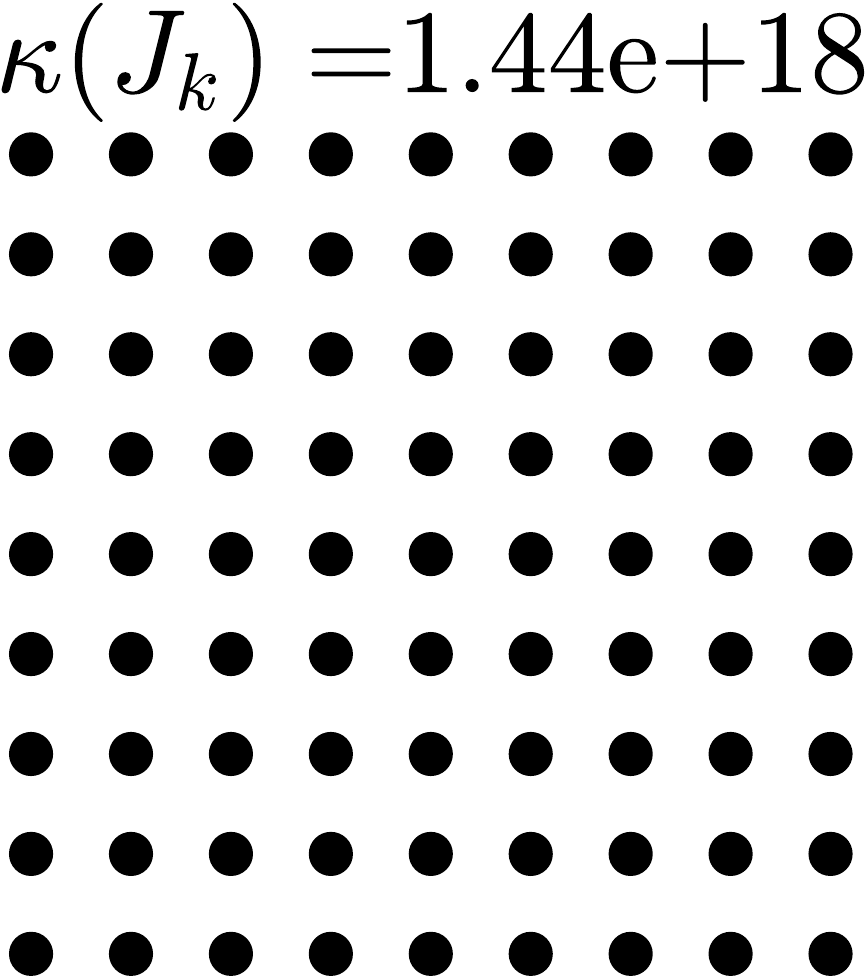}\\
			\hline
			\hline
			\includegraphics[scale=0.25]{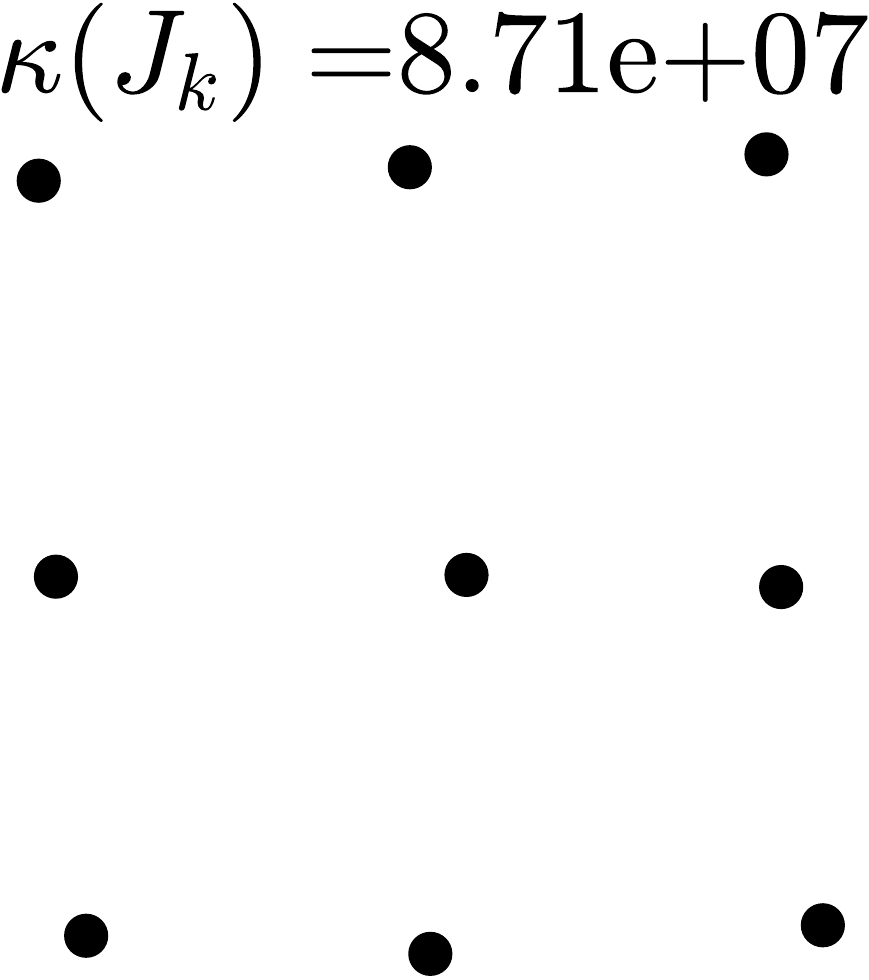}&\includegraphics[scale=0.25]{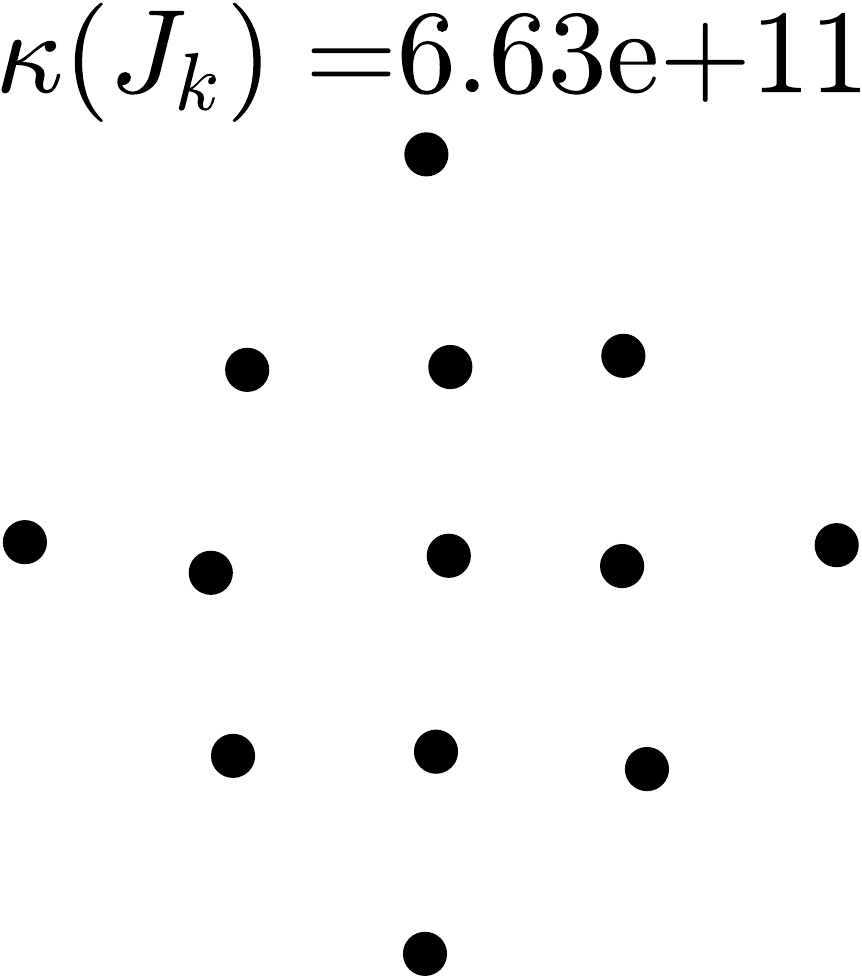} &\includegraphics[scale=0.25]{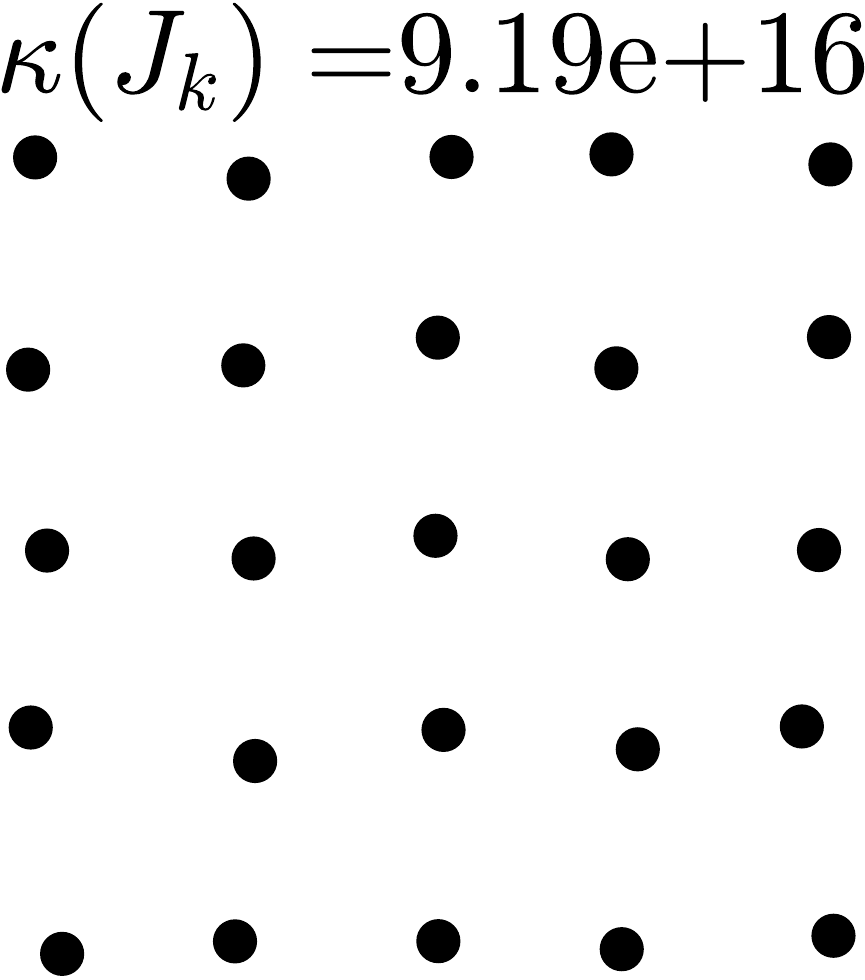} &\includegraphics[scale=0.25]{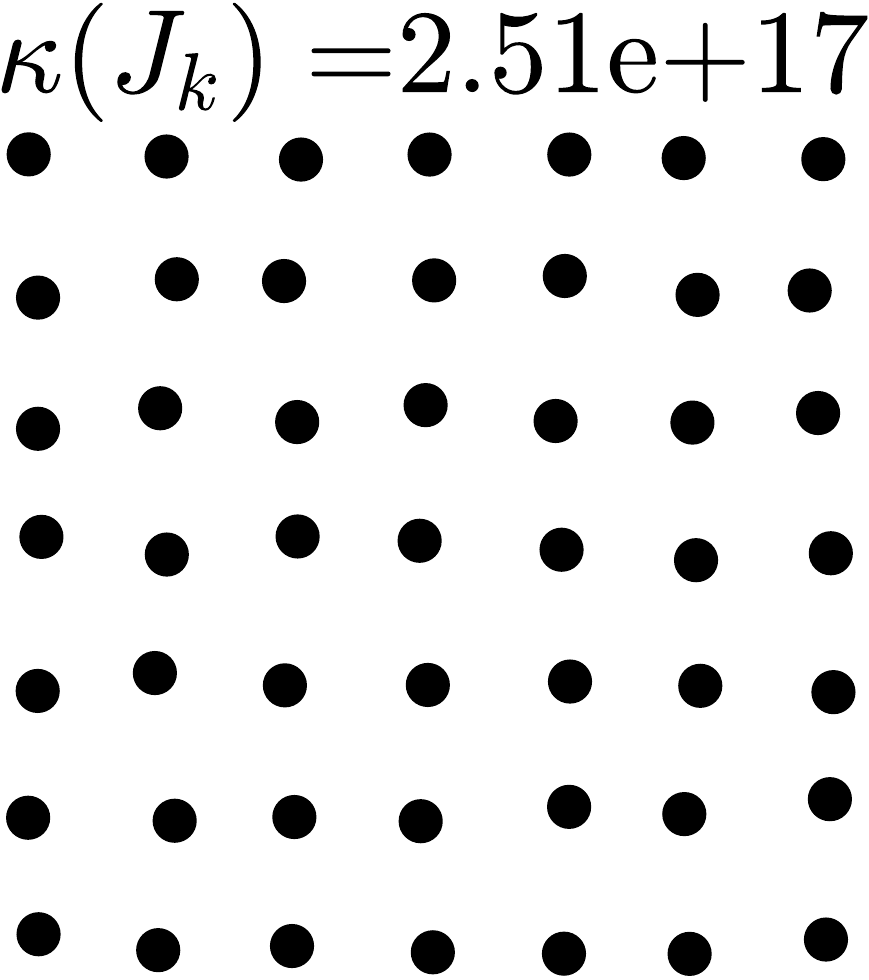}
			 &\includegraphics[scale=0.25]{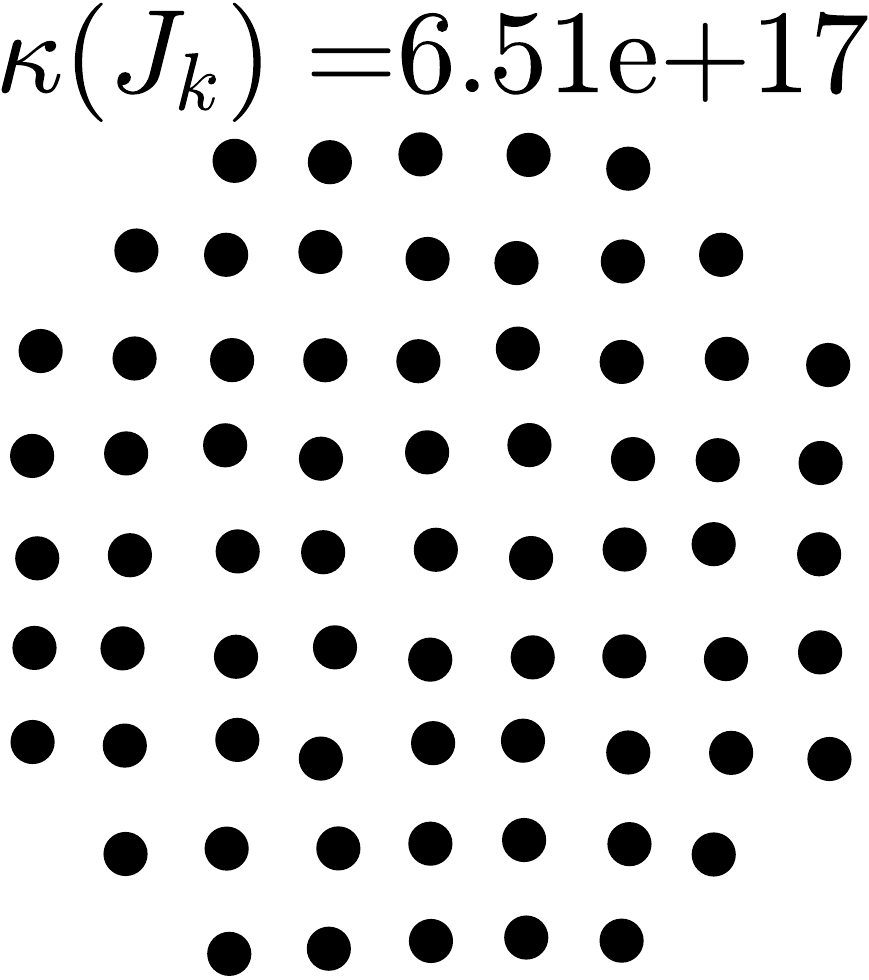} &\includegraphics[scale=0.25]{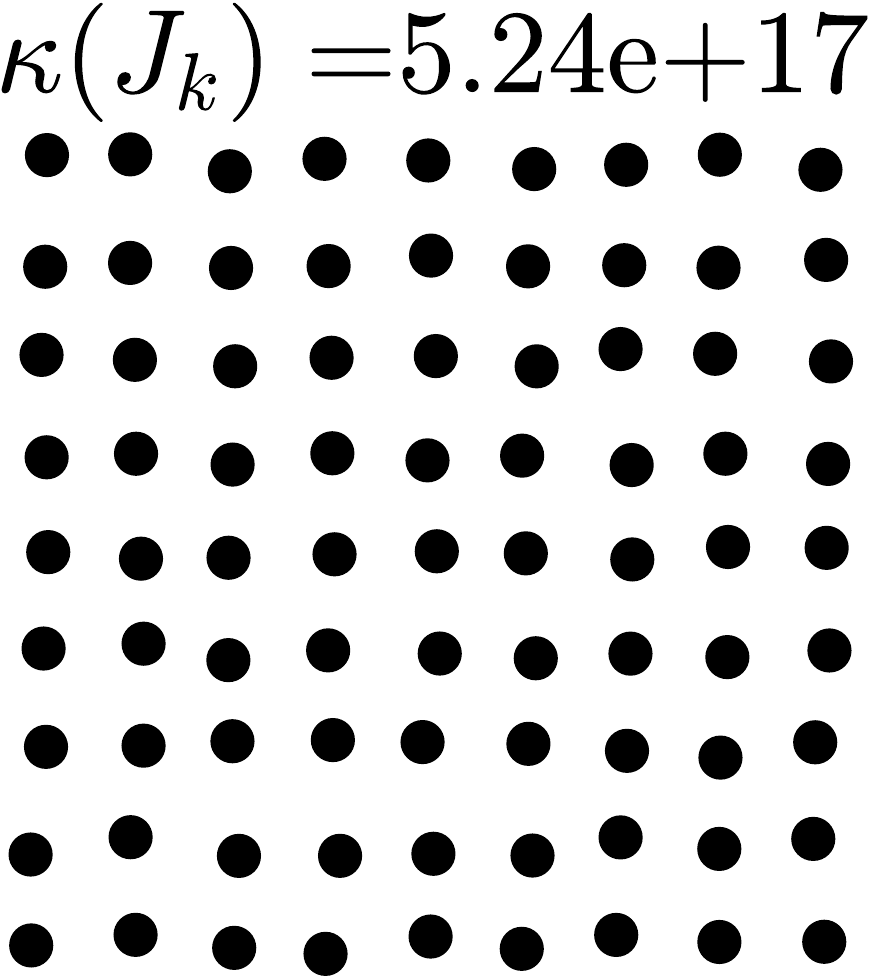}\\
			\hline
		\end{tabular}
	\end{center}
	\caption{Plots of some small stencils with the respective approximated condition number of $\JJ_k$, with $k=100$ and $h=\frac{2\pi}{8k}$. Top row: Stencils are taken from a regular square grid. Bottom row: perturbed position from stencils of top row.} 	\label{fig:stencils_square}
\end{figure}

\begin{figure}[ht]
	\begin{center}	
		\begin{tabular}{|c|c|c|c|c|c|}
			\hline
			\hline
			\includegraphics[scale=0.23]{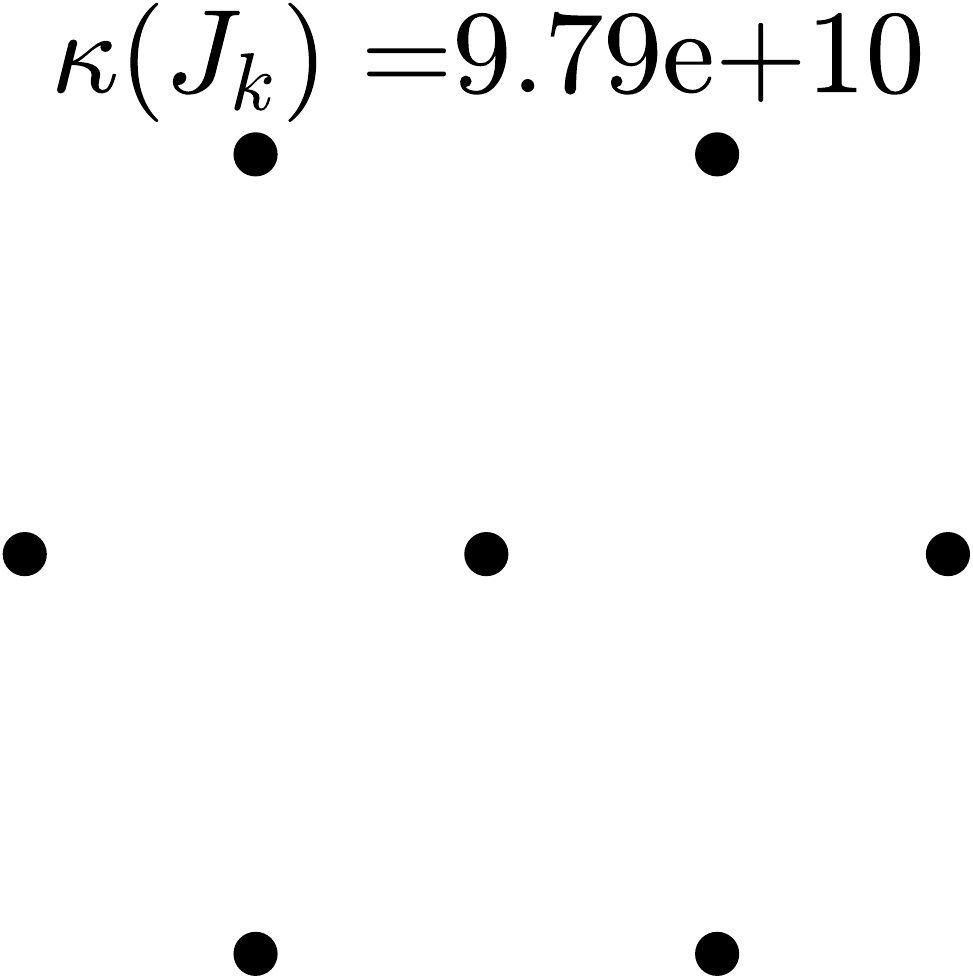}&\includegraphics[scale=0.23]{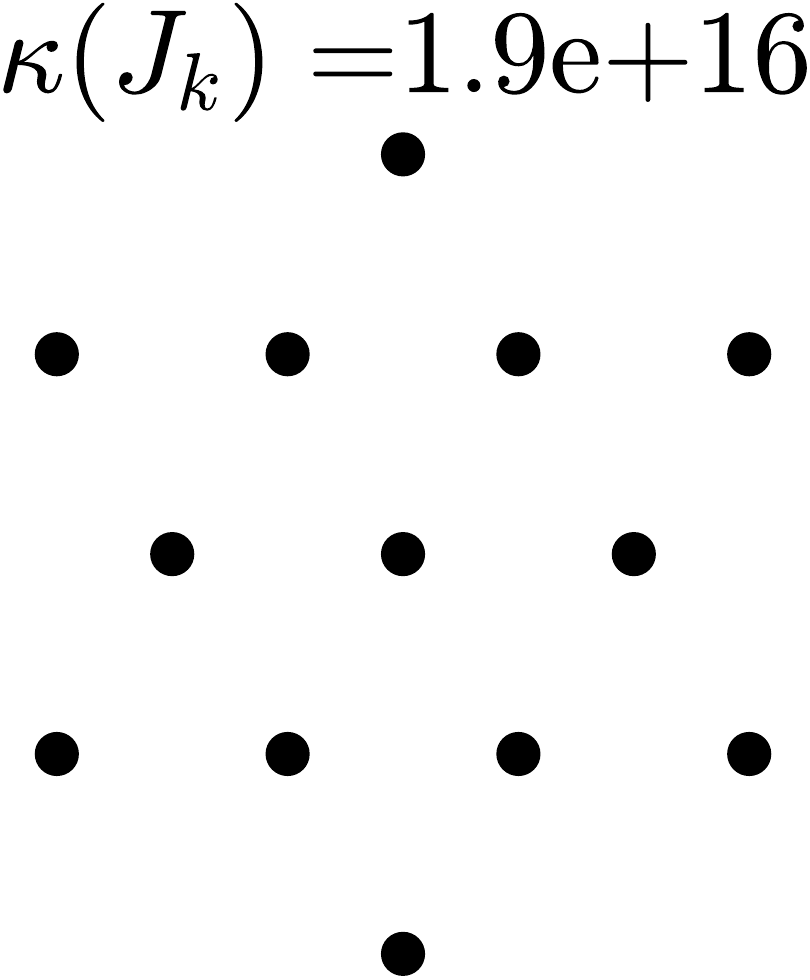} &\includegraphics[scale=0.23]{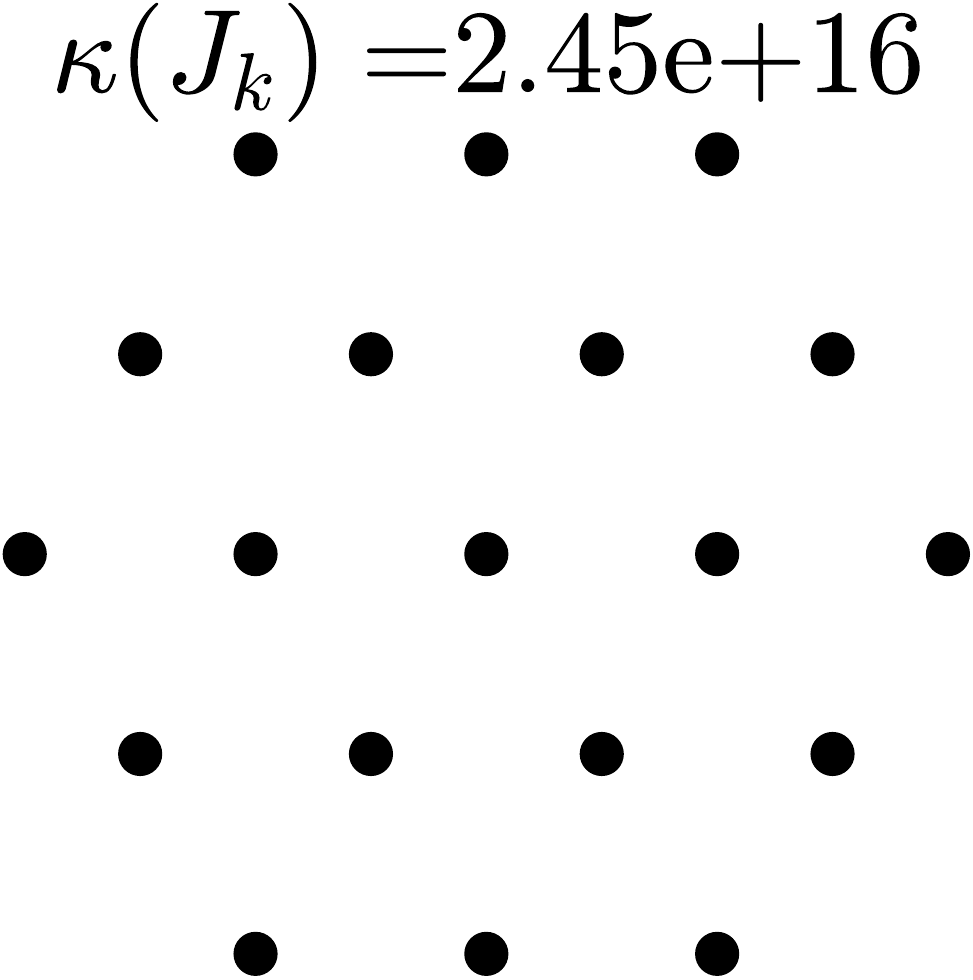} &\includegraphics[scale=0.23]{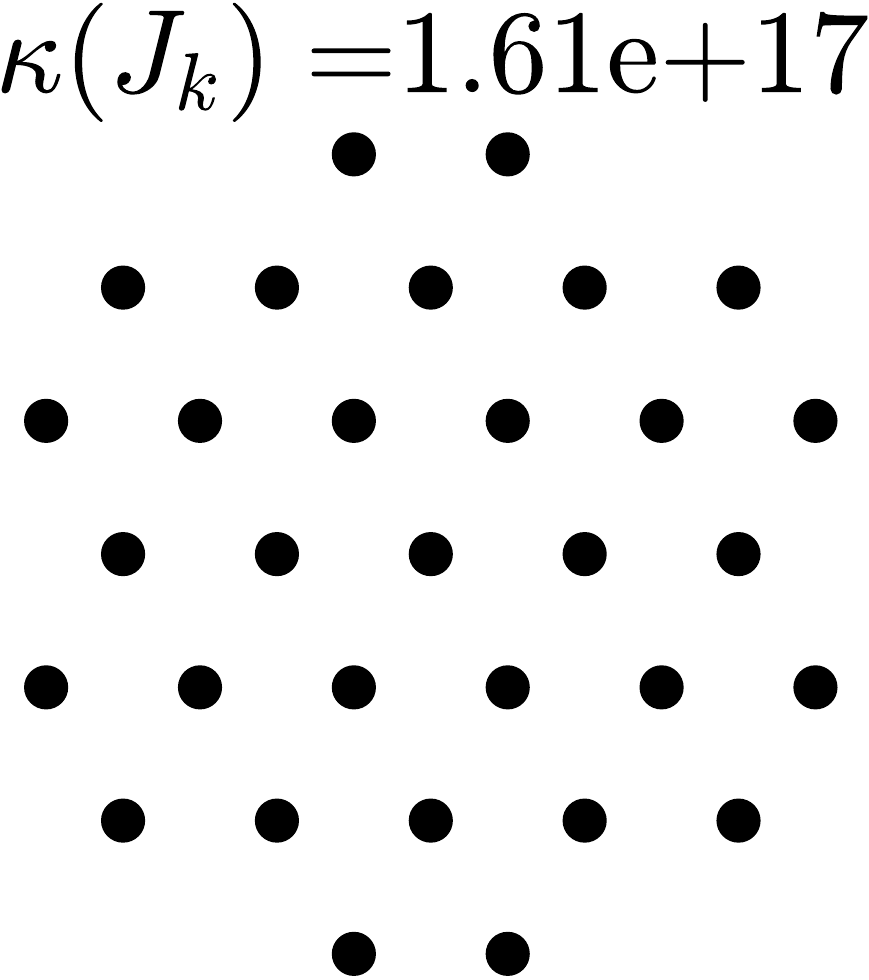}
			&\includegraphics[scale=0.23]{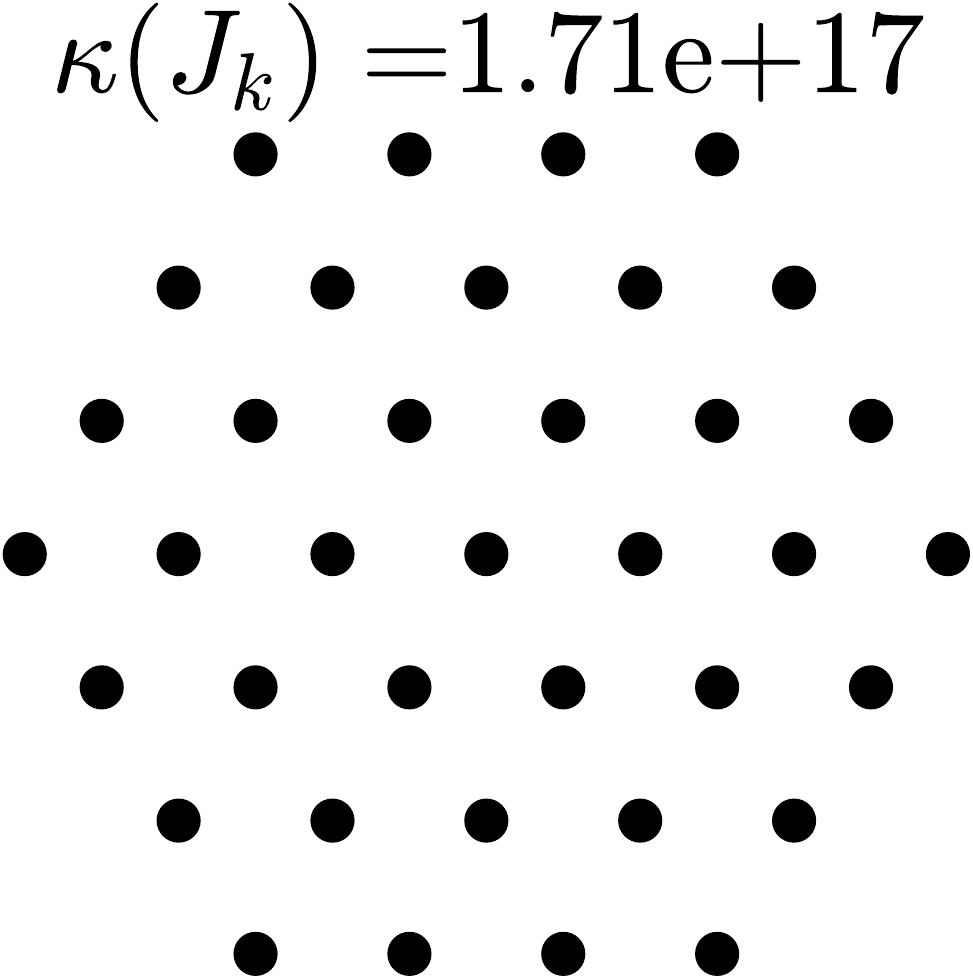}&\includegraphics[scale=0.23]{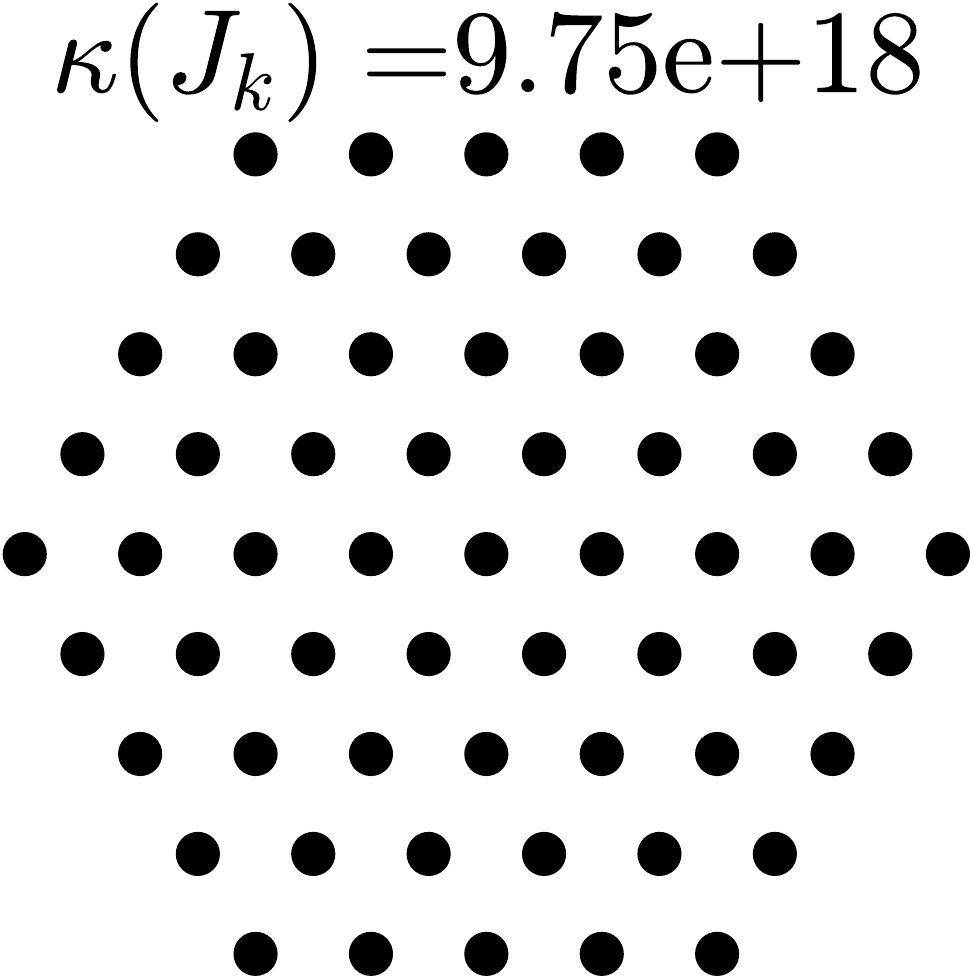}\\
			\hline
			\hline
			\includegraphics[scale=0.23]{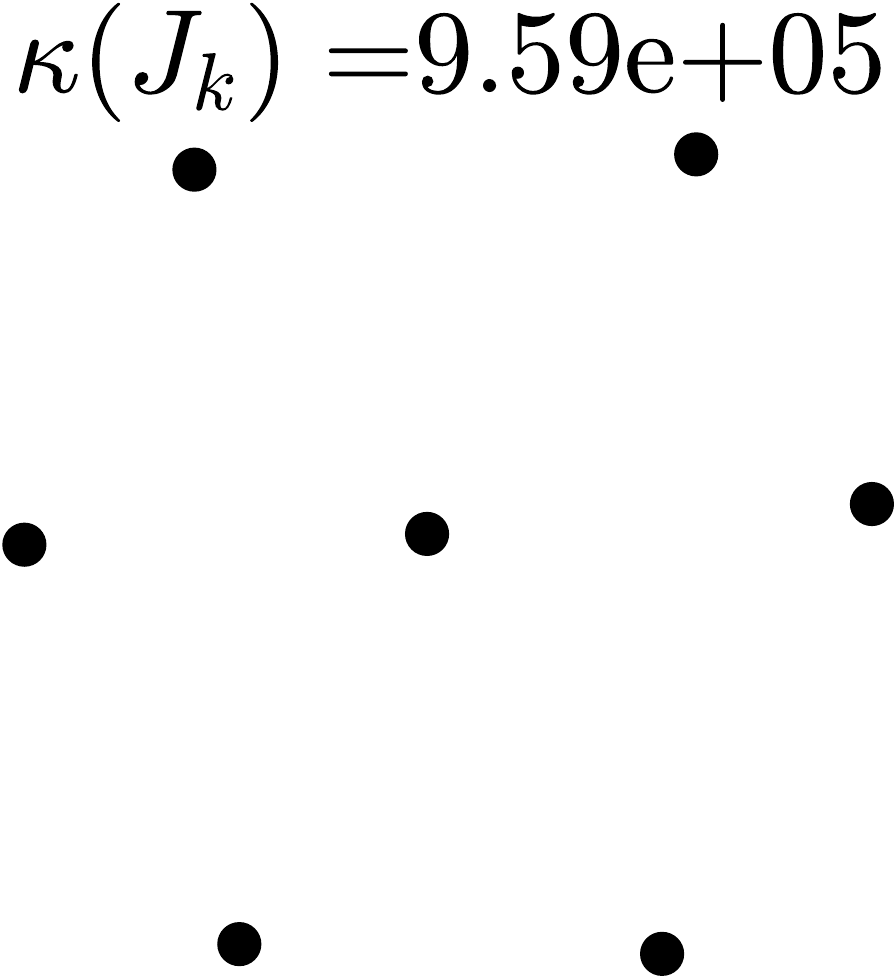}&\includegraphics[scale=0.23]{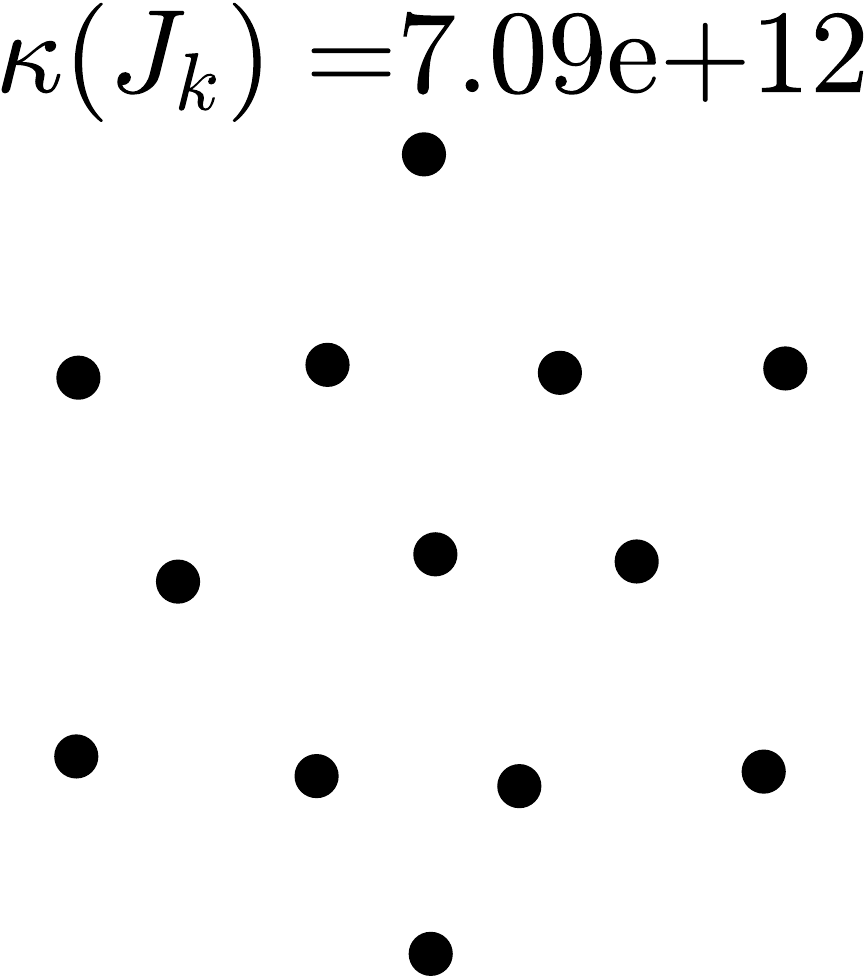} &\includegraphics[scale=0.23]{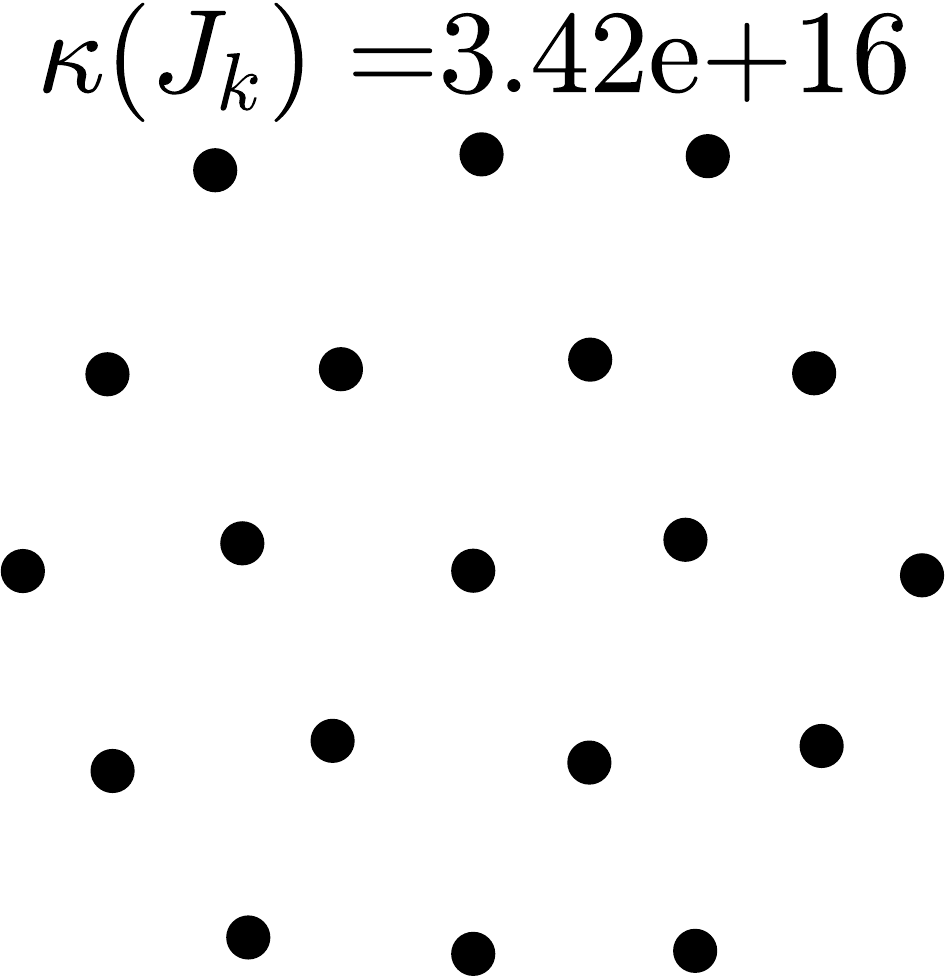} &\includegraphics[scale=0.23]{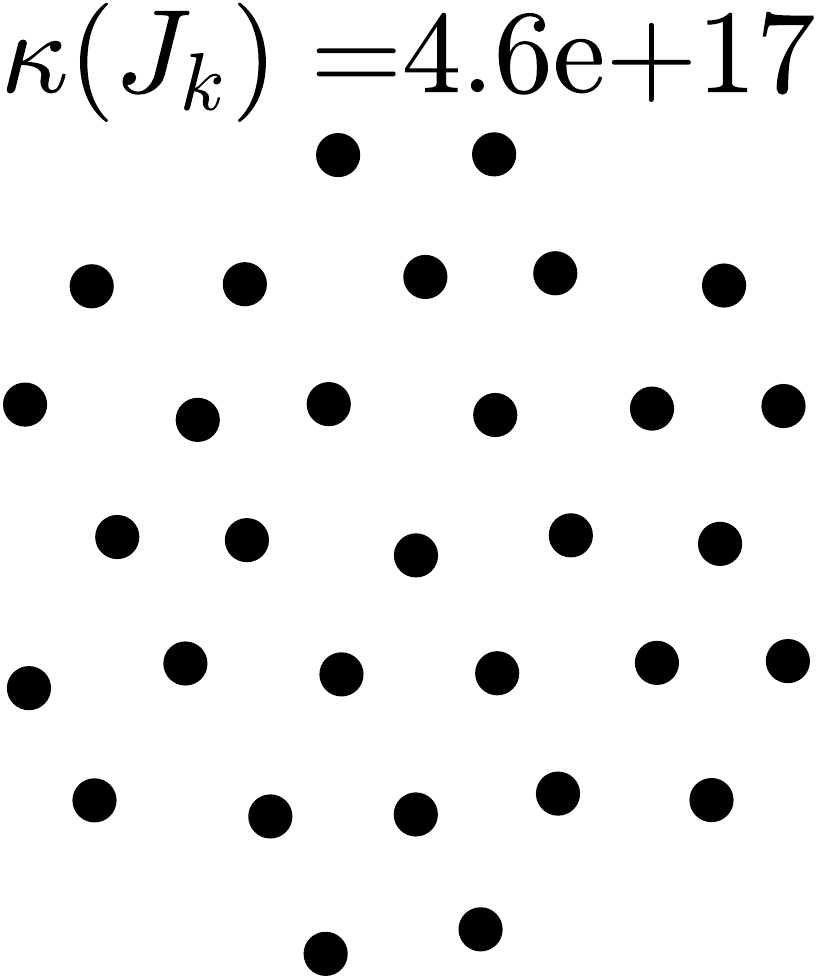}
			&\includegraphics[scale=0.23]{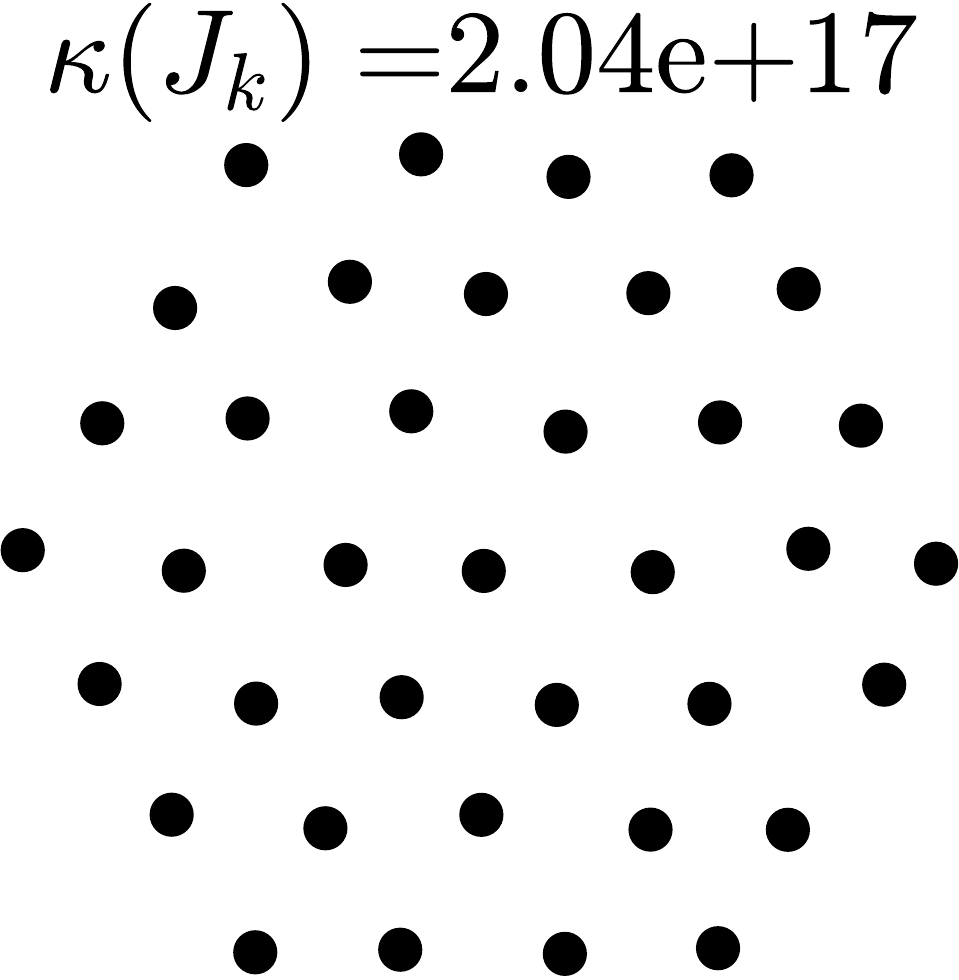}&\includegraphics[scale=0.23]{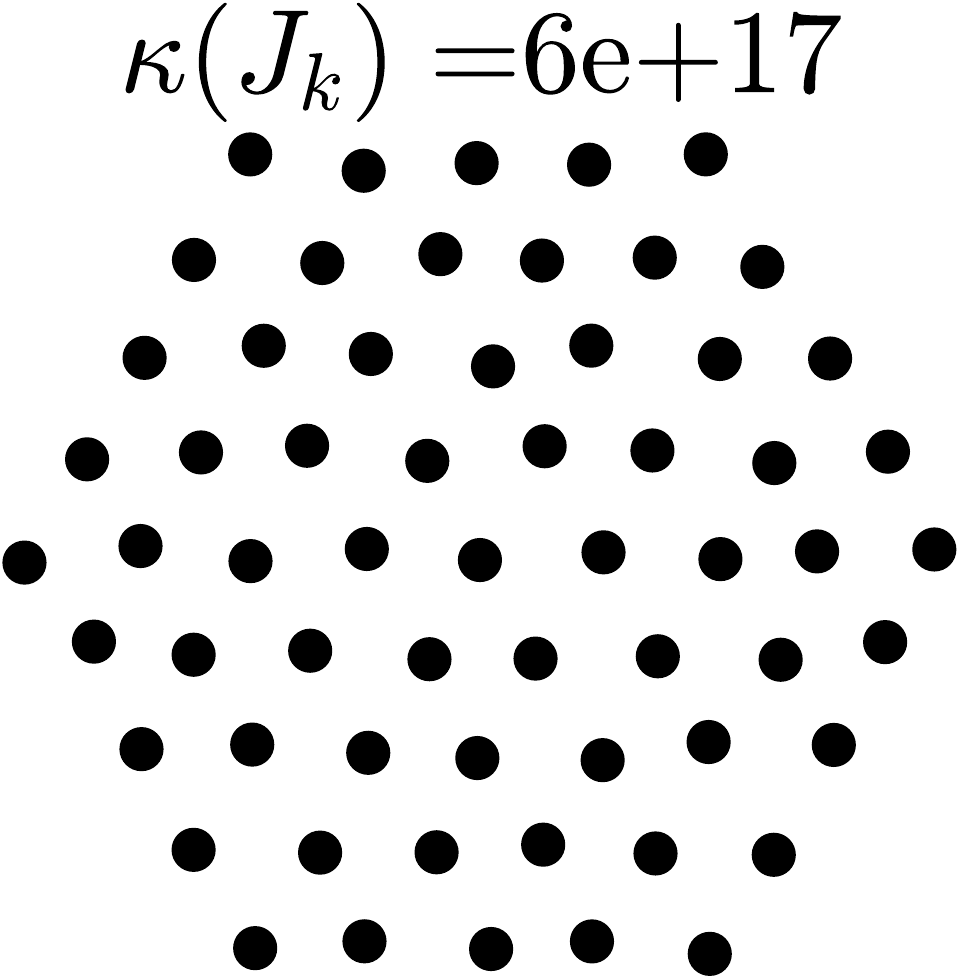}\\
			\hline
		\end{tabular}
	\end{center}
	\caption{Diagrams of some small stencils with the respective approximated condition number of $\JJ_k$, with $k=100$ and $h=\frac{2\pi}{6k}$. Top row: Stencils are taken from a regular hexagonal grid. Bottom row: perturbed position from stencils of top row.} \label{fig:stencils_hexagonal}
\end{figure}

As can be seen in Figures \ref{fig:stencils_square} and \ref{fig:stencils_hexagonal},
 the condition number of $\JJ_k$ may becomes very large, then we adopt to handle values into a computationally acceptable range and to use this fact for obtaining the regularization parameter $\beta$.
For small stencils we take $10^7 \leq\kappa_0\leq10^{14}$, 
with $\kappa_0=10^{7+\sqrt{n}}$ where $n$ is the size of the stencil,  and we take

\begin{equation}\label{eq:beta_recond}
    \beta=\frac{\lambda_1-\kappa_0  \lambda_n}{\kappa_0 -1},
\end{equation}
as regularization parameter, ensuring, from remark \ref{rmrk:better_coditioning_wideJk}, that $\kappa(\widetilde{\JJ}_k)\approx\kappa_0$, which is an adequate condition number to work in double precision. Now, $\kappa_0$ must be taken such that $\kappa_{\min}\leq \kappa_0\leq \kappa_{\max}$, thus 

\begin{equation}\label{betavalues}
   0\leq \frac{\lambda_1-\kappa_{\max}  \lambda_n}{\kappa_{\max} -1}\leq \beta \leq \frac{\lambda_1-\kappa_{\min}  \lambda_n}{\kappa_{\min} -1}.
\end{equation}

It is important to remark that $\beta$ is recalculated for every node, so its values may have a wide range; in our numerical tests and according to \eqref{betavalues}, has been about $0\le\beta\le 10^{-6}$.

We have noted empirically that the matrix $\JJ_k$ is worse conditioned for stencils with nodes collocated symmetrically on a regular grid, e.g. with square and hexagonal grids. See figures \ref{fig:stencils_square} and \ref{fig:stencils_hexagonal}, where we can observe that  severe ill-conditioning begins with symmetric stencils of $13$ nodes. However, in this case, with a small perturbation in the position of the nodes, its associated interpolation matrix $\JJ_k$ has a better condition number.

For numerical tests we considered the solutions for Helmholtz equation, $u_1$ and $u_2$, given by
\begin{equation}\label{eq:sol_imped1}
u_1(x,y)=\sqrt{k}H^{(1)}_0(k\sqrt{(x-2)^2+(y-2)^2}),
\end{equation}
and
\begin{equation}\label{eq:sol_imped2}
\small
\begin{array}{rcl}
u_2(x,y)&=&\sqrt{k}H^{(1)}_0(k\sqrt{(x+20)^2+(y+20)^2})+2\sqrt{k}H^{(1)}_0(k\sqrt{(x-20)^2+(y-20)^2})\\
&+&0.5\sqrt{k}H^{(1)}_0(k\sqrt{(x+20)^2+(y-20)^2})-\sqrt{k}H^{(1)}_0(k\sqrt{(x-20)^2+(y+20)^2}),
\end{array}
\end{equation}
\normalsize
where $H_0^{(1)}$ is the Hankel function of the first kind and $k$ is a constant wavenumber. The solution $u_1$ corresponds to the solution of the single source problem located at $\x_s=(2,2)$, whereas the solution $u_2$ corresponds to a solution of the problem with four single sources located at $\x_{s_1}=(-20,-20)$, $\x_{s_2}=(20,20)$, $\x_{s_3}=(-20,20)$, and $\x_{s_4}=(20,-20)$.
In the first test, we considered the function $u_2$ to see the behavior of the local truncation error of the approximation $(\partial_x)_{S,k}$. In addition we validate the conditioning of the matrix $\JJ_k$ and the better conditioning of $\widetilde{\JJ}_k$. 
With stencils $S$ as in Fig. \ref{fig:stencils_square} and Fig. \ref{fig:stencils_hexagonal} we approximate the solutions of the system $\JJ_k\bm{\alpha}=U$, with $U=u_2|_S$, by using  LDL$^T$, MDI, and ITMDI with 15 iterations ($M=15$ in \eqref{eq:ITMDI}). We used these results to compute the approximation $(\partial_x)_{S,k}u_2(\x)$. In tables \ref{table:Jk_bar_square} and \ref{table:Jk_bar_hexagonal} we can see a comparison of the relative local truncation error, given by $|(\partial_x)_{S,k}u_2(\x)-\partial_xu_2(\x)|/|\partial_x  u_2(\x)|$.
\begin{table}[ht]
	\begin{center}
		\begin{tabular}{c|cc|ccc}
			\hline
			\rowcolor{gray!20} Stencil size ($n$)  &  $\kappa(\JJ_k)$ & $\kappa(\widetilde{\JJ}_k)$ & Error LDL$^T$ & Error MDI  & Error ITMDI \\
			\hline\hline 
	         9 & 7.8e+14 & 7.8e+08 & 0.0036 & 0.0036 & 0.0036\\
	        13 & 1.34e+17 & 1.34e+11 & 0.000138 & 0.000139 & 0.000138\\
	        25 & 2.38e+18 & 2.38e+12 & 1.46e-07 & 7.93e-07 & 7.72e-07\\
	        49 & 1.61e+18 & 1.61e+12 & 6.44e-07 & 4.28e-08 & 5.79e-09\\
	        69 & 8.93e+17 & 8.93e+11 & 2.57e-06 & 1.14e-07 & 7.18e-09\\
	        81 & 1.44e+18 & 1.44e+12 & 3.98e-07 & 7.06e-08 & 1.45e-08
	       \\ \hline
		\end{tabular}
		\\
	\end{center} 
	\caption{This table shows values of the condition number of the matrices $\JJ_k$ and  $\widetilde{\JJ}_k$ for several stencil sizes $n$, and relative local truncation errors of the approximation $(\partial_{S,k})u_2(\x)\approx\partial_x u_2(\x)$. We used stencils as in Fig. \ref{fig:stencils_square}. We compared relative errors $|(\partial_{S,k})u_2(\x)-\partial_x u_2(\x)|/|\partial_x u_2(\x)|$ produced by using LDL$^T$, MDI and   ITMDI with $15$ iterations. The function $u_2$ as in \eqref{eq:sol_imped2}.}\label{table:Jk_bar_square}	
\end{table}

\begin{table}[ht]
	\begin{center}
		\begin{tabular}{c|cc|ccc}
			\rowcolor{gray!20} Stencil size ($n$)  &  $\kappa(\JJ_k)$ & $\kappa(\widetilde{\JJ}_k)$ & Error LDL$^T$ & Error MDI  & Error ITMDI \\
			\hline\hline 
		7 & 9.79e+10 & 9.79e+04 & 0.00311 & 0.00437 & 0.00311\\
		13 & 1.9e+16 & 1.9e+10 & 3.42e-05 & 3.7e-05 & 3.42e-05\\
		19 & 2.45e+16 & 2.45e+10 & 1.48e-08 & 8.36e-05 & 4.84e-05\\
		31 & 1.61e+17 & 1.62e+11 & 1.15e-06 & 4.33e-07 & 1.18e-07\\ 
		37 & 1.71e+17 & 1.71e+11 & 9.29e-08 & 3.23e-08 & 2.69e-08\\
		61 & 9.75e+18 & 9.78e+12 & 1.07e-07 & 9.56e-09 & 5.28e-09\\
		\hline
		\end{tabular}
		\\
	\end{center} 
		\caption{This table shows values of the condition number of the matrices $\JJ_k$ and  $\widetilde{\JJ}_k$, respect to the stencil size $n$, and relative local truncation errors of the approximation $(\partial_{S,k})u_2(\x)\approx\partial_x u_2(\x)$. We used stencils as in Fig. \ref{fig:stencils_hexagonal}. We compared relative errors $|(\partial_{S,k})u_2(\x)-\partial_x u_2(\x)|/|\partial_x u_2(\x)|$ produced by using LDL$^T$, MDI and   ITMDI with $15$ iterations. We used the $u_2$  in \eqref{eq:sol_imped2}.}\label{table:Jk_bar_hexagonal}		
\end{table}

\subsection{Pollution-effect and convergence}
To see the impact of the pollution effect in numerical solutions, the standard procedure is:  to compare the errors in several solutions obtained by increasing the wavenumber and keeping constant the Number of nodes Per Wavelength (NPW), i.e., the product $hk=\frac{2\pi}{Ng}$ should be constant \cite{chen9}.

\subsubsection{Test 1}

In this test we calculate the approximated solution for the problem

    \begin{equation}\label{eq:helmholtz_impedance_TEST1}
    \left\{
    \begin{array}{rcll}
    -\Delta u(\x)-k^2u(\x)&=&0, &\mbox{ in } \Omega \\
    \frac{\partial}{\partial\n}u(\x)+\ri k u(\x)&=&g(\x), & \mbox{ on } \Gamma=\partial\Omega
    \end{array}
    \right.
    \end{equation}
     with the known data $g(\x)$, $\Omega=(-0.5,0.5)\times(-0.5,0.5)$ and $k=\omega c^{-1}$, with constant wave speed $c\equiv 1$.   Results are verified with solutions $u_1$ and $u_2$ in \eqref{eq:sol_imped1} and \eqref{eq:sol_imped2}, respectively. 
    Tables \ref{table:fang1square_pollution},  \ref{table:fang1hexagon_pollution} and \ref{table:fang1hexagon19_pollution} show errors when the resolution is kept constant at $N_g=6$ NPW. In calculations we have taken the perturbed matrices, $\widetilde{\JJ}_k$, such that $\kappa_0=\kappa(\widetilde{\JJ}_k)=10^{-6}\kappa(\JJ_k)$. In the three cases, for uniform square and hexagonal grids, we see that the order of the error remains constant, i.e., $\|u_1-\widetilde{u}_1\|_{\infty}\sim \bO(1)$ and  $\|u_2-\widetilde{u}_2\|_{\infty}\sim \bO(1)$, as $h\rightarrow 0$ with $hk=\frac{2\pi}{6}$. Hence, in these examples, pollution effects are mitigated.

\begin{figure}[ht]
	\begin{center}
		\includegraphics[width=16.0cm]{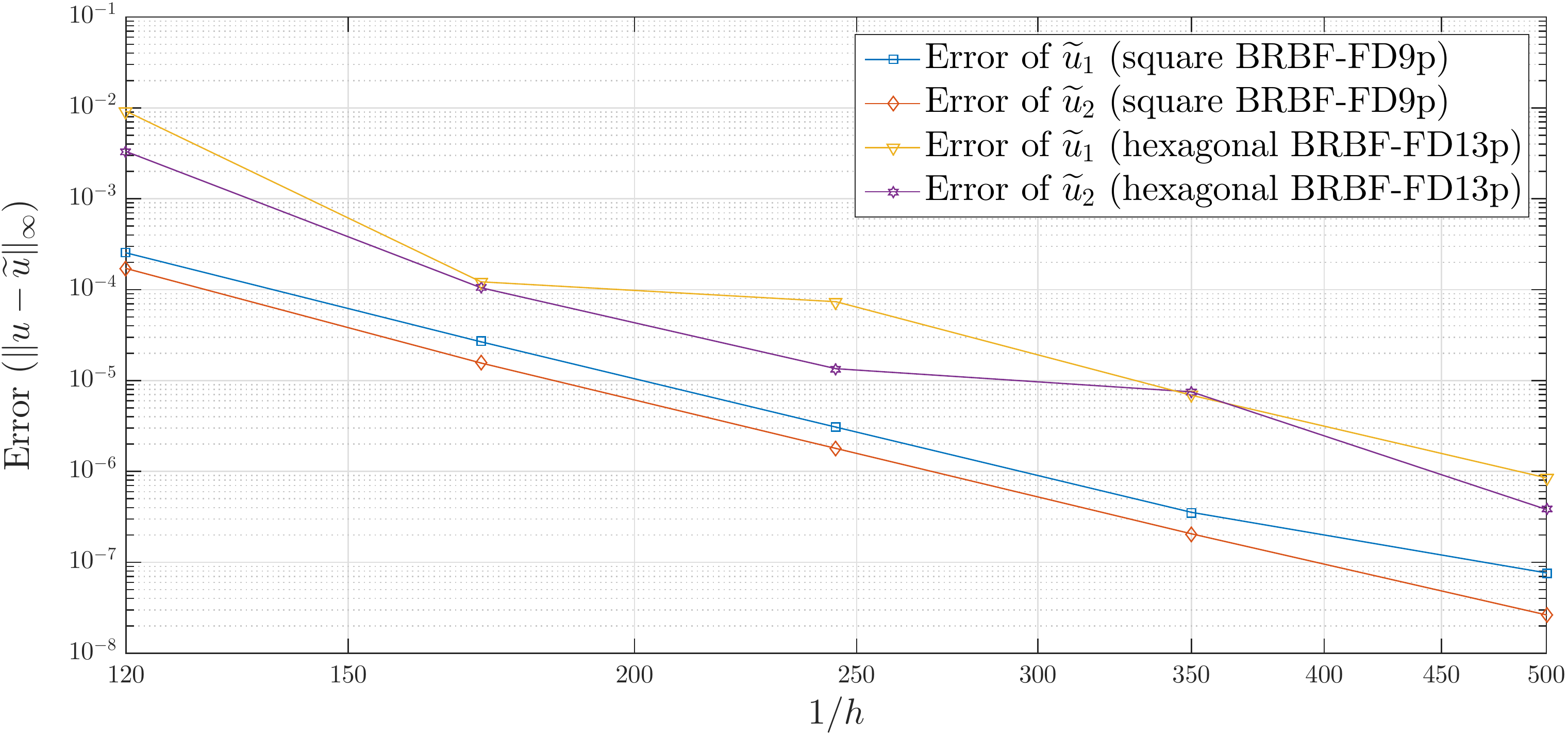}          
		\caption{Plot of results in tables \ref{table:fang1convergence_square} and \ref{table:fang1convergence_hexa}.}
		\label{fig:convergence_RBF}
	\end{center}
\end{figure}

     Results of convergence tests are summarized in tables \ref{table:fang1convergence_square} and \ref{table:fang1convergence_hexa} and Fig \ref{fig:convergence_RBF}. For these tests we choose $\widetilde{\JJ}_k=\beta\mathbf{I}+\JJ_k$, with $\beta$ according to \eqref{eq:beta_recond},  such that the condition number  $\kappa_0=\kappa(\widetilde{\JJ}_k)=10^{-4}\kappa(\JJ_k)$. 

\begin{table}[ht]
	\begin{center}
		\begin{tabular}{cccc|cc}
			\rowcolor{gray!20} $\frac{k}{2\pi}$  &  $\frac{1}{h}$ & Nodes ($N$) & $\kappa(\mathbf{H})$ & $\|u_1-\widetilde{u}_1\|_{\infty}$ & $\|u_2-\widetilde{u}_2\|_{\infty}$   \\
			\hline\hline 
	    	10 & 60 & 3721 & 1.40e+04 & 1.97e-04 & 1.79e-04       \\  
	    	20 & 120 & 14641 & 7.54e+04 & 1.95e-04  & 1.56e-04    \\
			40 & 240 & 58081 & 4.17e+05 & 1.96e-04 & 1.10e-04     \\
			80 & 480 & 231361 & 2.38e+06 & 1.98e-04  & 1.72e-04   \\
			120 & 720 & 519841 & 6.55e+06 & 1.96e-04  & 1.21e-04   \\
			\hline
		\end{tabular}
		\\
	\end{center} 
	\caption{Results for approximated solutions of \eqref{eq:helmholtz_impedance_TEST1}. We used a square uniform grid in $\Omega\cap\partial\Omega$. For inner nodes the stencil size is $n=13$, at boundary nodes $n_b=15$, the number of nodes per wavelength is kept constant with $N_g=6$.}\label{table:fang1square_pollution}

\end{table}
\begin{table}[ht]
	\begin{center}
		\begin{tabular}{cccc|cc}
			\rowcolor{gray!20} $\frac{k}{2\pi}$  &  $\frac{1}{h}$ & Nodes ($N$) & $\kappa(\mathbf{H})$ & $\|u_1-\widetilde{u}_1\|_{\infty}$ & $\|u_2-\widetilde{u}_2\|_{\infty}$   \\
			\hline\hline 
			10 & 60  & 4237   & 2.52e+04 & 1.32e-04 & 3.30e-05  \\  
	    	20 & 120 & 16752  & 1.09e+05 & 4.49e-05 & 3.35e-05  \\
	     	40 & 240 & 66861  & 5.25e+05 & 9.30e-05 & 3.59e-05   \\
			80 & 480 & 266680 & 3.46e+06 & 6.69e-05 & 5.62e-05\\
			120 & 720 & 599458 & 8.21e+06 & 9.05e-05 & 4.25e-05 \\
		 \hline
		\end{tabular}
		\\
	\end{center} 
	\caption{Results for approximated solutions of \eqref{eq:helmholtz_impedance_TEST1}. We used an hexagonal uniform grid in $\Omega\cap\partial\Omega$. For inner nodes the stencil size is $n=13$, at boundary nodes $n_b=25$, the number of nodes per wavelength is kept constant with $N_g=6$. }\label{table:fang1hexagon_pollution}
\end{table}

\begin{table}[ht]
	\begin{center}
		\begin{tabular}{cccc|cc}
			\rowcolor{gray!20} $\frac{k}{2\pi}$  &  $\frac{1}{h}$ & Nodes ($N$) & $\kappa(\mathbf{H})$ & $\|u_1-\widetilde{u}_1\|_{\infty}$ & $\|u_2-\widetilde{u}_2\|_{\infty}$   \\
			\hline\hline 
		    10 & 60 & 4237     & 8.18e+04 & 1.97e-05 & 1.23e-05   \\  
			20 & 120 & 16752   & 8.58e+05 & 2.22e-05 & 1.77e-05   \\
			40 & 240 & 66861   & 6.12e+05 & 1.94e-05 & 1.47e-05   \\
			80 & 480 & 266680  & 5.22e+06 & 1.78e-05 & 1.28e-05   \\
			120 & 720 & 599458 & 9.53e+06 & 1.85e-05  & 9.86e-06  \\
			\hline
		\end{tabular}
		\\  
	\end{center} 
	\caption{Results for approximated solutions of \eqref{eq:helmholtz_impedance_TEST1}. We used a square uniform grid in $\Omega\cap\partial\Omega$. For inner nodes the stencil size is $n=19$, at boundary nodes $n_b=25$, the number of nodes per wavelength is kept constant with $N_g=6$. } \label{table:fang1hexagon19_pollution}		
\end{table}

\begin{table}[ht]
	\begin{center}
		\begin{tabular}{cccc|cc}
			\rowcolor{gray!20} NPW=$N_g$  &  $\frac{1}{h}$ & Nodes ($N$) & $\kappa(\mathbf{H})$ & $\|u_1-\widetilde{u}_1\|_{\infty}$ & $\|u_2-\widetilde{u}_2\|_{\infty}$   \\
			\hline\hline 
	6.0 & 120.0 & 14641 & 7.38e+04 & 2.54e-04 & 1.71e-04 \\
	8.6 & 171.4 & 29584 & 1.03e+05 & 2.67e-05 & 1.56e-05  \\
	12.2 & 244.9 & 60025 & 1.47e+05 & 3.07e-06 & 1.79e-06 \\
	17.5 & 350.0 & 122500 & 2.06e+05 & 3.55e-07 & 2.06e-07 \\
	25.0 & 500.0 & 251001 & 2.75e+09 & 7.65e-08 & 2.64e-08  \\
			\hline
		\end{tabular}
		\\
	\end{center} 
	\caption{Here $\frac{k}{2\pi}=20$. In a square uniform grid, we took stencils of size: $n=9$ for inner nodes,  and $n_b=15$ for boundary nodes. By applying a linear regression we have that $\log_{10}(\|u_1-\widetilde{u}_1\|_{\infty})\approx5.76\log_{10}(h)+8.3$ and $\log_{10}(\|u_2-\widetilde{u}_2\|_{\infty})\approx6.132\log_{10}(h)+8.932$. }	
		\label{table:fang1convergence_square}	
\end{table}

\begin{table}[ht]
	\begin{center}
		\begin{tabular}{cccc|cc}
			\rowcolor{gray!20} NPW=$N_g$  &  $\frac{1}{h}$ & Nodes ($N$) & $\kappa(\mathbf{H})$ & $\|u_1-\widetilde{u}_1\|_{\infty}$ & $\|u_2-\widetilde{u}_2\|_{\infty}$   \\
			\hline\hline 
	    6.0 & 120.0 & 16752 & 1.34e+24 & 9.18e-03 & 3.31e-03\\
	    8.6 & 171.4 & 33960 & 9.84e+05 & 1.22e-04 & 1.05e-04\\
	    12.2 & 244.9 & 69338 & 1.69e+16 & 7.37e-05 & 1.35e-05\\
	    17.5 & 350.0 & 141753 & 4.04e+32 & 6.89e-06&7.51e-06\\
	    	25.0 & 500.0 & 289291 & 2.14e+20 & 8.50e-07 & 3.82e-07   \\
			\hline
		\end{tabular}
		\\
	\end{center} 
	\caption{Here $\frac{k}{2\pi}=20$. In a hexagonal grid, we took stencils of size:  $n=13$ for inner nodes, and  $n_b=25$ for boundary nodes. By applying a linear regression we have that $\log_{10}(\|u_1-\widetilde{u}_1\|_{\infty})\approx6.01\log_{10}(h)+10.01$ and $\log_{10}(\|u_2-\widetilde{u}_2\|_{\infty})\approx5.82\log_{10}(h)+9.33$. 
		} \label{table:fang1convergence_hexa}		
\end{table}

\subsubsection{Test 2}

In this example we consider the problem

\begin{equation}
\begin{cases}\label{eq:abc_problem}
-\Delta u(x,y)-k^2u(x,y)=0& \mbox{ in } \Omega\\
\frac{\partial u}{\partial \mathbf{n}}u(x,y)+\ri k u(x,y)=g(x,y)&  \mbox{ on } \partial\Omega,
\end{cases}
\end{equation}
with $\Omega=(0,1)\times(0,1)$, whose solution is given by the plane wave $u(x,y;k,\theta)=e^{\ri k(x\cos\theta+y\sin\theta)}$  when the data $g$ on the impedance boundary condition is given by
\begin{equation*}
g(x,y)=
\begin{cases}
\ri (k-k_2)e^{\ri k_1x} & \mbox{ if } x\in (0,1) \mbox{ and } y=0\\
\ri (k+k_1)e^{\ri (k_1+k_2y)}& \mbox{ if } x=1 \mbox{ and } y\in(0,1)\\
\ri (k+k_2)e^{\ri (k_1x+k_2)}& \mbox{ if } x\in (0,1) \mbox{ and } y=1\\
\ri (k-k_1)e^{\ri k_2y}& \mbox{ if }  x=0 \mbox{ and } y\in(0,1),
\end{cases}
\end{equation*}
with $k_1=k\cos\theta$ and $k_2=k\sin\theta$. This example is standard for testing  numerical dispersion of solvers for Helmholtz equation. The approximated solutions for this problem ware calculated using Gaussian RBF-FD on hexagonal grids with 7-stencil (GRBF-FD-7p), Bessel RBF-FD with 9-stencils (BRBF-FD-9p) and 13-stencils (BRBF-FD-13p) on uniform Cartesian grids and we compare with results reported in
\cite{chen9}. In all methods it has been fixed NPW$=2\pi$to keep constant resolution when the
wavenumber k is increasing. For GRBF-FD-7p we have used the approximations in \cite{londono2019}
with its respective shape parameter $\ve_{op}$. We use it to approximate the Laplace operator
and all partial derivative operators involved in the boundary condition. To solve the local interpolations in BRBF-FD9p and BRBF-FD13p we used condition numbers $\kappa_0=10^{-4}\kappa(\JJ_k)$ and $\kappa_0=10^{-6}\kappa(\JJ_k)$, respectively. 

\begin{figure}[ht]
	\begin{center}
		\begin{tabular}{cc}			
			\includegraphics[width=7.7cm]{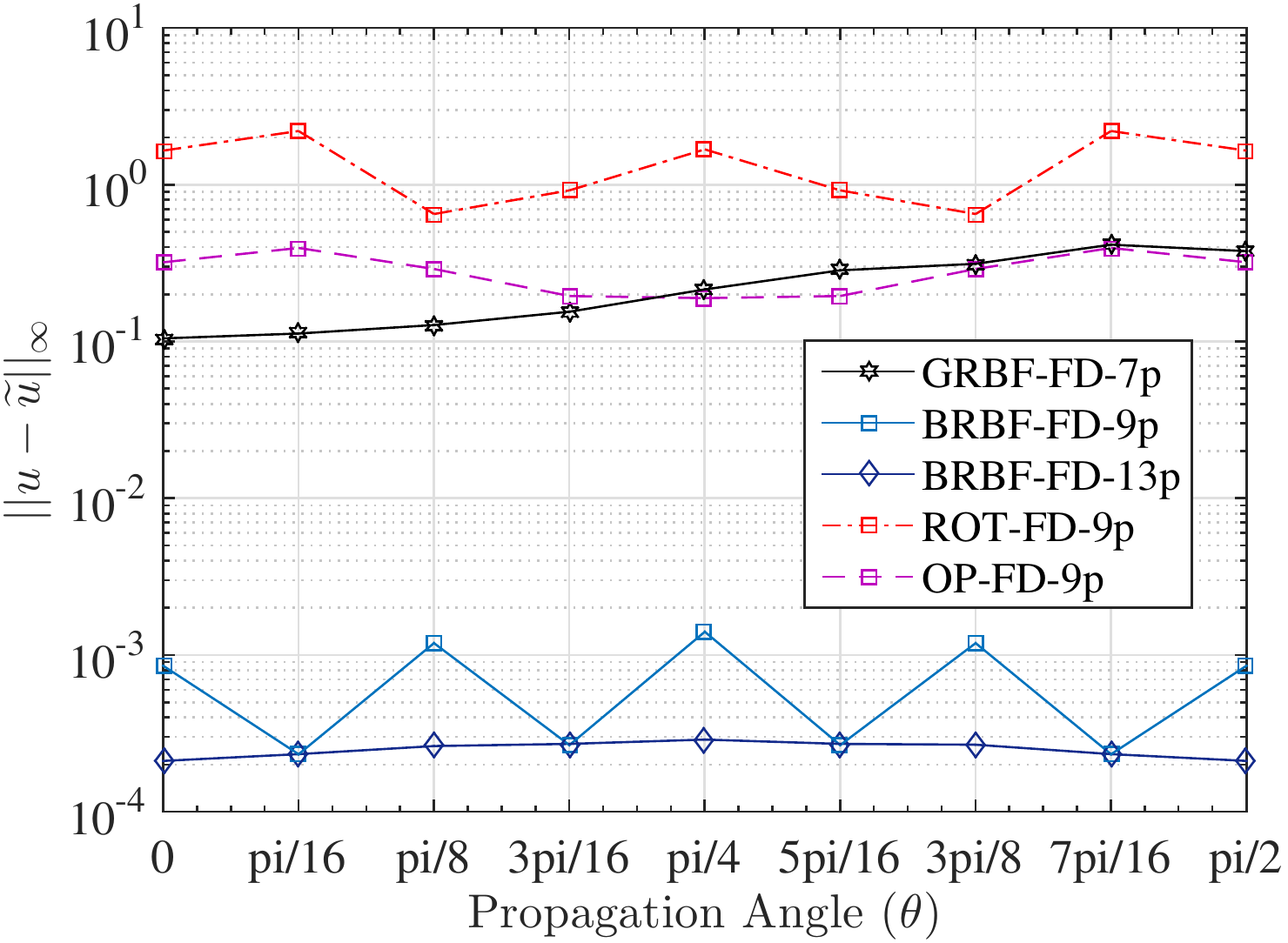} &          \includegraphics[width=7.7cm]{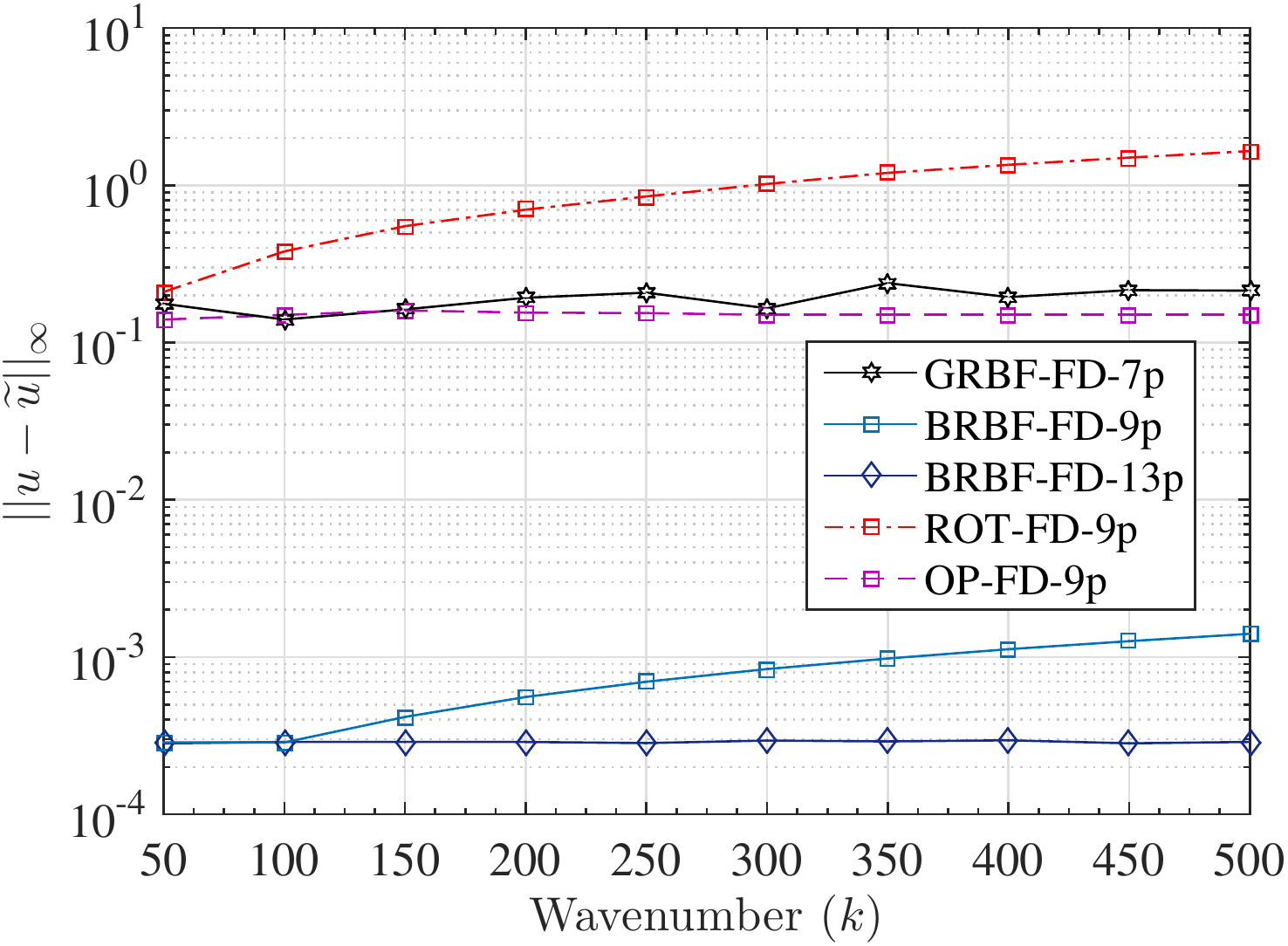}  			
		\end{tabular}       	
		\caption{Comparison of results among: BRBF-FD9p, BRBF-FD13p with square stencils,  GRBF-FD7p with hexagonal stencils, and ROT-FD9p and OP-FD9p from \cite{chen9}. (Left) Results for $k=500$ and $h=1/500$, varying the propagation angle. (Right) With  $\theta=\pi/4$ and $h=1/k$, varying the wavenumber $k$.}
		\label{fig:testpoll}
	\end{center}
\end{figure}

 Results and comparisons can be seen in Fig. \ref{fig:testpoll}. On the left frame we can see that errors of BRBF-FD9p and BRBF-FD13p are smaller than GRBF-FD7p, ROT-FD9p and OP-FD9p, at least in two orders of magnitude.  Besides, we see that there is less anisotropy in the error of BRBF-FD13p, where over all propagation angles we have improved the error at least three orders of magnitude. We can see that the behavior of GRBF-FD7p is similar to OP-FD9p. On the right frame we can see that the dispersion and pollution effects are mitigated with BRBF-FD13p because, when the wavenumber increases while we keep a fixed resolution with $kh=\frac{2\pi}{Ng}=1$, the error remains almost constant.   
 
\section{Examples with some Helmholtz problems}\label{sec:numerical_tests}

In this section we test our  BRBF-FD scheme by computing numerical solutions of some Helmholtz problems  with second and third order Absorbing Boundary Conditions (ABC) \cite{majda1}.

\subsection{Approximated fundamental solutions}

It is known that the problem $-\Delta u-k^2u=\delta(\x)$ in the free-space $\R^2$ has a unique solution when it is imposed the Sommerfeld radiation condition
\begin{equation*}
\lim_{\|\x\|\rightarrow\infty}\|\x\|^\frac{1}{2}\left(\frac{\partial}{\partial r}-\ri k\right)u(\x)=0.
\end{equation*}
Particularly, the associated Green's function, which is solution of $-\Delta u(\x)-k^2u(\x)=\delta(\x-\x_0)$,  \cite{Ciraolo2009_radiation_condition_for_the_2D_helmholtz_equation}  is given by $u(\x)=G(\x,\x_0)$, with the fundamental solution
\begin{equation*}
G(\x,\x_0)=\frac{\ri}{4}H_0^{(1)}(k\|\x-\x_0\|).
	\end{equation*}

\begin{figure}[ht]
	\includegraphics[scale=0.8]{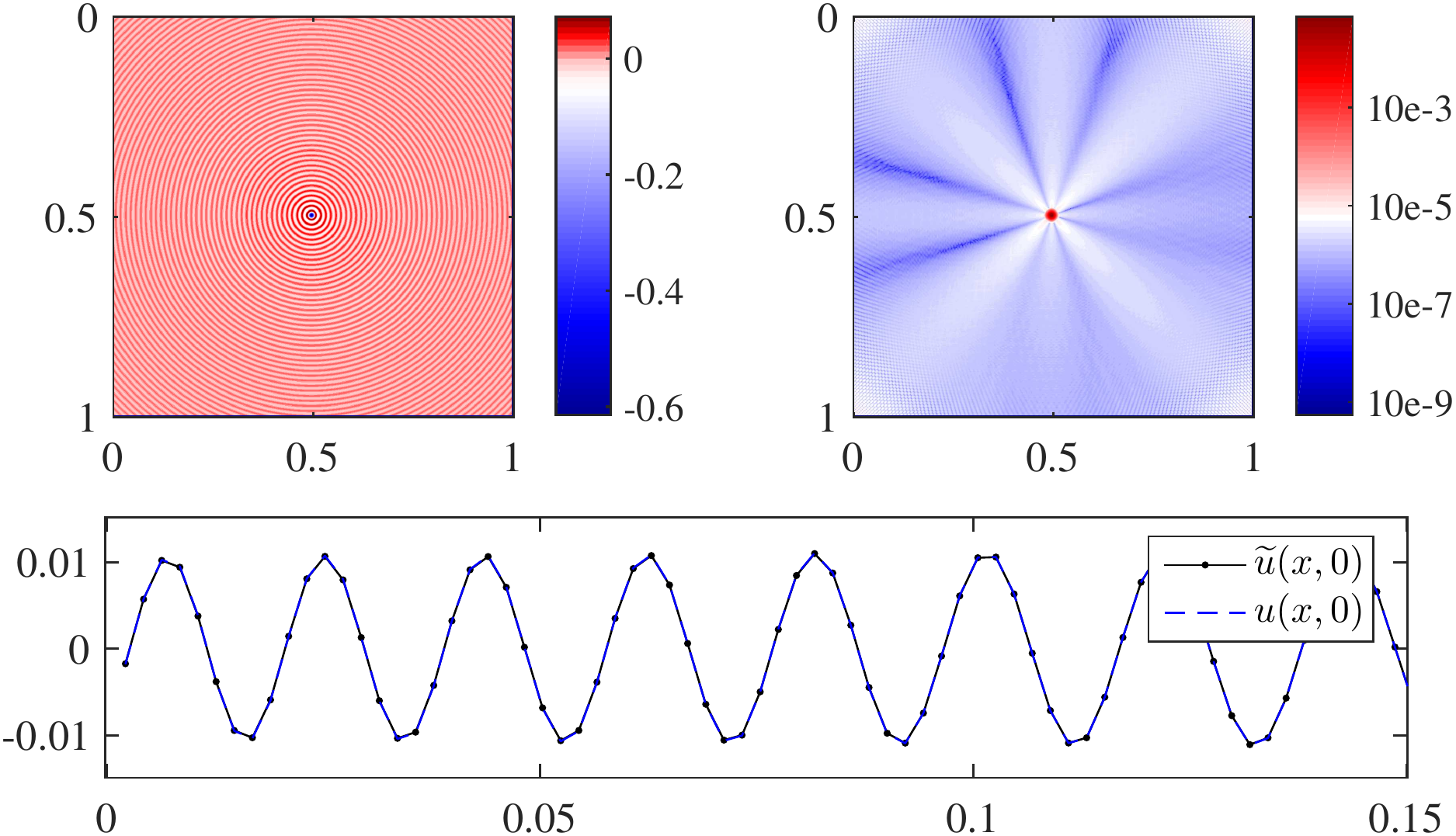}
	\caption{(Top-left), approximated solution $u(\x)$ (top-right) plot of  $|\widetilde{u}(\x)-u(\x)|$, (bottom) comparison on boundary values. Source at $\x_0=(0.5,0.5)$.}\label{fig:hankel_c_500}
\end{figure}

\begin{figure}[ht]
		\includegraphics[scale=0.8]{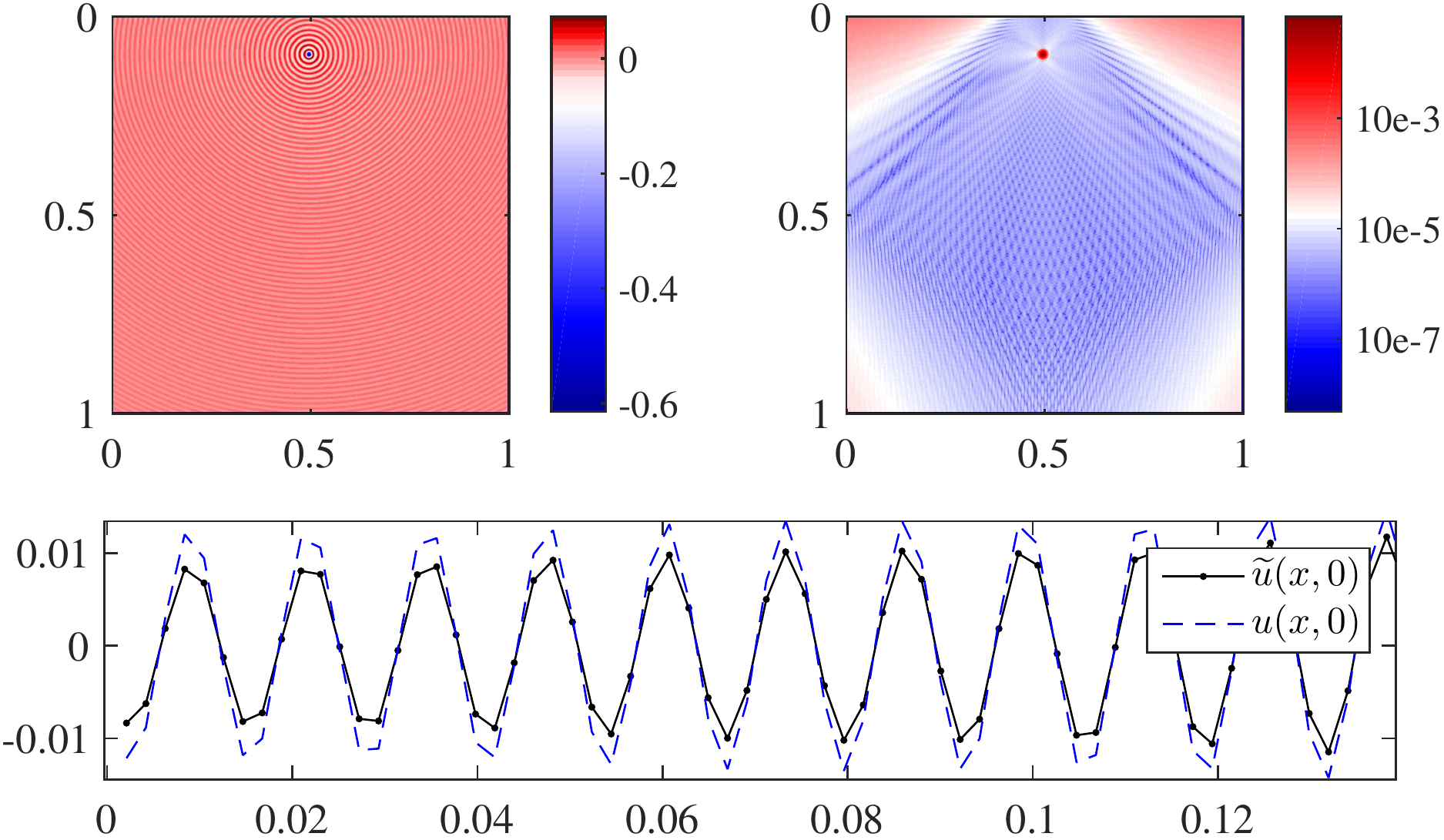}
		\caption{(Top-left) approximated solution $\widetilde{u}(\x)$. (Top-right) plot of  $|\widetilde{u}(\x)-u(\x)|$. (Bottom) Zoom of comparison on boundary values. The source is located at $\x_0=(0.5,0.1)$.}\label{fig:hankel_m_500}
\end{figure}

We compute approximated Green's functions in the free space truncated to a bounded domain $\Omega\subset\R^2$ for $\x\neq\x_0$, where $\x_0\in\Omega$,  through the boundary value problem
 \begin{equation}\label{eq:helmholtz_impedance}
 \left\{
 \begin{array}{rcll}
 -\Delta u(\x)-k^2u(\x)&=&\widetilde{\delta}(\x-\x_0), &\mbox{ in } \Omega\\
 \frac{\partial}{\partial\n}u(\x)+\ri k \mathcal{B}u(\x)&=&0, & \mbox{ on } \Gamma=\partial\Omega,
 \end{array}
 \right.
 \end{equation}
 where $\mathcal{B}=1+\frac{3}{4k^2}\frac{\partial^2}{\partial\bm{\tau}^2}-\frac{\ri}{4k^3}\frac{\partial^3 }{\partial \mathbf{n}\partial\bm{\tau}^2 } $. The boundary condition corresponds to the ABC in the Padé approximation \cite{majda1}. 
  The single source is given by the Gaussian function
\begin{equation}\label{eq:single_source}
\widetilde{\delta}(\x-\x_0)=\frac{1}{2\pi\sigma^2}e^{-\frac{\|\x-\x_0\|^2}{2\sigma^2}},
\end{equation}
where $\sigma$ is a value such that $\int_{\Omega}\widetilde{\delta}(\x-\x_0)d\x\approx1$.
To solve \eqref{eq:helmholtz_impedance} we have used the BRFB-FD9p scheme  in $\Omega=(0,1)\times(0,1)$, by using square $9-$stencils at inner nodes, and $19$-stencils at boundary nodes. The results
shown in Fig. \ref{fig:hankel_c_500} and Fig. \ref{fig:hankel_m_500} were calculated with  $k=500$ and $N_g=6$ \ NPW, i.e., $h=\frac{2\pi}{6k}$. We point out that results show a good accuracy at $\x\neq\x_0$. Besides, for the source located at the center of the square domain, the wavelength of the numerical solution on the boundary, matches very close to the exact one. This is a good indication that dispersion errors are not significant. However in Fig. \ref{fig:hankel_m_500} (bottom) we see that amplitude has a considerable discrepancy with respect to the exact one, this is due to the approximation of ABC.


\subsection{Heterogeneous medium}

\subsubsection{Smooth medium}
\begin{figure}[ht]
	\begin{center}
		\begin{tabular}{ll}		
			\includegraphics[scale=0.45]{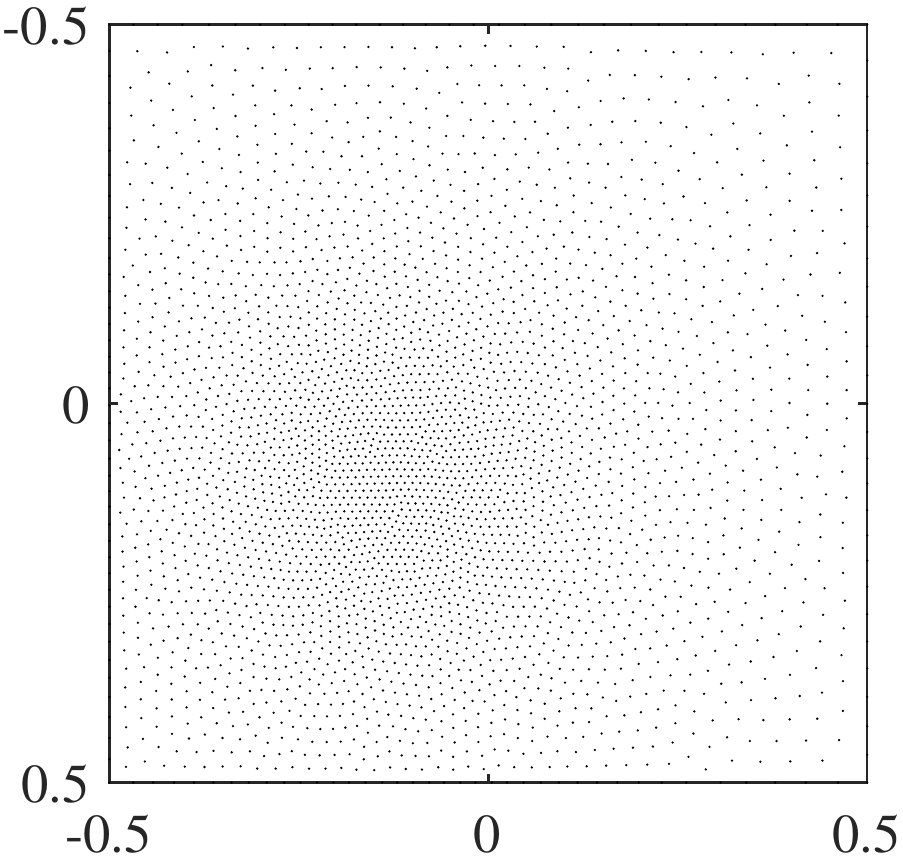} &\includegraphics[scale=0.6]{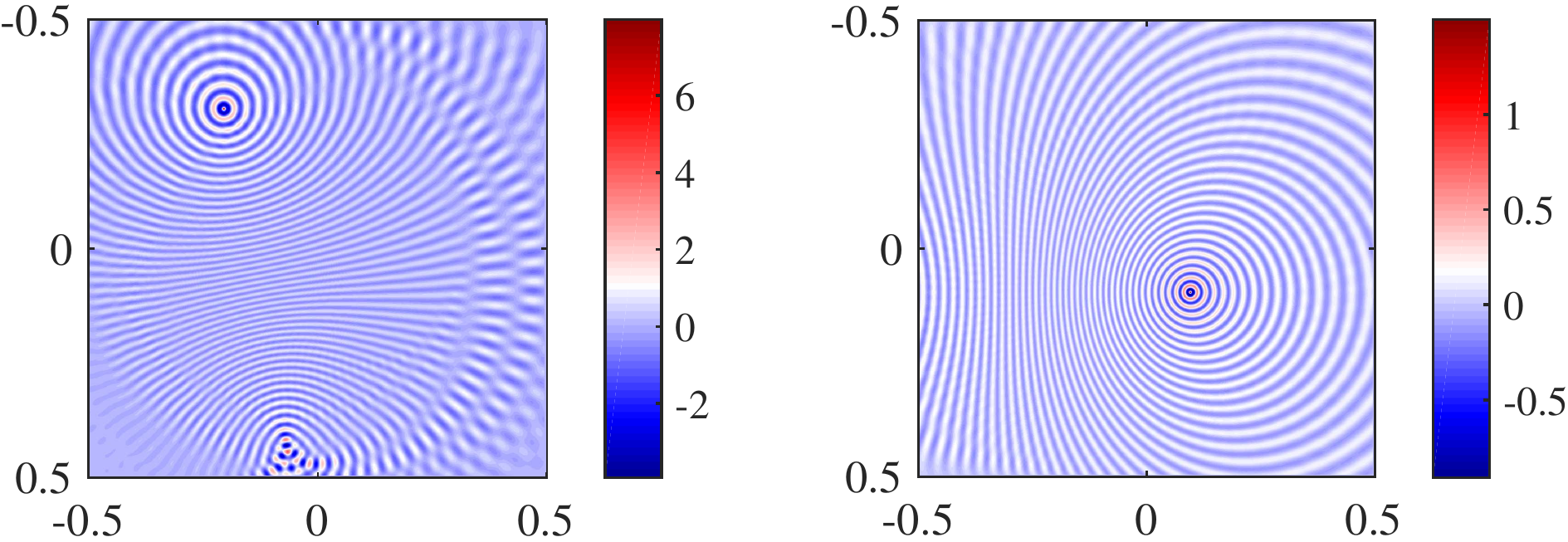} \\		
			\includegraphics[scale=0.45]{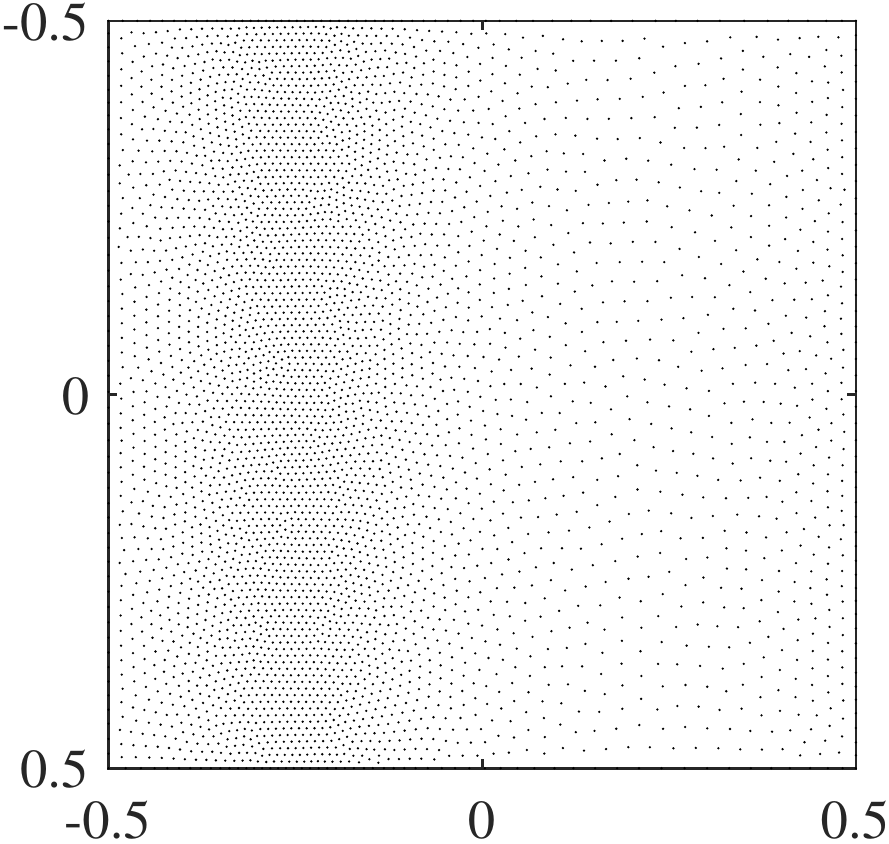} &\includegraphics[scale=0.6]{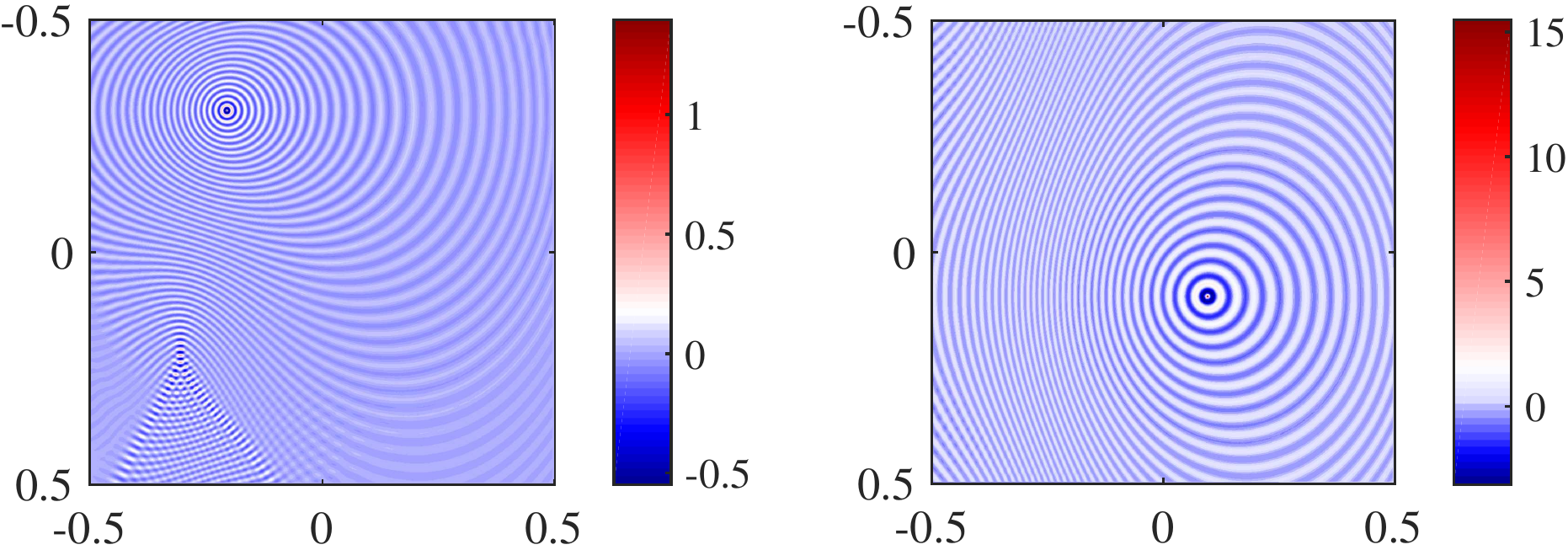} \\
		\end{tabular} \caption{Plots corresponding to velocity models in \eqref{eq:speed_fang3} (top row) and \eqref{eq:speed_fang4} (bottom row). (On left column) Sketch of node distributions according to their local wavelengths. Plots of the real part of solutions corresponding to single sources located at $(x,y)=(-0.2,-0.3)$ (central column) and $(x,y)=(0.1,0.1)$ (right column).}
		\label{fig:fang34}
	\end{center}
\end{figure}

For this qualitative test, we have calculated approximated solutions of the problem
\begin{equation}\label{eq:helmholtz_ABC2}
\left\{
\begin{array}{rcll}
-\Delta u(\x)-\omega^2c(\x)^{-2}u(\x)&=&\widetilde{\delta}(\x-\x_0), &\mbox{ in } \Omega\\
\frac{\partial}{\partial\n}u(\x)+\ri\omega c(\x)^{-1}\mathcal{B}u(\x)&=&0, & \mbox{ on } \Gamma=\partial\Omega
\end{array}
\right.
\end{equation}
where $\mathcal{B}$ is  the operator $\mathcal{B}=1+\frac{c(\x)^2}{2\omega^2}\frac{\partial^2}{\partial\bm{\tau}^2}$ corresponding to the ABC of second order and $\Omega=(-0.5,0.5)\times(-0.5,0.5)$. Here we perform two examples with the velocity functions
\begin{equation}
c(x,y)=3 - 2.5e^{  -( (x+0.125)^2 + (y-0.1)^2 )/0.8^2  }
\label{eq:speed_fang3}
\end{equation} 
and
\begin{equation}
c(x,y)=1+0.5\sin(2\pi x).
\label{eq:speed_fang4}
\end{equation}
For these velocity models, nodes distributions are sketched on the left column of Fig. \ref{fig:fang34}. On center and right columns it can be seen the real part of the approximated solution for two different single sources. Table \ref{table:fang34}  shows results for required times to assembly sparse matrices $\mathbf{H}=\mathbf{H}_{\Omega}+\mathbf{H}_{\Gamma}$ and for solution of the system $-\mathbf{H}U=\mathbf{F}$ by LU factorization.

\begin{table}[ht]
	\begin{center}
		\begin{tabular}{ccc|ccc}
			\rowcolor{gray!20} $\omega/2\pi$  & Nodes ($N$)  &  $\kappa(\mathbf{H})$  &Time (s) for $\mathbf{H}_{\Omega}$ & Time (s) for $\mathbf{H}_{\partial\Omega}$ & Time (s) for LU  \\
			\hline\hline 
			2.5 & 1072 & 1.09e+04 & 0.37 & 0.24 & 0.02 \\
			5 & 4404 & 1.46e+05 & 1.18 & 0.45 & 0.10  \\
			10 & 17563 & 6.25e+05 & 4.72 & 1.01 & 0.45  \\
			20 & 70585 & 1.93e+06 & 18.98 & 2.58 & 3.00 \\
			40 & 283458 & 1.27e+07 & 75.10 & 9.05 & 17.14   \\
			\hline
		\end{tabular}
		\\
	\end{center} 
	\caption{Results in computing solutions corresponding to the smooth velocity model in \eqref{eq:speed_fang3}.} \label{table:fang34}		
\end{table}

\subsubsection{Non smooth medium (Test in the 2004 BP velocity model)}
\begin{figure}[ht]
	\begin{center}
		\begin{tabular}{c}
			\includegraphics[scale=0.475]{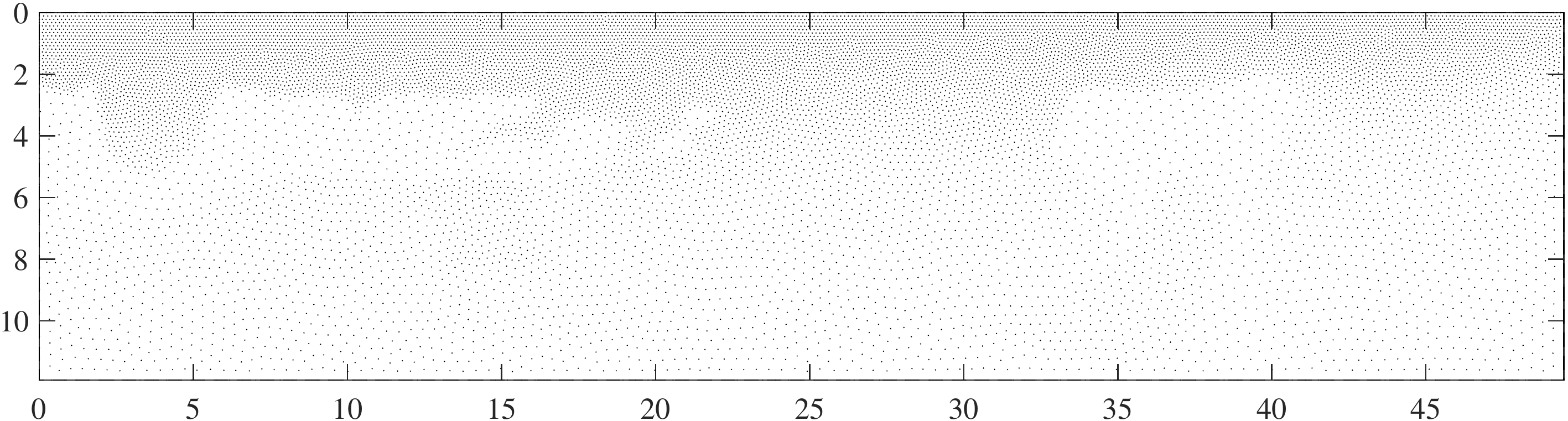}\\
			\includegraphics[scale=0.85]{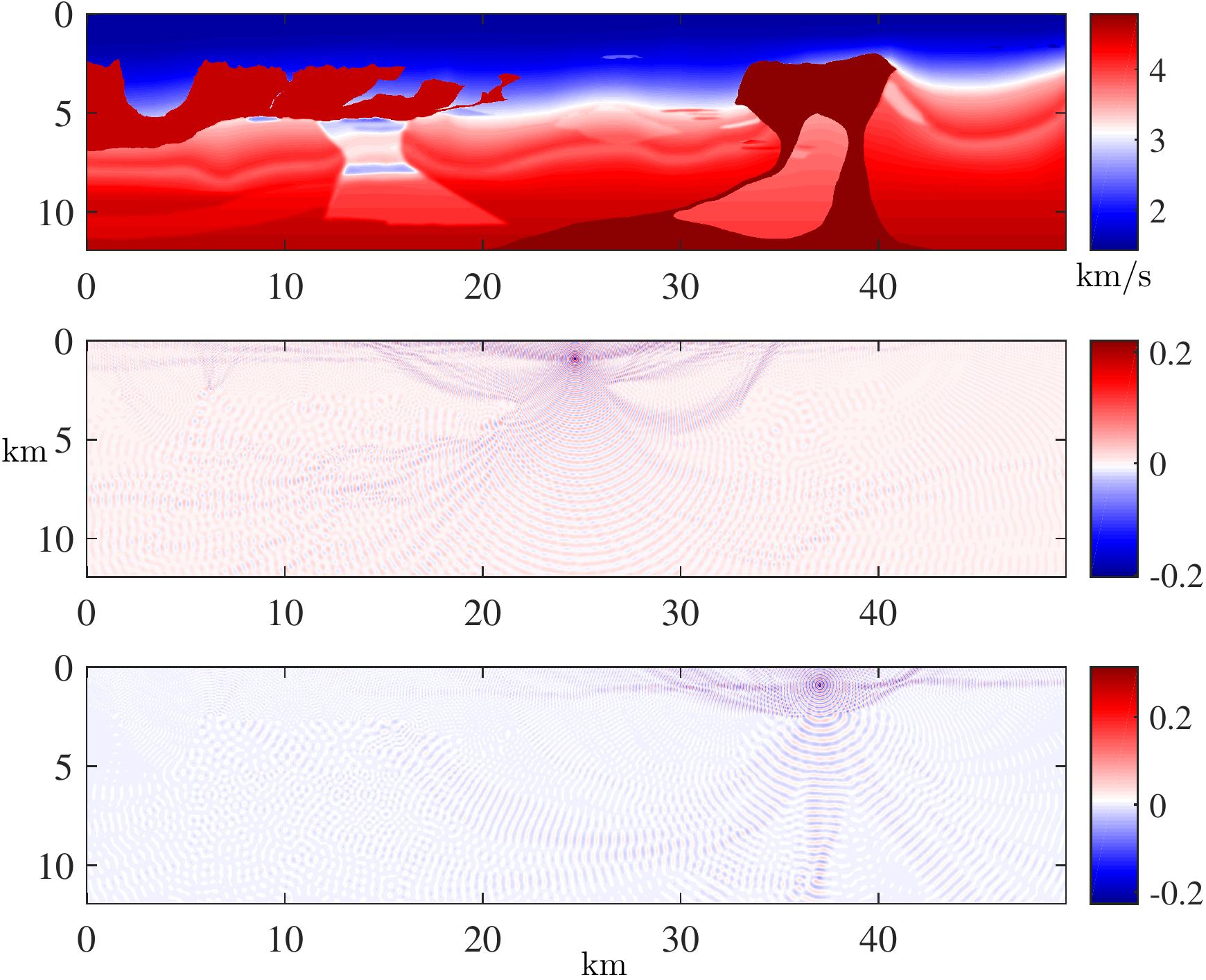}
		\end{tabular}
	\end{center}
	\caption{(Top) Sketch of node distribution for 2004 BP model, for a frequency at $\omega/2\pi=2$Hz with $Ng=10$ NPW. (Middle-top) velocity model.  (Middle-Bottom and bottom) plots of real part of the pressure wavefield, for a frequency of $\omega/2\pi=10$Hz, with two single sources located at $(x,y)=(24.74\,\mbox{km} ,1.00\,\mbox{km} )$ and $(x,y)=(37.11\,\mbox{km} ,1.00\,\mbox{km} )$. The domain is discretized with $N=345393$ nodes. The local interpolation is performed with $19$-stencils for inner nodes and $25$-stencils for boundary nodes.}\label{fig:2004bp_model}
\end{figure}

We consider the 2004 BP velocity benchmark, which is a popular model in research for velocity estimation methods in seismic imaging, which is presented as a challenge due to its complexity and large scale \cite{2004bp}. The  velocity model $c(\x)$ can be seen in the density plot on the middle-top in Fig. \ref{fig:2004bp_model}. Roughness of the velocity model such as hard interfaces and sharp transitions generates strong reflections that hinder the efficiency of known iterative methods due to the ill-conditioning of the matrix, when  the number of iterations increases \cite{ZEPEDA_NUNEZ2016347}. In addition, for large $\omega$, the interaction of high frequency waves with short wavelength structures such as discontinuities, increases the reflections, further deteriorating the convergence rate. In BRBF-FD schemes the local interpolation matrices $\JJ_k$ becomes dramatically ill-conditioned. However, in our tests, we have found empirically that in keeping the condition number of $\widetilde{\JJ}_k$ to the value $\kappa_0\approx 10^{(1+\sqrt{n_s})}$ remains the condition number of $\mathbf{H}$ in the range: $10^4\leq\kappa(\mathbf{H})\leq 10^9$. Hence, in this situation, it is feasible to perform LU factorization. 
 We see in Fig. \ref{fig:2004bp_model} that wavelengths of the wavefield have the expected behavior according to wave speed. In this test we solved \eqref{eq:helmholtz_ABC2} with ABC of second order.  Table  \ref{table:2004bp} resumes the computational complexity of the method by showing the execution times for key routines depending of frequency values and number of nodes.

\begin{table}[ht]
	\begin{center}
		\begin{tabular}{ccc|ccc}
			\rowcolor{gray!20} $\omega/2\pi$  & Nodes ($N$)  &  $\kappa(\mathbf{H})$  &Time for $\mathbf{H}_{\Omega}$ (s)  & Time for $\mathbf{H}_{\partial\Omega}$  (s) & Time for LU (s)  \\
			\hline\hline 
		2 & 38488  &  9.01e+04 & 20.59 & 4.72 & 0.99 \\
		4 & 153223 & 3.43e+05 & 74.27 & 10.52 & 5.42\\	
		6 & 345389 & 9.93e+05 & 173.78 & 19.75 & 16.04\\
		8 & 613529 & 7.22e+06 & 192.71 & 33.50 & 35.76\\	
			\hline
		\end{tabular}
		\\
	\end{center} 
	\caption{Results in computing solutions corresponding to the 2004 BP velocity model.} \label{table:2004bp}		
\end{table}

\section{Conclusions and future research}
In this paper we perform local interpolation with a  shape parameter-free oscillatory RBF based on Bessel functions of the first kind to obtain a higher order RBF-FD scheme for solving Helmholtz equation. In this approach the shape parameter is substituted for the local wavenumber. However due to the local interpolation the resulting matrices are extremely ill-conditioned even for stencils with a low number of nodes. We overcome this issue with a regularization method that introduce small perturbations in  the diagonal o the matrix. In some tests we have achieved convergence rates of third and sixth order, this is a performance in accordance with the state of the art methods for Helmholtz problems.

Among the pending problems and future research that may improve our scheme it is very important:

\begin{itemize}
    \item To  explore, in an  analytic way, the behavior of  the numerical solutions of  Helmholtz equation with discontinuous and piecewise constant coefficients.

    \item To improve the choice of regularization parameters.

    \item To exploit the meshless features of the BRBF-FD method in 3D Helmholtz problems.
    
\end{itemize}

\section*{Acknowledgments}
This work is supported by Colombian Oil Company ECOPETROL and COLCIENCIAS as a part of the research project grant No. 0266-2013.
    
    \section*{References}
	\bibliographystyle{apalike}
	\bibliography{references.bib}

	
	
	
	
	

\end{document}

\documentclass[12pt,preprint]{elsarticle}

\usepackage{amssymb}



\usepackage[utf8]{inputenc}
\usepackage{graphicx}
\usepackage{dcolumn}
\usepackage{dsfont}
\usepackage{amsthm}
\usepackage{epstopdf}
\usepackage{textcomp}
\usepackage{colortbl}
\usepackage{float}
\usepackage{bm}
\usepackage[format=plain,font={small},justification=centerlast]{caption}
\usepackage{subcaption}
\usepackage{color}
\usepackage{tabularx,ragged2e,booktabs}
\usepackage{MnSymbol}
\usepackage{wasysym}
\usepackage{anysize}
\usepackage{epstopdf}
\epstopdfsetup{update}
\usepackage{verbatim}
\usepackage{tikz}
\usepackage{pgfplots}
\usepackage{mathtools}
\usepackage{natbib}

\DeclareGraphicsExtensions{.pdf,.png,.jpg,.pdf}

\makeatletter
\newenvironment{figurehere}
{\def\@captype{figure}}
{}
\makeatother

\biboptions{comma,square}


\def\im{\mathop{\rm \od{\iota}}\nolimits}
\newcommand{\BR}{{\mathbb{R}}}
\newcommand{\p}{\partial}

\newcommand{\ts}[1]{\textstyle #1}
\newcommand{\bn}[1]{\mbox{\boldmath $#1$}}
\newcommand{\bc}{\begin{center}}
	\newcommand{\ec}{\end{center}}
\newcommand{\be}{\begin{equation}}
\newcommand{\ee}{\end{equation}}
\newcommand{\bea}{\begin{eqnarray}}
\newcommand{\eea}{\end{eqnarray}}
\newcommand{\ba}{\begin{array}}
	\newcommand{\ea}{\end{array}}

\newcommand{\JJ}{\mathbf{J}}
\newcommand{\II}{\mathbf{I}}

\newcommand{\N}{\mathds{N}}
\newcommand{\x}{\mathbf{x}}
\newcommand{\n}{\mathbf{n}}
\newcommand{\e}{\mathbf{e}}
\newcommand{\y}{\mathbf{y}}
\newcommand{\ds}{\displaystyle}
\renewcommand{\a}{\mathbf{a}}
\renewcommand{\r}{\mathbf{r}}
\newcommand{\h}{\mathbf{h}}
\renewcommand{\v}{\mathbf{v}}
\renewcommand{\u}{\mathbf{u}}
\newcommand{\bO}{\mathcal{O}}
\newcommand{\Pcal}{\mathcal{P}}

\newcommand{\Z}{\mathds{Z}}
\newcommand{\Q}{\mathds{Q}}
\newcommand{\R}{\mathds{R}}
\newcommand{\C}{\mathds{C}}
\renewcommand{\P}{\mathds{P}}
\newcommand{\phiep}{\varphi_{\epsilon}}
\newcommand{\mt}{\mathcal{T}}
\newcommand{\mb}{\mathcal{B}}
\newcommand{\ml}{\mathcal{L}}
\newcommand{\ve}{\varepsilon}
\newcommand{\om}{\Omega}
\newcommand{\td}{\tilde}
\newcommand{\al}{\alpha}
\newcommand{\vp}{\varphi}

\newtheorem{thm}{Theorem}
\newtheorem{lemma}{Lemma}
\newtheorem{definition}{Definition}
\newtheorem{rmrk}{Remark}
\newtheorem{pro}{Proposition}

\journal{---}

\begin{document}
	
	\begin{frontmatter}
		
		
		\title{Radial basis function-generated finite differences with Bessel weights for the 2$D$ Helmholtz equation}
		

		
		\author{Mauricio A. Londoño}
		\ead{alejandro.londono@udea.edu.co}
		\author{Hebert Montegranario}
		\ead{hebert.montegranario@udea.edu.co}
		
		\address{Instituto de Matemáticas \\
			Universidad de Antioquia\\
			Calle 67 53-108, Medell\'in, Colombia}
		
		\begin{abstract}
			In this paper we obtain approximated numerical solutions for the 2D Helmholtz equation using a radial basis function-generated finite difference scheme, where  weights are calculated by taking the oscillatory radial basis function given in terms of Bessel functions of first kind. The problem of local interpolation to obtain weights is ill-conditioned and we overcome this difficulty by means of regularization of the interpolation matrix by perturbing the diagonal, where the condition number of the perturbed matrix is controlled according a prescript value. We perform different tests from which observe convergence, algorithm complexity, besides we see that the pollution-effects are mitigated. 
		\end{abstract}
		
		\begin{keyword}
			RBF-FD\sep Helmholtz equation \sep shape parameter \sep pollution effect \sep Oscillatory RBF.
			
			
		\end{keyword}
		
	\end{frontmatter}
	
	
	\section{Introduction}
    
    In this work we consider the 2$D$ Helmholtz equation
    
    \begin{equation}\label{eq:helmholtz1}
    \left\{
    \begin{array}{rcll}
    -\Delta u(\x)-\omega^2c(\x)^{-2}u(\x)&=&f(\x), &\mbox{ in } \Omega\\
    b\frac{\partial}{\partial\n}u(\x)+i\omega c(\x)^{-1}\mathcal{B}u(\x)&=&g(\x), & \mbox{ on } \Gamma=\partial\Omega
    \end{array}
    \right.
    \end{equation}
where $\omega$ is the angular frequency, $c(\x)>0$ is the sound speed of the continuous media, $f(\x)$ is the source term, $\n$ is unitary normal vector to the boundary $\Gamma$, $b$ takes values zero or one, $\mathcal{B}$ is a certain linear operator and $g(\x)$ is certain exact data on $\Gamma$.

    \section{Oscillatory RBF}

Trefftz methods for the Helmholtz problems are schemes of type finite elements where test and trial functions are local solutions of the differential equation to solve. In \cite{trefftz_paper} it can be seen a wide survey about. Inspired by Trefftz methods, in this section we work with a family of oscillatory RBF which is solution of the homogeneous Helmholtz equation. Besides,  given the oscillatory behavior of solutions of Helmholtz equation, it makes sense to consider a such family, whose members are given in terms of Bessel functions of the first kind.

Now, we are going to describe briefly the class of oscillatory radial basis functions (RBF), which are given by
\begin{equation}\label{eq:bessel_RBF}
\varphi_k^{(d)}(r)=\frac{J_{d/2-1}(kr)}{(kr)^{d/2-1}}, \ \ \ d=1,2,\ldots,
\end{equation}
 whose detailed study is presented in \cite{FORNBERG2006_oscillatory_RBF}. Here $J_{\alpha}(r)$ is denoting the Bessel function of the first kind and order $\alpha$. Two remarkable properties of these oscillatory RBF are:

\begin{itemize}

	\item the non-singularity of the interpolation matrix for arbitrarily scattered data in up to $d$ dimensions, when $d>1$,

	\item and that the Laplace eigenvalue problem $\Delta \varphi + k^2\varphi =0$ has as bounded solutions, at the origin, the functions given in \eqref{eq:bessel_RBF}, thus any interpolant of the form

	\begin{equation}\label{eq:oscilla_interpolant}
	s(\x)=\sum_{j=1}^n\alpha_j\varphi^{(d)}_k(||\x-\x_j||
	\end{equation}

	will satisfy too $\Delta s+k^2 s=0$.  

\end{itemize}

In the case $d=2$ the oscillatory RBF \eqref{eq:bessel_RBF} have been little bit studied due to the function \eqref{eq:oscilla_interpolant} does not have local maximum at points where this is negative, which put a restriction to be used for general 2D interpolation. But in this work a such feature becomes a strength, since that we are just interpolate solutions of Helmholtz problems, which locally can be seen as plane waves satisfying the homogeneous Helmholtz equation. In early works, as in \cite{LIN2012_oscillatory_radial_Basis_function_helmholtz}, oscillatory RBF based on Bessel functions have been employed to solve the 2D Helmholtz equation with constant wavenumber within the approach of global collocation method.
They presented the RBF  
\begin{equation}\label{eq:RBFlin2012}
\phi_{C,k}(r)=J_0(k\sqrt{r^2+C^2}),
\end{equation}
which has two shape parameters but fixed $k$ to the wavenumber. The ill-conditioning of the interpolation matrix arises from \eqref{eq:RBFlin2012} is overcome by way of a regularized singular value decomposition method. 

For our interest, the Helmholtz problem 2D with large wavenumber, we take the special case $d=2$. So we work with the oscillatory Bessel-RBF
\begin{equation}\label{eq:bessel_RBF2D}
\phi_{k}(r)=J_0(kr),
\end{equation}
where in the approach RBF-FD the shape parameter $k$ will be evaluated at the wave number $k(\x)=\omega/c(\x)$ corresponding to the center of the stencil. As is well known, to assemble the sparse matrix, which discretizes the Helmholtz problem is necessary to solve a small linear equation system at each node.  As we will see, interpolation matrices are ill-conditioned and we will overcome this issue with the method of diagonal increments (MDI) \cite{method_of_diagonal_increments}, \cite{SARRA2014_regularized_positive_definite} adding to the diagonal entries a small regularization parameter $\beta>0$, thus we solve, instead of the linear system  $\mathbf{b}=\mathbf{A}\mathbf{y}$,  the equation
\begin{equation}
\mathbf{b}=(\mathbf{A}+\beta\mathbf{I})\widetilde{\mathbf{y}},
\end{equation}
where $\II$ is the identity matrix. We will give a brief explanation about MDI in the section \ref{sec:MDI}, where it shows that the matrix $\widetilde{\mathbf{A}}=\mathbf{A}+\beta\II$ is better conditioned than $\mathbf{A}$ and  $\mathbf{y}\approx\widetilde{\mathbf{y}}$. Before we describe the goodness of discretizing Helmholtz problems with Bessel-RFB.


\section{Bessel-RBF-FD}

Suppose that $u$ is a solution of the Helmholtz equation 
$\Delta u(\x)+k(\x)^2u(\x)=0$, for $\x\in\Omega$ and for $\x\in\partial\Omega$ it satisfies certain boundary condition. If $X=\{\x_i\}_{i=1}^N\subset\Omega\cup\partial\Omega$ is a set of nodes, for $\x_i\in X\cap\Omega$ we take a stencil $S_i=\{\x_j^i\}_{j=1}^{n_i}\subset X$ based on $\x_i$, with $\x_1^i=\x_i$. For $\x\in\Omega_i=\mbox{ConvexHull}(S_i)$ we define, with $k_i=k(\x_i)$, the interpolant 
\begin{equation}\label{eq:interpolant_bessel}
\widetilde{u}(\x)=\sum_{j=1}^{n_i}\alpha_j^i J_0(k_i\|\x-\x_j^i\|).
\end{equation}
With the local interpolation matrix, $\mathbf{J}_{k_i}=(J_0(k_i\|\x_l^i-\x_j^i\|))_{1\leq l,j\leq n_i}$, which is positive definite, and forcing the condition $\widetilde {u}(\x_l^i)=u(\x_l^i)$, then from \eqref{eq:interpolant_bessel} we have the linear equation
 \begin{equation}\label{eq:system_interpolation_ill}
 U_i=\mathbf{J}_{k_i}\bm{\alpha}_i,
 \end{equation}
where $U_i=\left(\begin{array}{cccc}
u(\x_1^i)&u(\x_2^i)&\cdots&u(\x_{n_i}^i)
\end{array}\right)^t$ and $\bm{\alpha}_i=\left(\begin{array}{cccc}
\alpha_1^i&\alpha_2^i&\cdots&\alpha_{n_i}^i
\end{array}\right)^t$.
In view that $\phi_k$, defined in \eqref{eq:bessel_RBF2D}, satisfies the homogeneous Helmholtz equation, then
\begin{eqnarray}
\Delta\widetilde{u}(\x)|_{\x=\x_i}&=&\sum_{j=1}^{n_i}\alpha_j^i \Delta J_0(k_i\|\x-\x_j^i\|)|_{\x=\x_i}\\
&=&-k_i^2\sum_{j=1}^n\alpha_j^i J_0(k_i\|\x_i-\x_j^i\|)\\
&=&-k_i^2\widetilde{u}(\x_i).
\end{eqnarray}
Hence the interpolant \eqref{eq:interpolant_bessel} satisfies the homogeneous Helmholtz equation. On the other hand, applying  $\Delta_{S_i,k_i}$ to the solution $u$, we have
\begin{eqnarray}
\Delta_{S_i,k_i}u(\x_i)&=&\Delta \mathbf{J}_{k_i}\mathbf{J}_{k_i}^{-1}U_i\\
&=&-k^2_i\mathbf{e}_1U_i\\
&=&-k_i^2u(\x_i)
\end{eqnarray}
Note that $\Delta_{S_i,k_i}u(\x_i)-\Delta_{S_i,k_i}\widetilde{u}(\x_i)=-k_i^2(u(\x_i)-\widetilde{u}(\x_i))$, thus, for solutions of the homogeneous Helmholtz problem the local truncation error for the Laplace operator has a theoretical error depending of wavenumber at $\x_i$ and of the error of the local interpolant. The error of the approximated  solutions is produced by the interpolant \eqref{eq:interpolant_bessel} and by the ill-conditioning of the matrix $\JJ_{k_i}$, in solving the linear system \eqref{eq:system_interpolation_ill}. Next we will deal with solutions fo these systems.

\section{Method of diagonal increments (MDI)}\label{sec:MDI}

The interpolation matrix $\JJ_k$ is ill-conditioned, especially for certain node distributions. In literature there are several methods to overcome the ill-conditioning when the shape parameter is small \cite{stable_general_bessel_rbf}, ever in our case we are taking the shape parameter as the wavenumber $k$, which can be large. So we have chosen the MDI. For our case $\JJ_k$ will be considered  ill-conditioned when the condition number, $\kappa(\JJ_k)$, satisfies  $\kappa(\JJ_k)>10^{15}$, which hinders that the solution $\bm{\alpha}$ of $U_i=\mathbf{J}_k\bm{\alpha}$ be accurately calculated, in double precision, 	through Cholesky facorization. When $\JJ_k$ is ill-conditioned we solve the better conditioned problem $U_i=(\mathbf{J}_k+\beta\mathbf{I})\widetilde{\bm{\alpha}}$ instead, where $\II$ is  the identity matrix and $\beta$  a small positive real number.  Next, we will give some important aspects about the spectrum of $\JJ_k+\beta\II$ in order to give an estimate of the error $\|\bm{\alpha}-\widetilde{\bm{\alpha}}\|_2$.

\begin{rmrk}\label{rmrk:convergence_neumann_series_jk}

	$\JJ_k$ is positive definite, thus its spectrum is real and positive. If $\{\lambda_m\}_{m=1}^n$ is its spectrum,  with $\lambda_1\geq\lambda_2\geq\cdots\geq\lambda_n$, then $\{\lambda_m+\beta\}_{m=1}^n$ is the specturm of $\widetilde{\JJ}_k=\JJ_k+\beta\II$ and $\{\frac{\beta}{\lambda_m+\beta}\}_{m=1}^n$ is the spectrum of $\beta\widetilde{\JJ}_k^{-1}$, hence its spectral norm is $\|\beta\widetilde{\JJ}_k^{-1}\|_2=\frac{\beta}{\lambda_n+\beta}$. The  above implies that the Neumann series  $\sum_{m=0}^{\infty}(\beta\widetilde{\JJ}_k^{-1})^m$ converges and the equality 
	\begin{equation}
	(\II-\beta\widetilde{\JJ}_k^{-1})^{-1}=\sum_{m=0}^{\infty}(\beta\widetilde{\JJ}_k^{-1})^m
	\end{equation} holds.

\end{rmrk}

\begin{rmrk}\label{rmrk:better_coditioning_wideJk}

	If $\{\lambda_m\}_{m=1}^n$ is the spectrum  of $\JJ_k$,  then the condition number of $\JJ_k$ is given by $\kappa(\JJ_k)=\frac{\lambda_1}{\lambda_n}$ and $\kappa(\widetilde{\JJ}_k)=\frac{\lambda_1+\beta}{\lambda_n+\beta}$, which implies that 
	\begin{equation}
	\kappa(\widetilde{\JJ}_k)<\kappa({\JJ}_k)
	\end{equation}
	whit this, the matrix $\widetilde{\JJ}_k$ is better conditioned than $\JJ_k$.

\end{rmrk}

\begin{rmrk}\label{rmrk:aprrox_jk_by_neumann_series}

Given that $\widetilde{\JJ}_k=\JJ_k+\beta\II$, then $\JJ_k^{-1}=\widetilde{\JJ}_k^{-1}(\II-\beta\widetilde{\JJ}_k^{-1})^{-1}$ and by using the Neumann series we have
\begin{eqnarray}
\JJ_k^{-1}&=&\widetilde{\JJ}_k^{-1}\sum_{m=0}^{\infty}(\beta\widetilde{\JJ}_k^{-1})^m\\
&=&\frac{1}{\beta}\sum_{m=1}^{\infty}(\beta\widetilde{\JJ}_k^{-1})^m
\end{eqnarray}

\end{rmrk}

\begin{rmrk}

	If $\bm{\alpha}$ is the true solution of the equation $\JJ_k\bm{\alpha}=U$ and $\widetilde{\bm{\alpha}}$ is the solution for the equation $\widetilde{\JJ}_k\widetilde{\bm{\alpha}}=U$ then, from the remark \ref{rmrk:aprrox_jk_by_neumann_series}, \begin{equation}\label{eq:alpha_series}
	\bm{\alpha}=\frac{1}{\beta}\sum_{m=1}^{\infty}(\beta\widetilde{\JJ}_k^{-1})^mU
	\end{equation} and
	\begin{equation}
	\bm{\alpha}=\sum_{m=0}^{\infty}(\beta\widetilde{\JJ}_k^{-1})^m\widetilde{\bm{\alpha}}.
	\end{equation}
\end{rmrk}

If we truncate the series in \eqref{eq:alpha_series} up to order $M$ we obtain an approximation of the true   solution $\bm{\alpha}$, we denote it by  
\begin{equation}\label{eq:approx_alpha_barM}
\widetilde{\bm{\alpha}}_M=\frac{1}{\beta}\sum_{m=1}^{M}(\beta\widetilde{\JJ}_k^{-1})^mU.
\end{equation}
With this, the error of the approximation $\widetilde{\bm{\alpha}}_M$ can be bounded by using the formula
\begin{eqnarray}
\bm{\alpha}-\widetilde{\bm{\alpha}}_M&=&\frac{1}{\beta}\sum_{m=M+1}^{\infty}(\beta\widetilde{\JJ}_k^{-1})^mU\\
&=&\frac{1}{\beta}(\II-\beta\widetilde{\JJ}_k^{-1})^{-1}(\beta\widetilde{\JJ}_k^{-1})^{M+1}U.
\end{eqnarray}
Taking the Euclidean norm, then we have in terms of the spectral norm,
\begin{equation}
\|\bm{\alpha}-\widetilde{\bm{\alpha}}_M\|_2\leq\| (\beta\widetilde{\JJ}_k^{-1})^{M+1}\|_2\|(\II-\beta\widetilde{\JJ}_k^{-1})^{-1} \|_2\|U\|_2.
\end{equation}
From the remark \ref{rmrk:convergence_neumann_series_jk} we can conclude that
\begin{equation}\label{eq:error_bound_MDI}
\|\bm{\alpha}-\widetilde{\bm{\alpha}}_M\|_2\leq\frac{\beta}{\lambda_n}\left(\frac{\beta}{\lambda_n+\beta}\right)^{M}\|U\|_2.
\end{equation}
In particular,
\begin{equation}\label{eq:error_bound_R}
	\|\bm{\alpha}-\widetilde{\bm{\alpha}}\|_2\leq\frac{\beta}{\lambda_n}\|U\|_2.
\end{equation}

An iterative proceeder to compute \eqref{eq:approx_alpha_barM}, with the better conditioned matrix $\beta\widetilde{\JJ}_k$, can be obtained just by noting that, with $\widetilde{\bm{\alpha}}=\widetilde{\JJ}_k^{-1}U$,
  \begin{eqnarray}\label{eq:approx_alpha_barM_iter}
  \widetilde{\bm{\alpha}}_M&=&\frac{1}{\beta}\sum_{m=1}^{M}(\beta\widetilde{\JJ}_k^{-1})^mU\\
  &=&\sum_{m=1}^{M}(\beta\widetilde{\JJ}_k^{-1})^{m-1}\widetilde{\bm{\alpha}}\\
  &=&\widetilde{\bm{\alpha}}+\beta\widetilde{\JJ}_k^{-1}\left(\widetilde{\bm{\alpha}}+\beta\widetilde{\JJ}_k^{-1}\left(\widetilde{\bm{\alpha}}+\cdots\right)\right).
  \end{eqnarray}
  Hence we can compute $\widetilde{\bm{\alpha}}_M$  as: 
  \begin{eqnarray}
  	\label{eq:iter_first_step}\widetilde{\bm{\alpha}}_0&=&\widetilde{\JJ}_k^{-1}U\\
  	 \label{eq:iter_m_step}\widetilde{\bm{\alpha}}_{m}&=&\widetilde{\bm{\alpha}}_0+\beta\widetilde{\JJ}^{-1}_k\widetilde{\bm{\alpha}}_{m-1}, \ \ \ \mbox{\ \  for } m=1,2,\ldots,M
  \end{eqnarray}
Since $0<\frac{\beta}{\lambda_n+\beta}<1$ the convergence order is controlled theoretically, however when $\lambda_n$ is near to the machine epsilon (e.g., in double precision the machine epsilon is approximately 2.22e-16) the ratio $\frac{\beta}{\lambda_n+\beta}$ is very close to 1 and in this case the convergence may be too slow, but with a small $M$ (up to $M=5$ it works) is enough to improve the error in \eqref{eq:error_bound_R}.

\section{Numerical tests}
For this section we have tried to reproduce examples and tests presented in \cite{Fang2017} due to the solver presented by authors is recent and very accurate, also is they show a detailed local description of solutions of Helmholtz equation in terms of its phases and ray direction local.
\subsection{Local truncation error}

For small stencils we take $10^7 \leq\kappa_0\leq10^{14}$, with $\kappa_0=10^{7+\sqrt{n}}$ where $n$ is the size of the stencil,  and we take $\beta=\frac{\lambda_1-\kappa_0  \lambda_n}{\kappa_0 -1}$ ensuring, from remark \ref{rmrk:better_coditioning_wideJk}, that $\kappa(\widetilde{\JJ}_k)\approx\kappa_0$, which is an adequate condition number to work in double precision. 

\begin{figure}
	\begin{center}	
		\begin{tabular}{|c|c|c|c|}
			\hline
			\hline
			\includegraphics[scale=0.5]{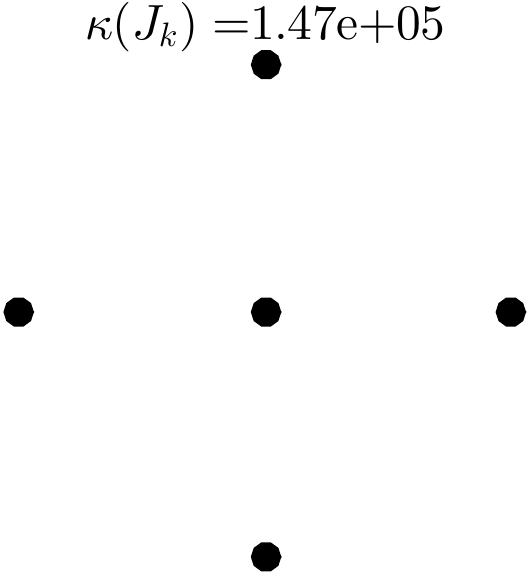}&\includegraphics[scale=0.5]{./figures/s9.pdf} &\includegraphics[scale=0.5]{./figures/s13.pdf} &\includegraphics[scale=0.5]{./figures/s25.pdf}\\
			\hline
			\hline
			\includegraphics[scale=0.5]{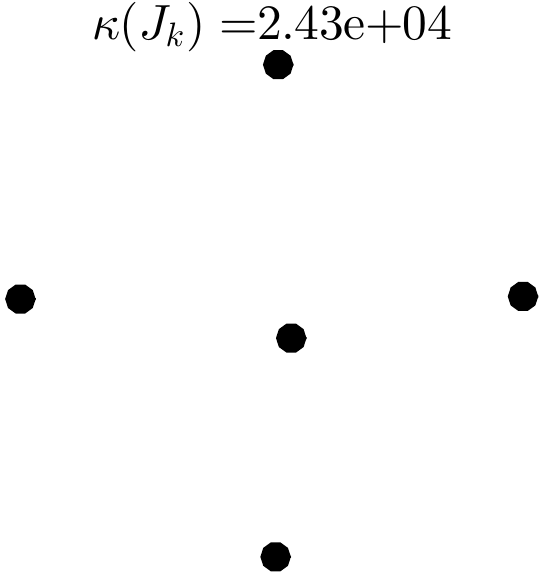}&\includegraphics[scale=0.5]{./figures/s9p.pdf} &\includegraphics[scale=0.5]{./figures/s13p.pdf} &\includegraphics[scale=0.5]{./figures/s25p.pdf}\\
			\hline
		\end{tabular}
	\end{center}
	\caption{Plots of some small stencils with the respective approximated condition number of $\JJ_k$, with $k=100$ and $h=\frac{2\pi}{6k}$. Top row: Stencils are taken from a regular square grid. Bottom row: perturbed position from stencils of top row.} 	\label{fig:stencils_square}
\end{figure}

\begin{figure}
	\begin{center}	
		\begin{tabular}{|c|c|c|c|c|}
			\hline
			\hline
			\includegraphics[scale=0.45]{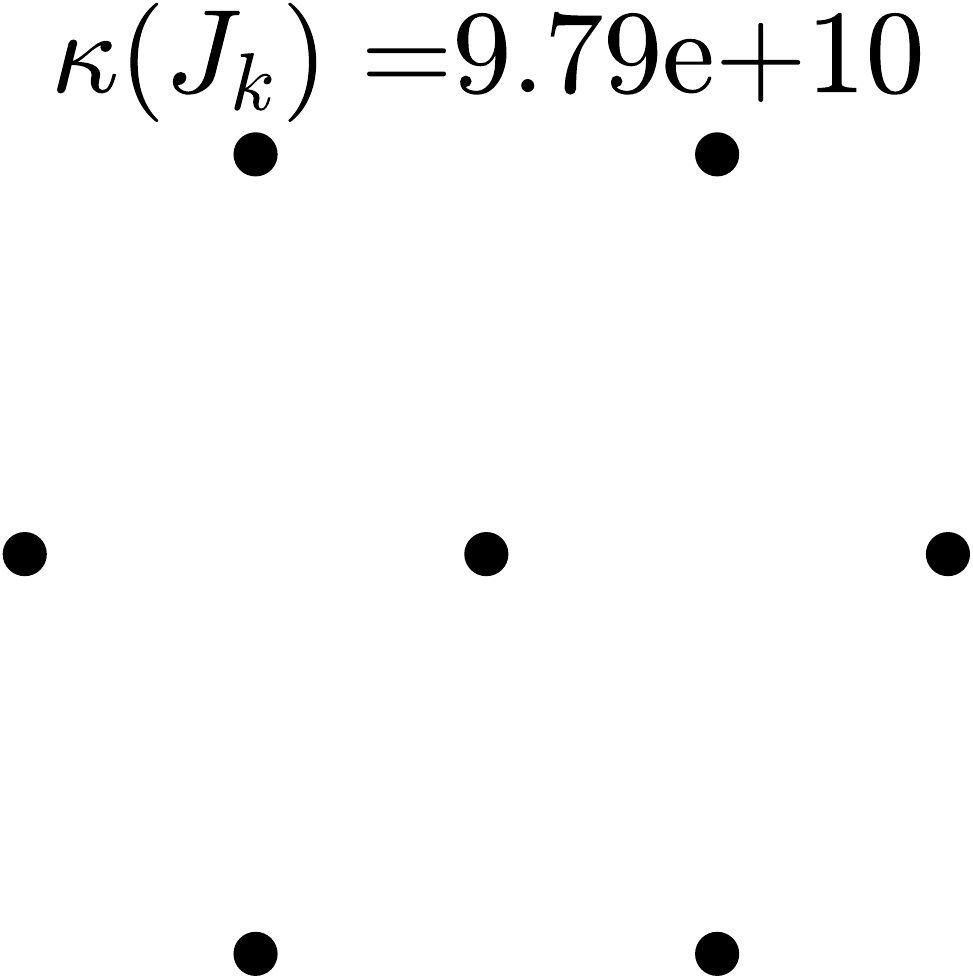}&\includegraphics[scale=0.45]{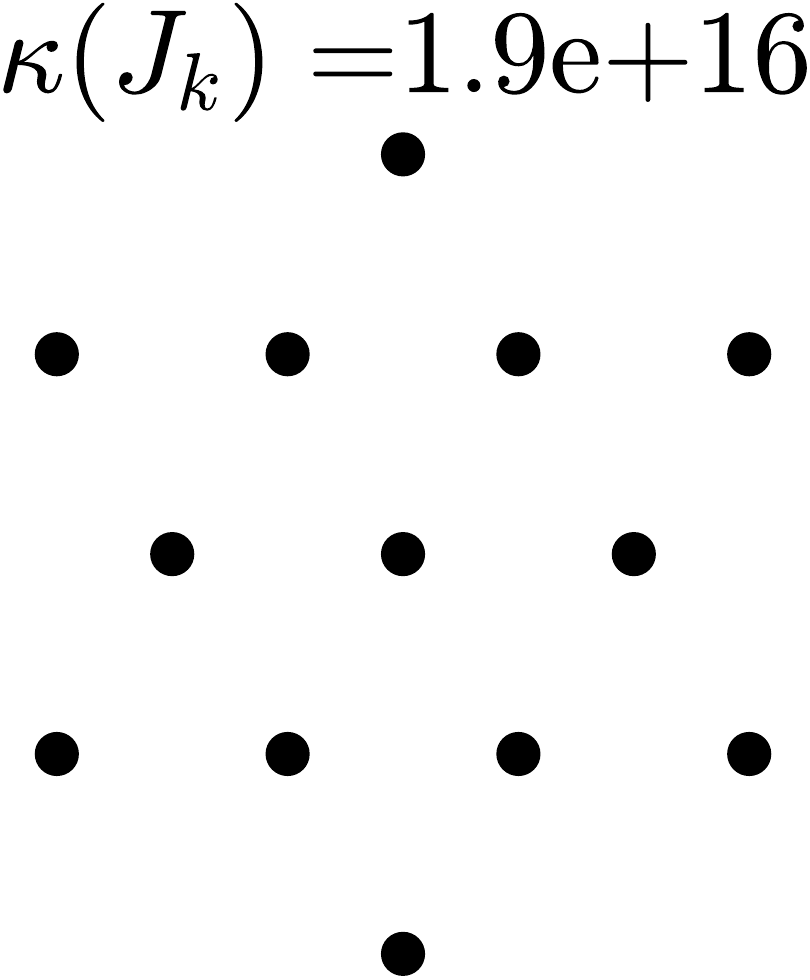} &\includegraphics[scale=0.45]{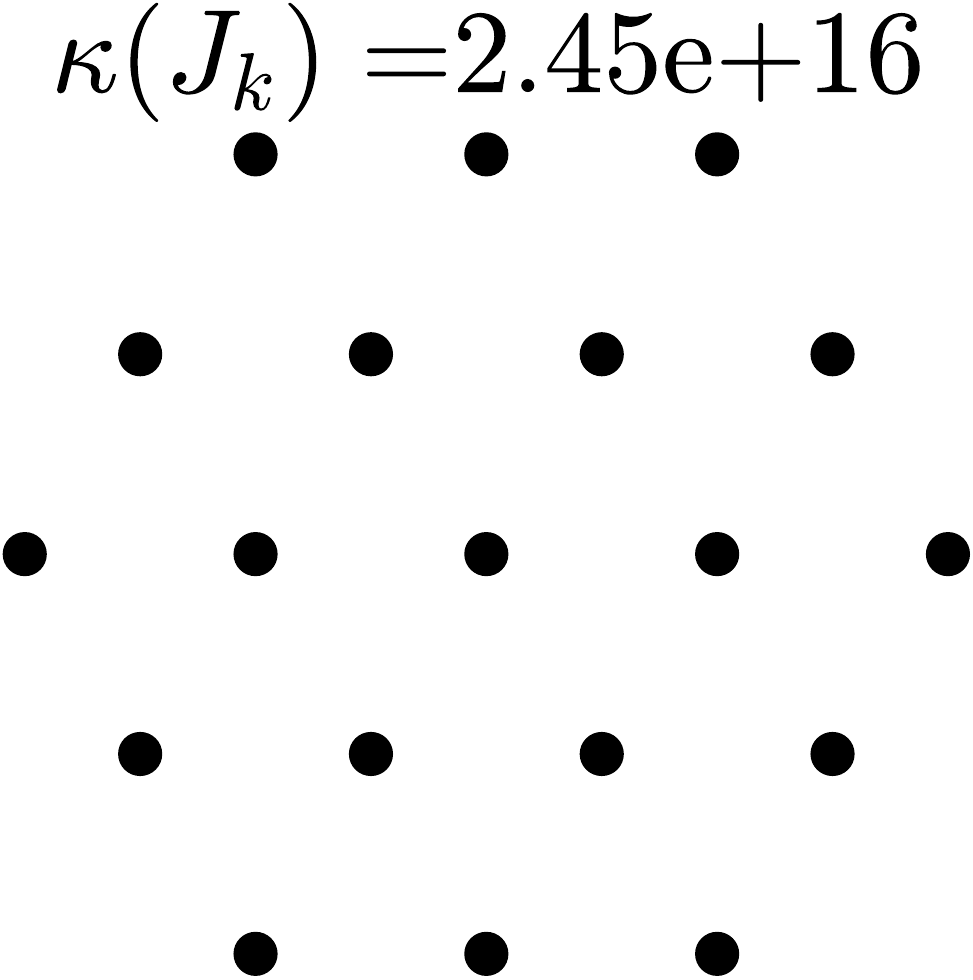} &\includegraphics[scale=0.45]{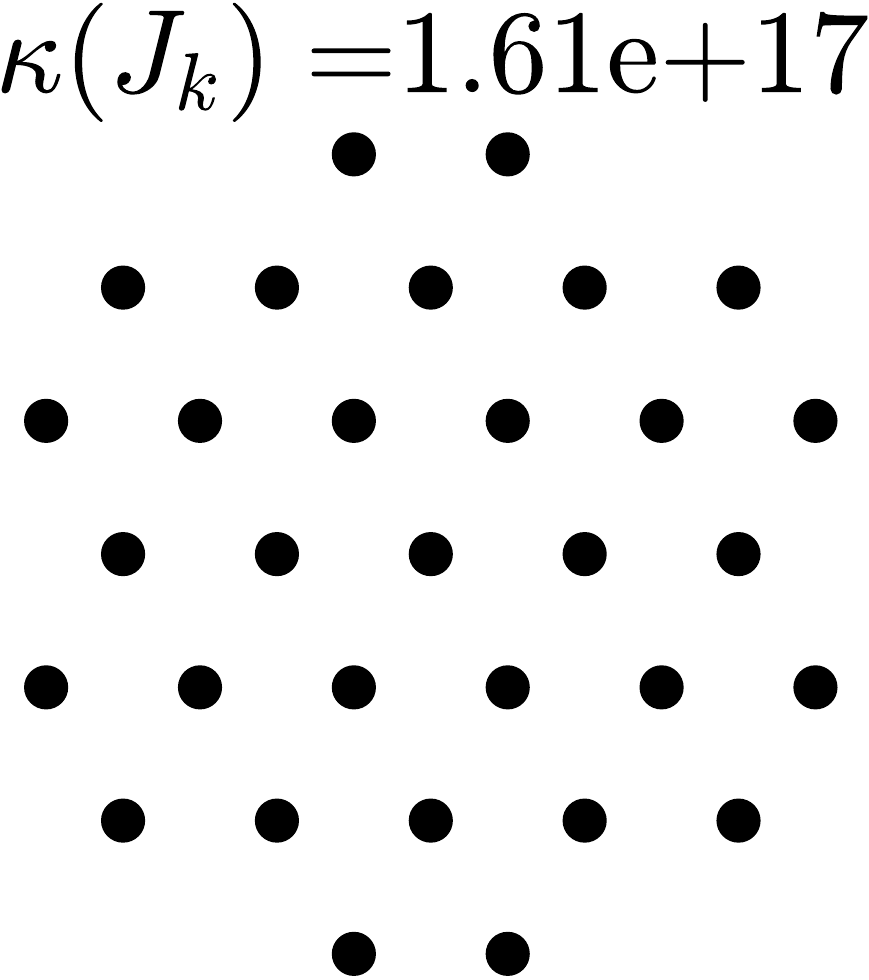}
			&\includegraphics[scale=0.45]{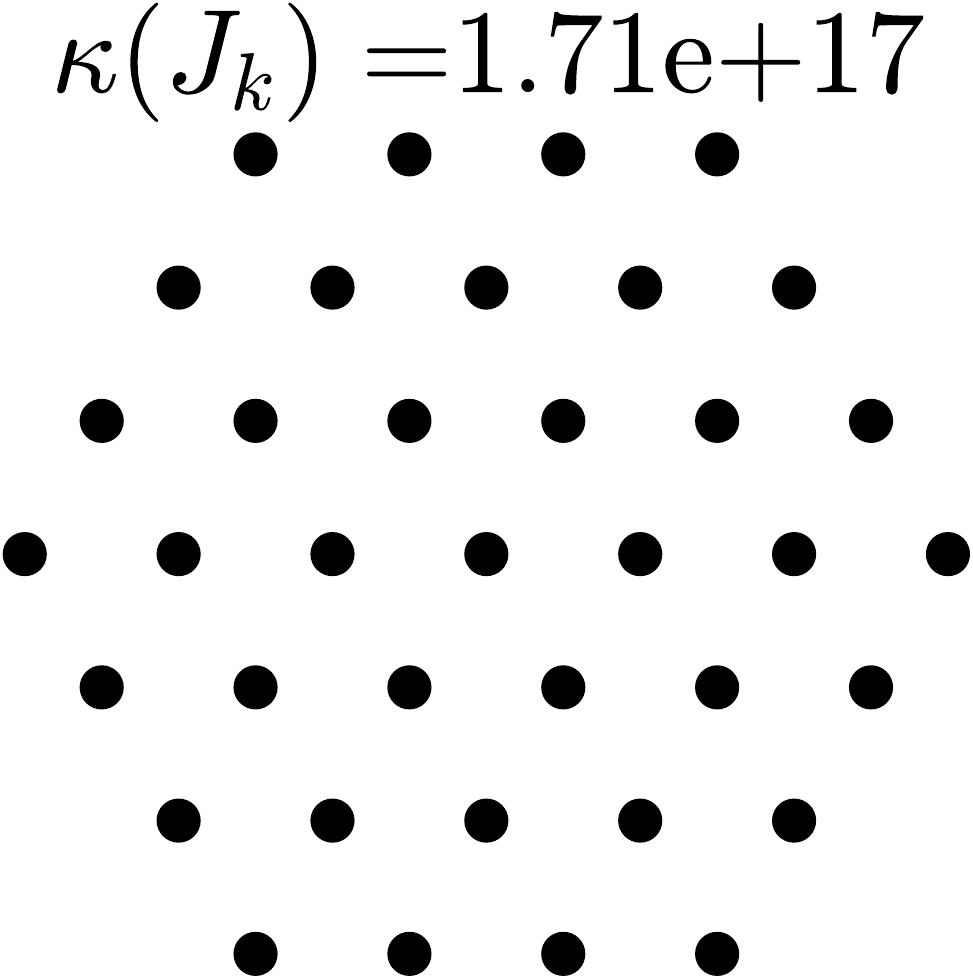}\\
			\hline
			\hline
			\includegraphics[scale=0.45]{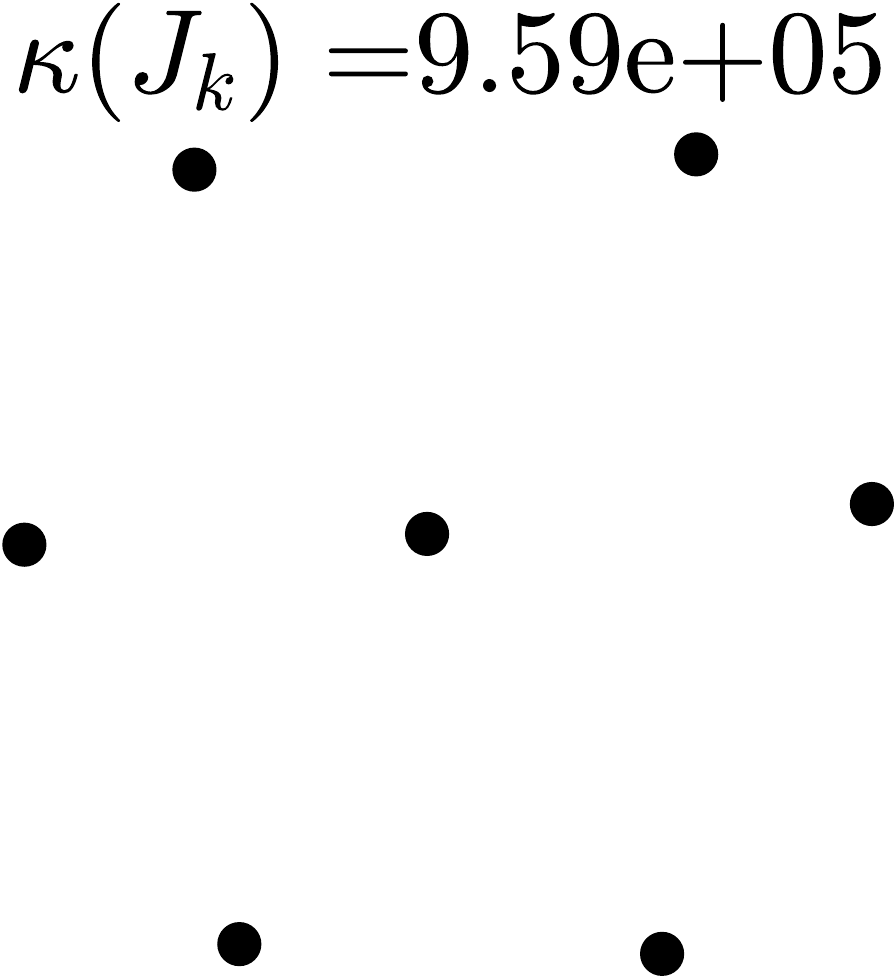}&\includegraphics[scale=0.45]{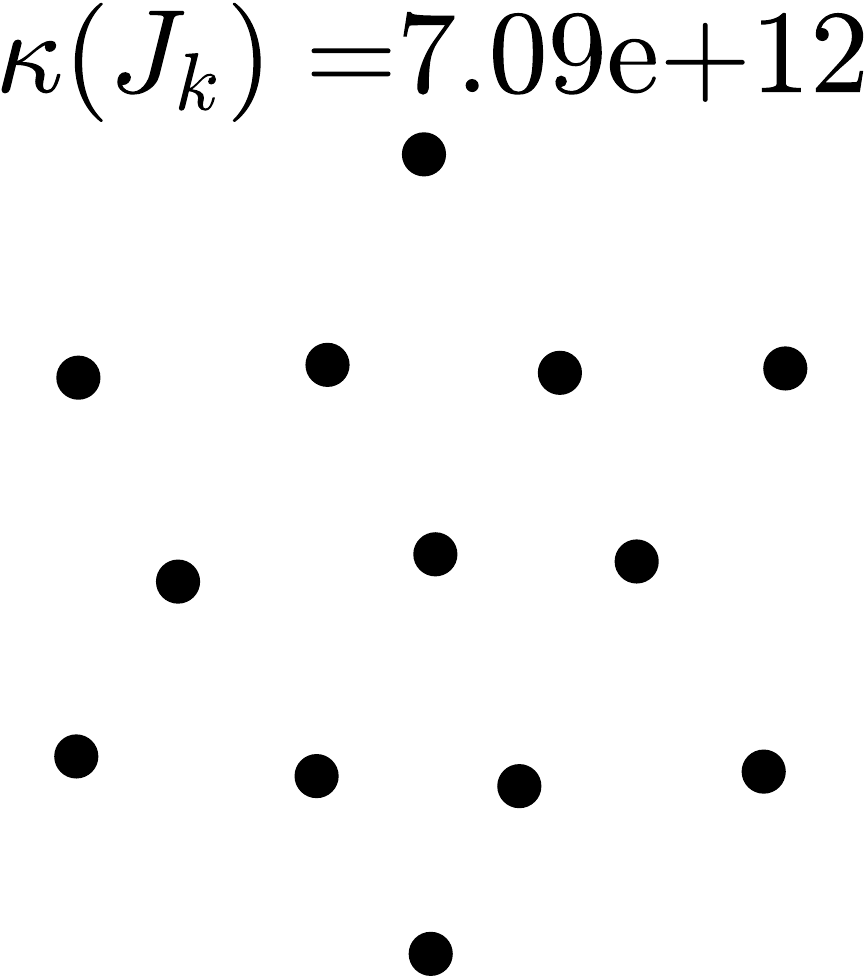} &\includegraphics[scale=0.45]{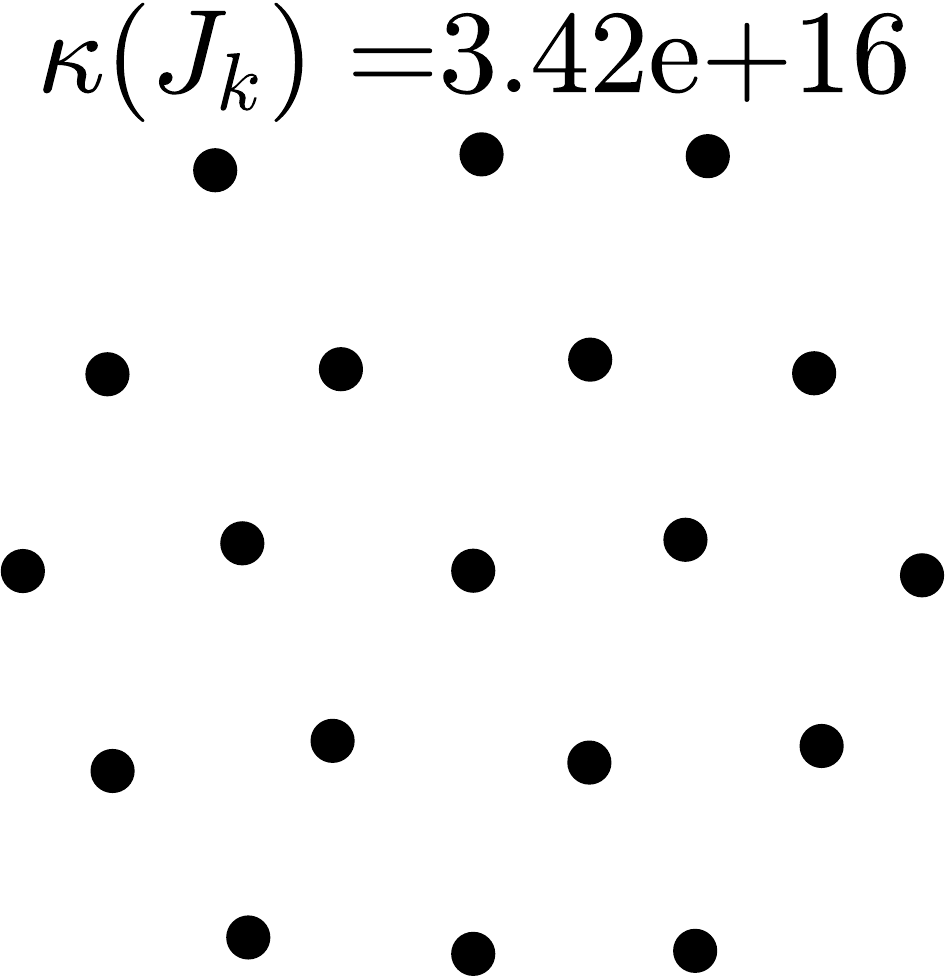} &\includegraphics[scale=0.45]{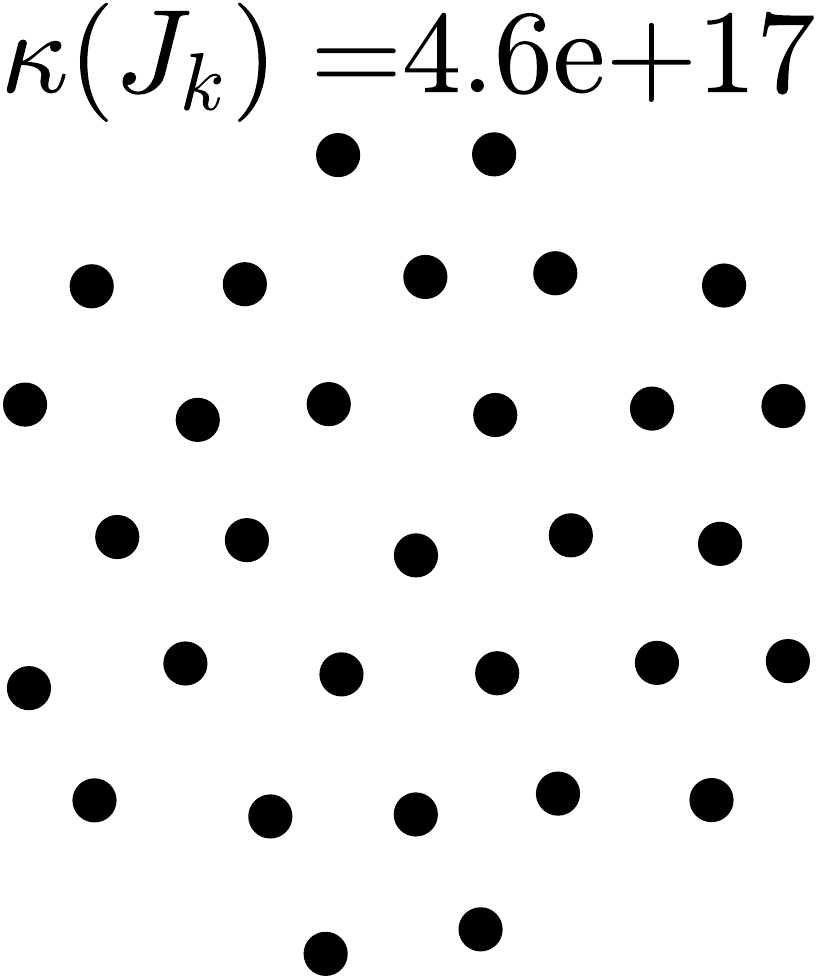}
			&\includegraphics[scale=0.45]{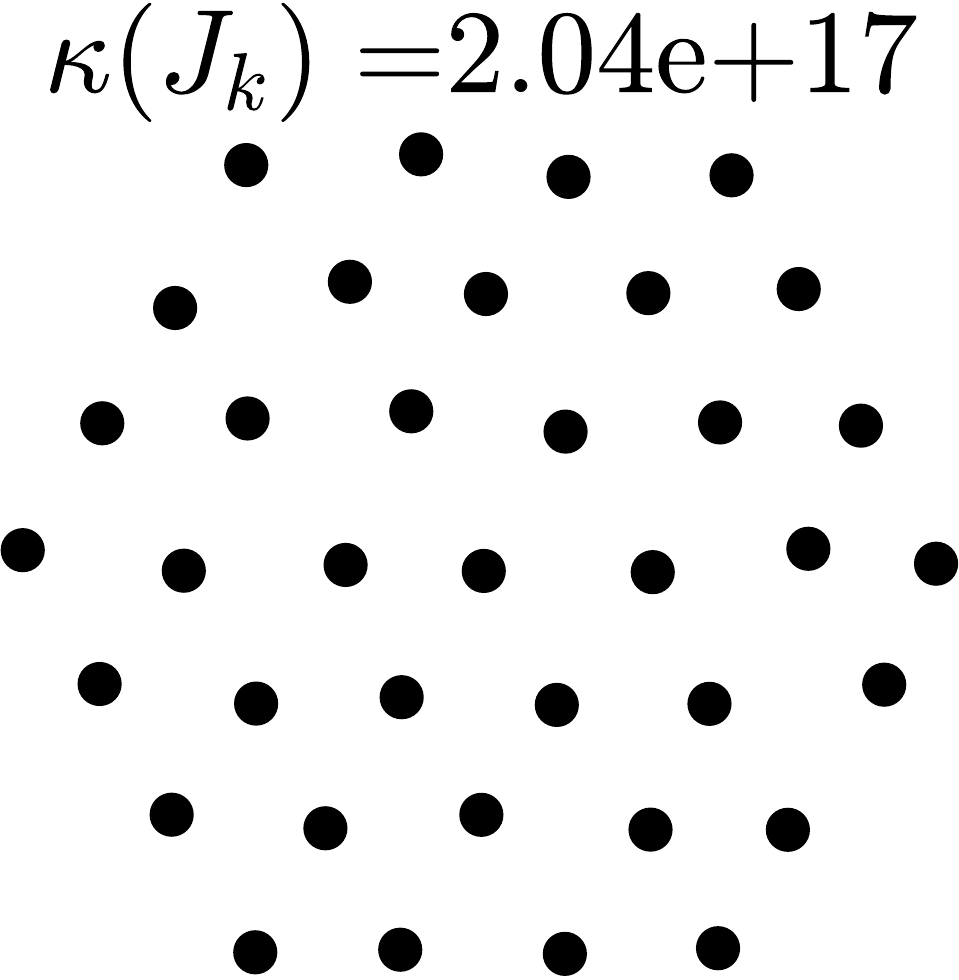}\\
			\hline
		\end{tabular}
	\end{center}
	\caption{Plots of some small stencils with the respective approximated condition number of $\JJ_k$, with $k=100$ and $h=\frac{2\pi}{6k}$. Top row: Stencils are taken from a regular hexagonal grid. Bottom row: perturbed position from stencils of top row.} \label{fig:stencils_hexagonal}
\end{figure}

We have noted empirically that the matrix $\JJ_k$ is worse conditioned for stencil with nodes put symmetrically on a regular grids, e.g. with square and hexagonal grids. See figures \ref{fig:stencils_square} and \ref{fig:stencils_hexagonal}, where we can observe that first size of severe ill-conditioning occur with symmetric stencils of $13$ nodes. However, in this case, with a small perturbation in the position of the nodes, its associated interpolation matrix $\JJ_k$ has a better condition number.

For numerical tests we consider the Helmholtz solutions

\begin{equation}\label{eq:sol_imped1}
u_1(x,y)=\sqrt{k}H^{(1)}_0(k\sqrt{(x-2)^2+(y-2)^2}),
\end{equation}
and
\begin{equation}\label{eq:sol_imped2}
\begin{array}{rcl}
u_2(x,y)&=&\sqrt{k}H^{(1)}_0(k\sqrt{(x+20)^2+(y+20)^2})+2\sqrt{k}H^{(1)}_0(k\sqrt{(x-20)^2+(y-20)^2})\\
&+&0.5\sqrt{k}H^{(1)}_0(k\sqrt{(x+20)^2+(y-20)^2})-\sqrt{k}H^{(1)}_0(k\sqrt{(x-20)^2+(y+20)^2}),
\end{array}
\end{equation}
where $H_0^{(1)}$ is the Hankel function of the first kind and $k$ is the constant wavenumber. $u_1$ corresponds to solution of the single source problem located at $\x_s=(2,2)$ whereas $u_2$ corresponds to solution of the problem with four single sources located at $\x_{s_1}=(-20,-20)$, $\x_{s_2}=(20,20)$, $\x_{s_3}=(-20,20)$ and $\x_{s_1}=(20,-20)$.
In the first test, we consider the exact solution to the Helmholtz equation \eqref{eq:sol_imped2} to get the local truncation error of the approximation $\Delta_{S,k}$, also we validate the ill-conditioning of the matrix $\JJ_k$ and the better conditioning of $\widetilde{\JJ}_k$. We compare 
to interpolate data produced with the function  $u(x,y)=\dfrac{\partial u_2}{\partial x}(x,y)$, at stencils nodes centered at $(x_c,y_c)=(0.5,0.5)$. 
With stencils $S$ as in Fig. \ref{fig:stencils_square} and Fig. \ref{fig:stencils_hexagonal} we solve the system $\JJ_k\bm{\alpha}=U$ with $U=u|_S$. In tables \ref{table:Jk_bar_square} and \ref{table:Jk_bar_hexagonal} it can be seen how is the behavior of the condition number $\kappa(\JJ_k)$ and $\kappa(\widetilde{\JJ}_k)$ for symmetric nodes distributions respectively. We have compared the relative errors of calculate approximations of $\alpha$: with direct method, the MDI and MDIIT, whit error: Error D, Error RD and Error IT as in tables \ref{table:Jk_bar_square} and \ref{table:Jk_bar_hexagonal}.
\begin{table}[ht]
	\begin{center}
		\begin{tabular}{c|cc|ccc}
			\hline
			\rowcolor{gray!20} Stencil size ($n$)  &  $\kappa(\JJ_k)$ & $\kappa(\widetilde{\JJ}_k)$ & Error D & Error RD  & Error IT \\
			\hline\hline 
	        5 & 1.47e+05 & 1.72e+09 & 3.22e-13 & 8.78e+04 & 3.22e-13\\
	        9 & 4.68e+12 & 1.00e+10 & 1.18e-08 & 1.60e-03 & 1.57e-03\\
	        13 & 6.75e+15 & 4.03e+10 & 5.67e-08 & 1.31e-03 & 1.48e-07\\
	        25 & 9.42e+16 & 1.00e+12 & 3.44e-06 & 3.34e-06 & 2.79e-06\\
	        37 & 9.59e+16 & 1.21e+13 & 1.80e-05 & 1.43e-04 & 7.26e-08\\ \hline
		\end{tabular}
		\\
	\end{center} 
	\caption{This table shows values of the condition number of the matrices $\JJ_k$ and  $\widetilde{\JJ}_k$, and local truncation errors of the approximation $\Delta_{S,k}u(\x)\approx\Delta u(\x)$. We have used stencils as in Fig. \ref{fig:stencils_square}. We compare errors $|\Delta_{S,k}u(\x)-\Delta u(\x)|$ produced by using Cholesky (Error D), by \eqref{eq:iter_first_step} (Error RD) and   \eqref{eq:iter_m_step} (Error IT) in using $5$ iterations applied to solution. The used function is $u=u_2$ as in \eqref{eq:sol_imped2}.}\label{table:Jk_bar_square}	
\end{table}

\begin{table}[ht]
	\begin{center}
		\begin{tabular}{c|cc|ccc}
			\rowcolor{gray!20} Stencil size ($n$)  &  $\kappa(\JJ_k)$ & $\kappa(\widetilde{\JJ}_k)$ & Error D & Error RD  & Error IT \\
			\hline\hline 
	   	7 & 2.06e+09 & 4.42e+09 & 3.48e-10 & 2.19e-02 & 5.06e-02\\
		13 & 6.09e+15 & 4.03e+10 & 6.26e-09 & 4.15e-06 & 9.01e-08\\
		19 & 3.57e+16 & 2.29e+11 & 1.12e-06 & 2.75e-06 & 2.77e-06\\
		31 & 2.03e+17 & 3.70e+12 & 5.47e-06 & 8.83e-06 & 7.71e-06\\
		37 & 1.18e+17 & 1.21e+13 & 3.99e-07 & 4.33e-05 & 3.40e-05\\ \hline
		\end{tabular}
		\\
	\end{center} 
		\caption{This table shows values of the condition number of the matrices $\JJ_k$ and  $\widetilde{\JJ}_k$, and local truncation errors of the approximation $\Delta_{S,k}u(\x)\approx\Delta u(\x)$. We have used stencils as in Fig. \ref{fig:stencils_hexagonal}. We compare errors $|\Delta_{S,k}u(\x)-\Delta u(\x)|$ produced by using Cholesky (Error D), by \eqref{eq:iter_first_step} (Error RD) and   \eqref{eq:iter_m_step} (Error IT) in using $5$ iterations applied to solution. The used function is $u=u_2$ as in \eqref{eq:sol_imped2}.}\label{table:Jk_bar_hexagonal}		
\end{table}

\subsection{Pollution-effect and convergence}

\subsubsection{Test 1}

In this test we calculate the approximated solution for the problem

    \begin{equation}\label{eq:helmholtz_impedance}
    \left\{
    \begin{array}{rcll}
    -\Delta u(\x)-\omega^2u(\x)&=&0, &\mbox{ in } \Omega\\
    \frac{\partial}{\partial\n}u(\x)+i\omega u(\x)&=&g(\x), & \mbox{ on } \Gamma=\partial\Omega
    \end{array}
    \right.
    \end{equation}
     with the known data $g(\x)$, $\Omega=(-0.5,0.5)\times(-0.5,0.5)$.  Results are verified with solutions $u_1$ and $u_2$ in \eqref{eq:sol_imped1} and \eqref{eq:sol_imped2}. 
     
     We see in tables \ref{table:fang1square_pollution} and \ref{table:fang1hexagon_pollution} how it behaves the error when the resolution is kept constant with $N_g=6$ nodes per wavelength. We have fixed the condition number of the local interpolation matrices, $\widetilde{\JJ}_k$, to $\kappa_0=10^{(7+\sqrt{n})}$. In both cases, for uniform square and hexagonal grids, we see that the errors are $\bO(\omega)$, but with a small positive slope and in some cases, increasing the stencil size, the error fits to $m\omega+b$ with $m$ small and negative, see table \ref{table:fang1hexagon19_pollution}. Although pollution-effects does not eliminated this is kept relatively constant, in accord to order error and does not increase dramatically.

     Results of convergence tests are summarized in tables \ref{table:fang1convergence_square} and \ref{table:fang1convergence_hexa}. For the local interpolation matrices $\widetilde{\JJ}_k$ we have chosen, empirically, $\beta$ such that the condition number  $\kappa_0=10^{(7+(hk)^{-1}+\sqrt{n})}$.

\begin{table}
	\begin{center}
		\begin{tabular}{cccc|cc}
			\rowcolor{gray!20} $\frac{\omega}{2\pi}$  &  $\frac{1}{h}$ & $N$ & $\kappa(\mathbf{H})$ & $\|u_1-\widetilde{u}_1\|_{\infty}$ & $\|u_2-\widetilde{u}_2\|_{\infty}$   \\
			\hline\hline 
	    	10 & 60 & 3721 & 1.36e+04 & 1.76e-04 & 1.72e-04  \\  
	    	20 & 120 & 14641 & 7.46e+04 & 2.54e-04  & 1.71e-04   \\
			40 & 240 & 58081 & 4.19e+05 & 5.15e-04 & 2.82e-04    \\
			80 & 480 & 231361 & 2.34e+06 & 1.04e-03 & 5.90e-04\\
			120 & 720 & 519841 & 5.25e+06 & 1.57e-03 & 8.65e-04  \\
			\hline
		\end{tabular}
		\\
	\end{center} 
	\caption{Results for approximated solutions of \eqref{eq:helmholtz_impedance}. This approximation was calculated from a square uniform grid of $\Omega\cap\partial\Omega$. For inner nodes the stencil size is $n=9$, at boundary nodes $n=15$, the number of nodes per wave length is kept constant with $N_g=6$. Solution for the system $-\mathbf{H}\widetilde{U}=G$ is obtained for LU factorization. Errors scale respect to the frequency as $\|u_1-\widetilde{u}_1\|_{\infty}\sim(1.3e-05)(\omega/2\pi)+1.7e-05$ and $\|u_2-\widetilde{u}_2\|_{\infty}\sim(6.6e-06)(\omega/2\pi)+5.8e-05$ } \label{table:fang1square_pollution}		
\end{table}

\begin{table}
	\begin{center}
		\begin{tabular}{cccc|cc}
			\rowcolor{gray!20} $\frac{\omega}{2\pi}$  &  $\frac{1}{h}$ & $N$ & $\kappa(\mathbf{H})$ & $\|u_1-\widetilde{u}_1\|_{\infty}$ & $\|u_2-\widetilde{u}_2\|_{\infty}$   \\
			\hline\hline 
			10 & 60  & 4237   & 1.96e+04 & 3.95e-05 & 3.30e-05  \\  
	    	20 & 120 & 16752  & 1.09e+05 & 4.49e-05 & 3.35e-05  \\
	     	40 & 240 & 66861  & 5.25e+05 & 9.30e-05 & 3.59e-05   \\
			80 & 480 & 266680 & 3.46e+06 & 6.69e-05 & 5.62e-05\\
			120 & 720 & 599458 & 8.21e+06 & 9.05e-05 & 4.25e-05 \\
		 \hline
		\end{tabular}
		\\
	\end{center} 
	\caption{For inner points the stencil size is $n=13$ with an hexagonal uniform grid, at boundary points the stencil size is $n_b=25$, the number of nodes per wave length is kept constant with $N_g=6$. Solution for the system $-\mathbf{H}\widetilde{U}=G$ is obtained for LU factorization. Error scale respect to the frequency as $\|u_1-\widetilde{u}_1\|_{\infty}\sim(3.9e-07)(\omega/2\pi)+4.7e-05$ and $\|u_2-\widetilde{u}_2\|_{\infty}\sim(1.4e-07)(\omega/2\pi)+3.3e-05$ } \label{table:fang1hexagon_pollution}		
\end{table}

\begin{table}
	\begin{center}
		\begin{tabular}{cccc|cc}
			\rowcolor{gray!20} $\frac{\omega}{2\pi}$  &  $\frac{1}{h}$ & $N$ & $\kappa(\mathbf{H})$ & $\|u_1-\widetilde{u}_1\|_{\infty}$ & $\|u_2-\widetilde{u}_2\|_{\infty}$   \\
			\hline\hline 
		    10 & 60 & 4237     & 8.18e+04 & 1.97e-05 & 1.23e-05   \\  
			20 & 120 & 16752   & 8.58e+05 & 2.22e-05 & 1.77e-05   \\
			40 & 240 & 66861   & 6.12e+05 & 1.94e-05 & 1.47e-05   \\
			80 & 480 & 266680  & 5.22e+06 & 1.78e-05 & 1.28e-05   \\
			120 & 720 & 599458 & 9.53e+06 & 1.85e-05  & 9.86e-06  \\
			\hline
		\end{tabular}
		\\
	\end{center} 
	\caption{Results for approximated solutions of \eqref{eq:helmholtz_impedance}. This approximation was calculated from a square uniform grid of $\Omega\cap\partial\Omega$. For inner nodes the stencil size is $n=19$, at boundary nodes $n=25$, the number of nodes per wave length is kept constant with $N_g=6$. Solution for the system $-\mathbf{H}\widetilde{U}=G$ is obtained for LU factorization. Errors scale respect to the frequency as $\|u_1-\widetilde{u}_1\|_{\infty}\sim( -2.5e-08)(\omega/2\pi)+2.1e-05$ and $\|u_2-\widetilde{u}_2\|_{\infty}\sim(-4.4e-08)(\omega/2\pi)+ 1.6e-05$ } \label{table:fang1hexagon19_pollution}		
\end{table}

\begin{table}
	\begin{center}
		\begin{tabular}{cccc|cc}
			\rowcolor{gray!20} NPW=$N_g$  &  $\frac{1}{h}$ & $N$ & $\kappa(\mathbf{H})$ & $\|u_1-\widetilde{u}_1\|_{\infty}$ & $\|u_2-\widetilde{u}_2\|_{\infty}$   \\
			\hline\hline 
		6.0 & 120.0 & 14641 & 7.80e+04 & 2.54e-04 & 1.70e-04 \\
		8.6 & 171.4 & 29584 & 1.05e+05 & 2.67e-05 & 1.56e-05 \\
		12.2 & 244.9 & 60025 & 1.47e+05 & 3.07e-06  & 1.79e-06 \\
		17.5 & 350.0 & 122500 & 1.99e+05 &3.55e-07  & 2.06e-07 \\
		25.0 & 500.0 & 251001 & 2.99e+05& 4.44e-08  & 2.55e-08 \\
			\hline
		\end{tabular}
		\\
	\end{center} 
	\caption{With $\omega/2\pi=20$, stencil size for inner nodes $n=9$ and stencil size for boundary nodes $n_b=15$ taken in a square uniform grid. Errors scale respect to $h$  as $\|u_1-\widetilde{u}_1\|_{\infty}\sim( 1.74e+09) h^{6.2}$ and  $\|u_2-\widetilde{u}_2\|_{\infty}\sim (3.2e+09) h^{6.4}$.} \label{table:fang1convergence_square}		
\end{table}

\begin{table}
	\begin{center}
		\begin{tabular}{cccc|cc}
			\rowcolor{gray!20} NPW=$N_g$  &  $\frac{1}{h}$ & $N$ & $\kappa(\mathbf{H})$ & $\|u_1-\widetilde{u}_1\|_{\infty}$ & $\|u_2-\widetilde{u}_2\|_{\infty}$   \\
			\hline\hline 
	    	6.0 & 120.0 & 16752 & 3.39e+17 & 4.43e-04 & 3.45e-04   \\
	    	8.6 & 171.4 & 33960 & 2.29e+15 & 2.07e-04  & 1.01e-04  \\
			12.2 & 244.9 & 69338 & 2.10e+07 & 8.76e-05 & 7.58e-05  \\
			17.5 & 350.0 & 141753 & 2.65e+06 & 7.67e-06 & 1.04e-06 \\
	    	25.0 & 500.0 & 289291 & 2.14e+20 & 8.50e-07 & 3.82e-07   \\
			\hline
		\end{tabular}
		\\
	\end{center} 
	\caption{With $\omega/2\pi=20$, stencil size for inner nodes $n=13$ and stencil size for boundary nodes $n_b=19$, taken in an hexagonal uniform grid. The errors scale respect to $h$  as $\|u_1-\widetilde{u}_1\|_{\infty}\sim(1613 ) h^{3.1}$ and  $\|u_2-\widetilde{u}_2\|_{\infty}\sim (5.5e+04) h^{3.9}$.} \label{table:fang1convergence_hexa}		
\end{table}

\subsubsection{Test 2}

\begin{figure}[ht]
	\begin{center}
		\begin{tabular}{cc}			
			\includegraphics[width=7.7cm]{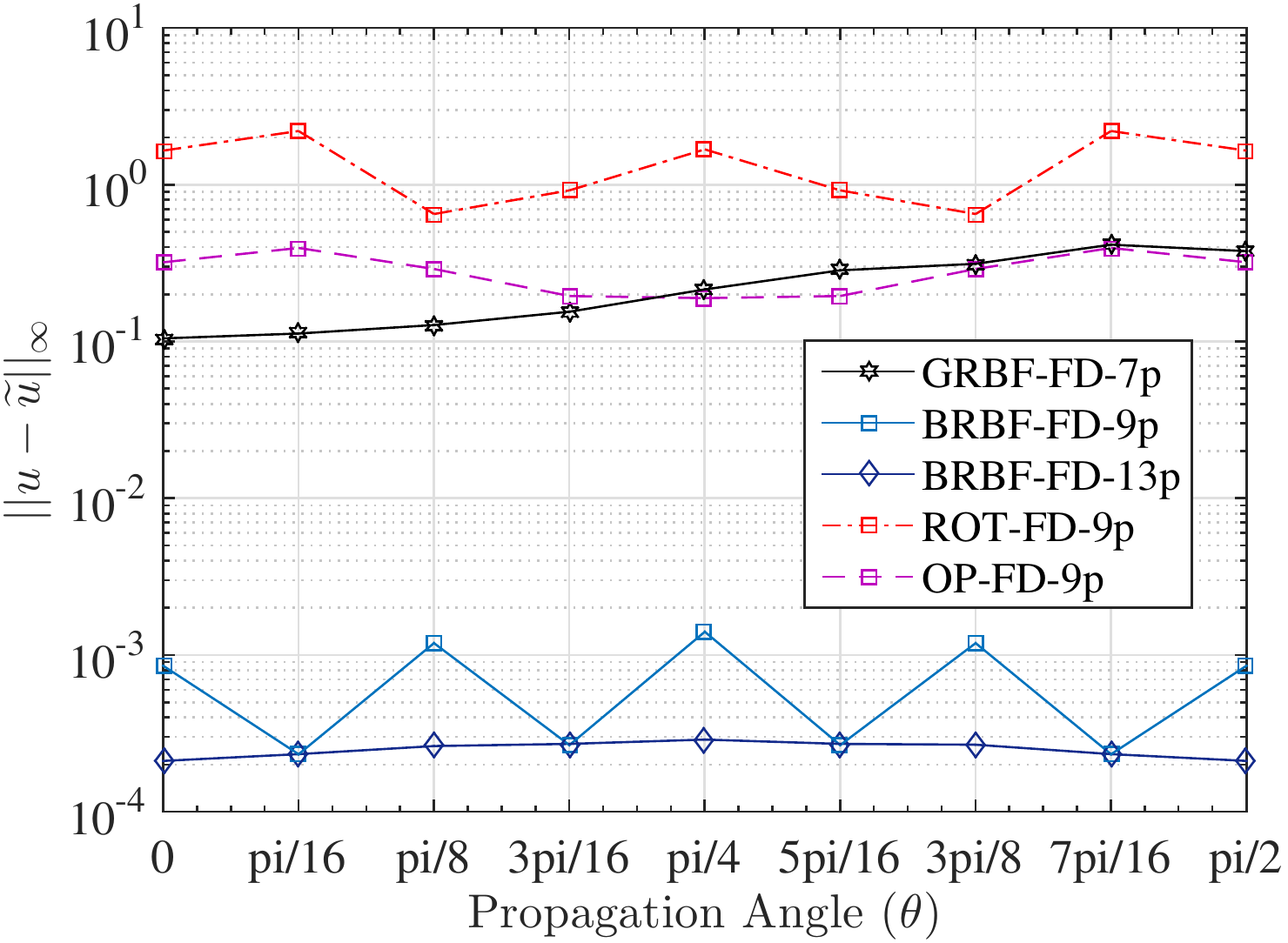} &          \includegraphics[width=7.7cm]{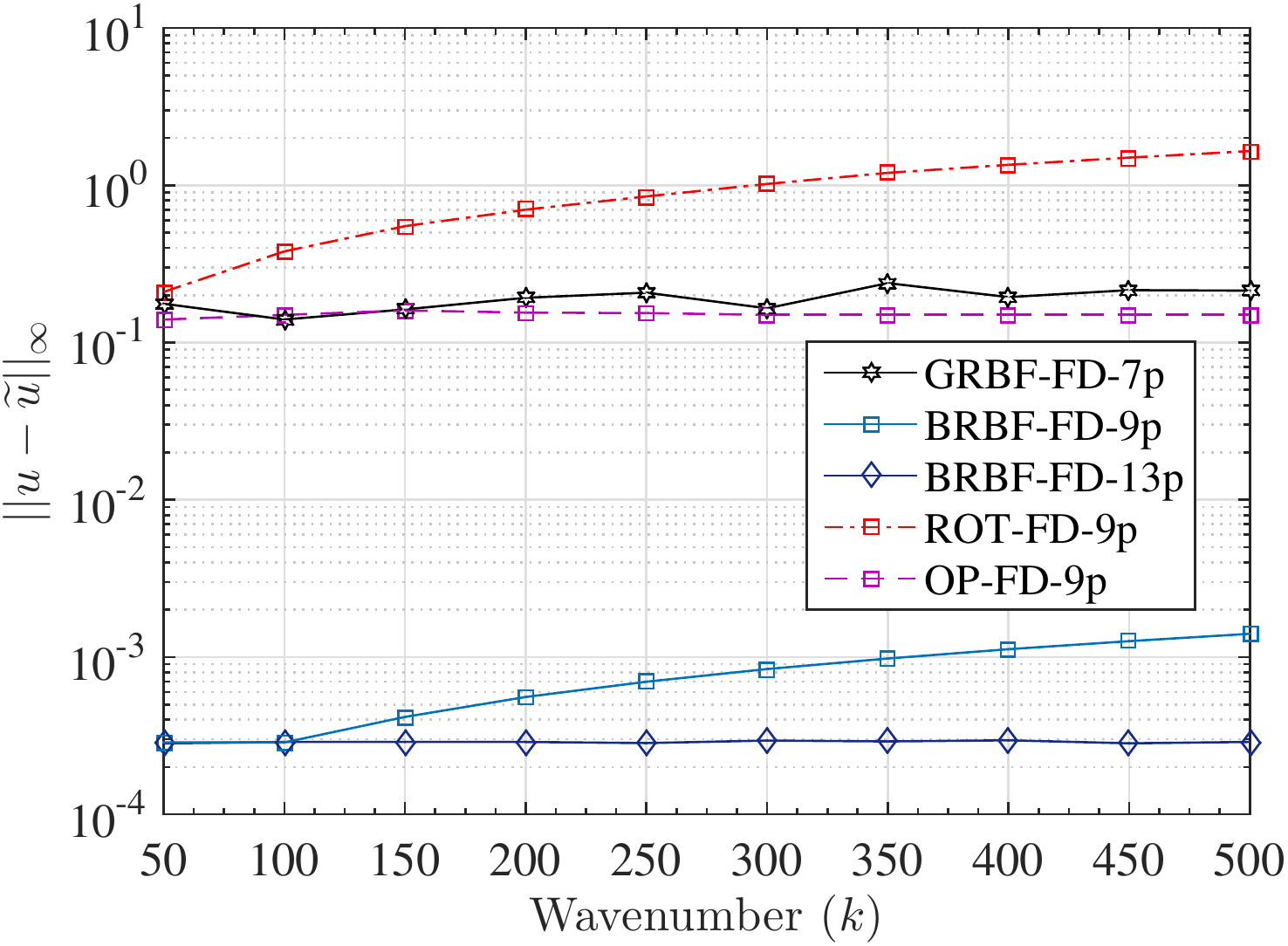}  			
		\end{tabular}       	
		\caption{Comparison of results between BRBF-FD-9p, BRBF-FD-13p,  GRBF-FD-7p and those reported in \cite{chen9}. (Left) Results for $k=500$ and $h=1/500$ varying the propagation angle. (Right) With $k$ varying, $\theta=\pi/4$ and $h=1/k$.}
		\label{fig:testpoll}
	\end{center}
\end{figure}

In this example we consider the problem

\begin{equation}
\begin{cases}\label{eq:abc_problem}
-\Delta u(x,y)-k^2u(x,y)=0& \mbox{ in } \Omega\\
\frac{\partial u}{\partial \mathbf{n}}u(x,y)+iku(x,y)=g(x,y)&  \mbox{ on } \partial\Omega,
\end{cases}
\end{equation}
with $\Omega=(0,1)\times(0,1)$, which has analytic solution given by the plane wave $u(x,y;k,\theta)=e^{ik(x\cos\theta+y\sin\theta)}$  when the data $g$ in the impedance boundary condition is given by
\begin{equation}
g(x,y)=
\begin{cases}
i(k-k_2)e^{ik_1x} & \mbox{ if } x\in (0,1) \mbox{ and } y=0\\
i(k+k_1)e^{i(k_1+k_2y)}& \mbox{ if } x=1 \mbox{ and } y\in(0,1)\\
i(k+k_2)e^{i(k_1x+k_2)}& \mbox{ if } x\in (0,1) \mbox{ and } y=1\\
i(k-k_1)e^{ik_2y}& \mbox{ if }  x=0 \mbox{ and } y\in(0,1)
\end{cases}
\end{equation}
with $k_1=k\cos\theta$ and $k_2=k\sin\theta$. This example is a standard benchmark to test the numerical dispersion of solvers for Helmholtz equation. The approximated solutions for this problem was calculated using Gaussian RBF-FD on hexagonal grids with 7-stencil (GRBF-FD-7p), Bessel RBF-FD with 9-stencils (BRBF-FD-9p) and 13-stencils (BRBF-FD-13p) on uniform Cartesian grids and we compare with results reported in \cite{chen9}. In all methods it has fixed $\mbox{NPW}=2\pi$ to keep constant resolution. For GRBF-FD-7p  we have used the approximations in \eqref{eq:wlap} with shape parameter $\ve_{op}$  given by \eqref{eq:optimal_shape_p} to approximate the Laplace operator and all partial derivative operators involved in the boundary condition. To solve the local interpolations in BRBF-FD-9p and BRBF-FD-13p we have used the condition number $\kappa_0=10^{7+1/hk+\sqrt{n_s}}$.  Our results can be seen in Fig. \ref{fig:testpoll}. On left plot we can see that there is less anisotropy in the error of BRBF-FD-13p. On right plot we can see that the dispersion and pollution effects are mitigated. Besides, we can see that the behavior of the error in the case GRBF-FD-7p is similar to reported in \cite{chen9}; in cases BRBF-FD-9p and BRBF-FD-13p the error it improvements at least with two accuracy orders.  We can stand out results for Bessel-RBF-13p, since even in increasing the wavenumber by fixing the resolution $kh=1$, the error remains constant. We can say about GRBF-FD-7p that although its accuracy is rather modest compared to BRBF-FD, its low complexity and easy implementation it makes a feasible option.  

The aim of this chapter is to show approximated solutions for some Helmholtz problems to test the performance of  Gaussian RBF-FD and Bessel RBF-FD schemes. We have compared with the fundamental solution of Helmholtz equation and also, we have computed solutions with smooth and non-continuous coefficients to have a qualitative sight. At last section we show results of inverted data with FWI,  applied in three well-known benchmarks: Marmousi, 2D Oversthrust and 2004 BP velocity model. 

All codes were typed in Matlab R2016a and run in a laptop \textcolor{black}{ with Core i7 processor at 2.8 Ghz with 12 Gbytes of RAM.}
\section{Numerical results for some Helmholtz problems}\label{sec:numerical_tests}

In this section we test the BRBF-FD scheme to computing solutions with second and third order ABC.

\subsection{Approximated fundamental solutions}

It is known that the problem $-\Delta u-k^2u=\delta(\x)$ in the free-space $\R^2$ has a unique solution when it imposes the the Sommerfeld radiation condition
\begin{equation}
\lim_{\|\x\|\rightarrow\infty}\|\x\|^\frac{1}{2}\left(\frac{\partial}{\partial r}-ik\right)u(\x)=0.
\end{equation}
Particularly, the associated Green's function, which is solution of $-\Delta u(\x)-k^2u(\x)=\delta(\x-\x_0)$,  \cite{Ciraolo2009_radiation_condition_for_the_2D_helmholtz_equation}  is given by $u(\x)=G(\x,\x_0)$, with
\begin{equation}
G(\x,\x_0)=\frac{i}{4}H_0^{(1)}(k\|\x-\x_0\|)
	\end{equation}

\begin{figure}[ht]
	\includegraphics[scale=0.8]{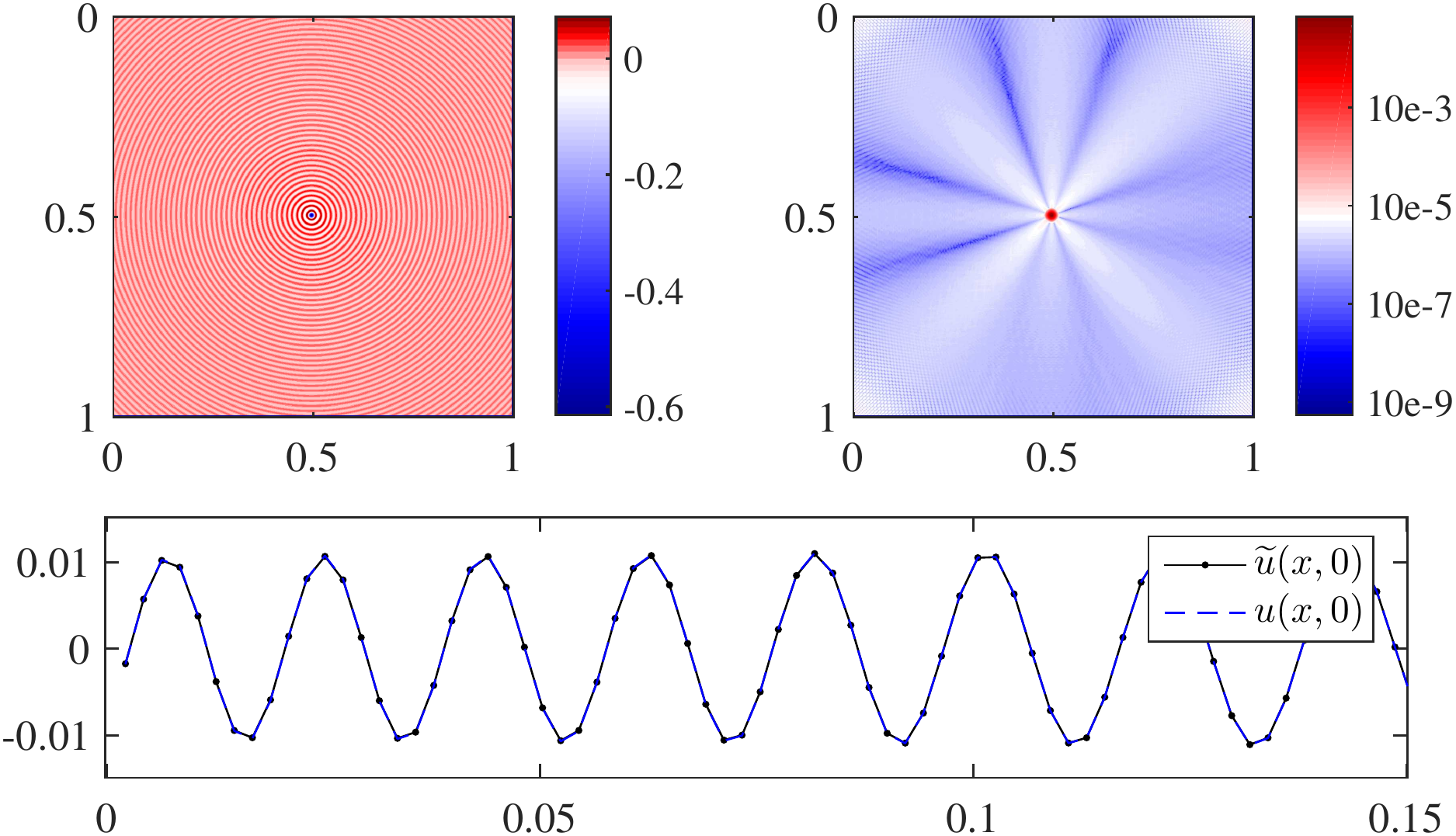}
	\caption{(Top-left), approximated solution $u(\x)$ (top-right) plot of  $|\widetilde{u}(\x)-u(\x)|$, (bottom) comparison on boundary values. Source at $\x=(0.5,0.5)$}\label{fig:hankel_c_500}
\end{figure}

\begin{figure}[ht]
		\includegraphics[scale=0.8]{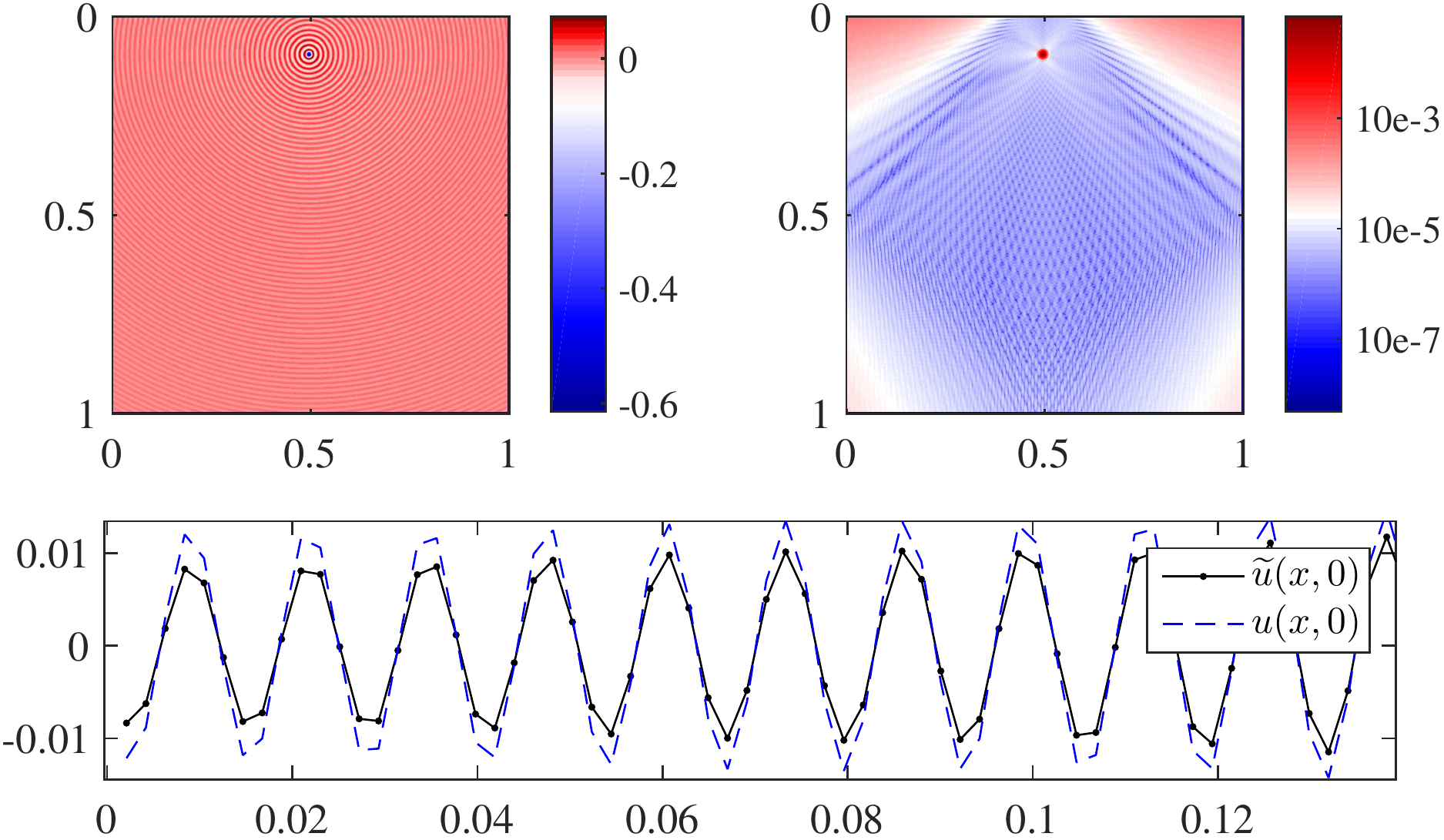}
		\caption{(Top-left), approximated solution $u(\x)$ (top-right) plot of  $|\widetilde{u}(\x)-u(\x)|$, (bottom) comparison on boundary values. Source at $\x=(0.5,0.1)$}\label{fig:hankel_m_500}
\end{figure}

We compute approximated Green's functions in the free space truncated to a bounded domain $\Omega\subset\R^2$ for $\x\neq\x_0$, where $\x_0\in\Omega$,  through the boundary problem
 \begin{equation}\label{eq:helmholtz_impedance}
 \left\{
 \begin{array}{rcll}
 -\Delta u(\x)-k^2u(\x)&=&\widetilde{\delta}(\x-\x_0), &\mbox{ in } \Omega\\
 \frac{\partial}{\partial\n}u(\x)+ik \left(1+\frac{3}{4k^2}\frac{\partial^2}{\partial\bm{\tau}^2}-\frac{i}{4k^3}\frac{\partial^3 }{\partial \mathbf{n}\partial\bm{\tau}^2 }\right) u(\x)&=&0, & \mbox{ on } \Gamma=\partial\Omega,
 \end{array}
 \right.
 \end{equation}
 where the boundary condition corresponds to the ABC in the Padé approximation. 
  The single source is given by the Gaussian function
\begin{equation}\label{eq:single_source}
\widetilde{\delta}(\x-\x_0)=\frac{1}{2\pi\sigma^2}e^{-\frac{\|\x-\x_0\|}{2\sigma^2}}
\end{equation}
for a value of $\sigma$ such that $\int_{\Omega}\widetilde{\delta}(\x-\x_0)d\x\approx1$.
We have used the BRFB-FD scheme with a square regular grid in $\Omega=(0,1)\times(0,1)$ with square $9$-stencils at inner nodes and $19$-stencils at boundary nodes. The results
shown in Fig. \ref{fig:hankel_c_500} and Fig. \ref{fig:hankel_m_500} where calculated with wavenumber $k=500$ with NPW$=6$, i.e., $h=2\pi/6k$. We point out that results show a good accuracy at $\x\neq\x_0$ and on the boundary the wavelength of the numerical solution matches very good with the exact one when the source is put at the center of the square domain, this is a good indication that dispersion errors are not significant, however in Fig. \ref{fig:hankel_m_500} (bottom) we see that amplitude has a considerable discrepancy respect to the exact one, this is due to the approximated ABC.


\subsection{Heterogeneous medium}

\subsubsection{Smooth medium}
\begin{figure}[ht]
	\begin{center}
		\begin{tabular}{ll}		
			\includegraphics[scale=0.45]{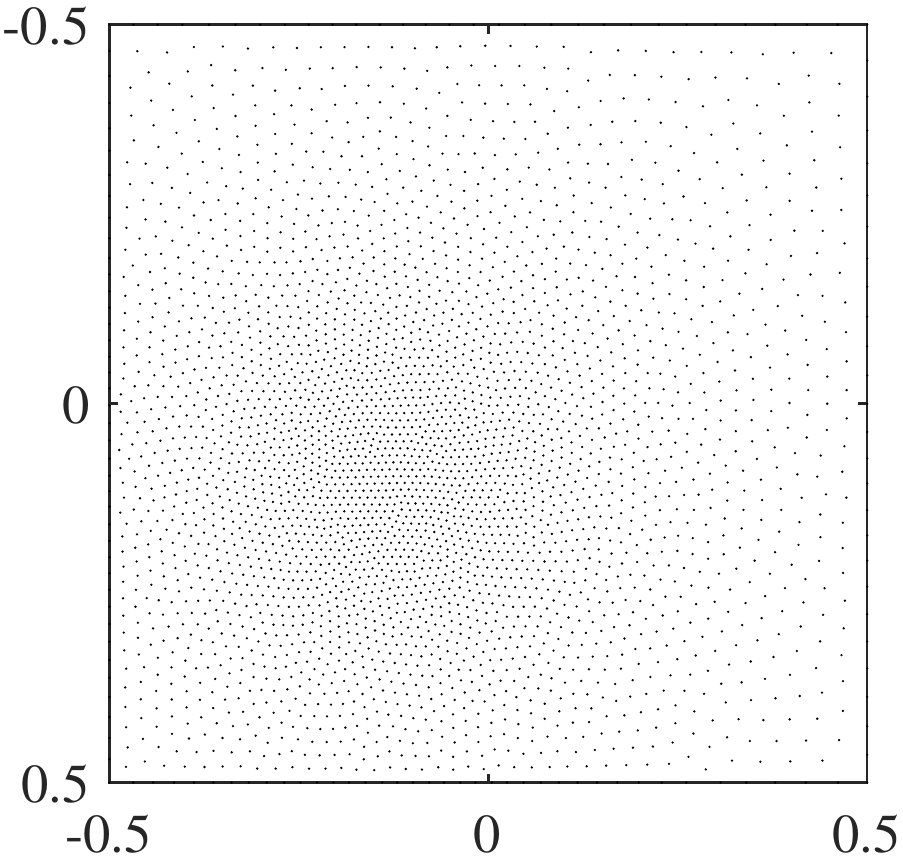} &\includegraphics[scale=0.6]{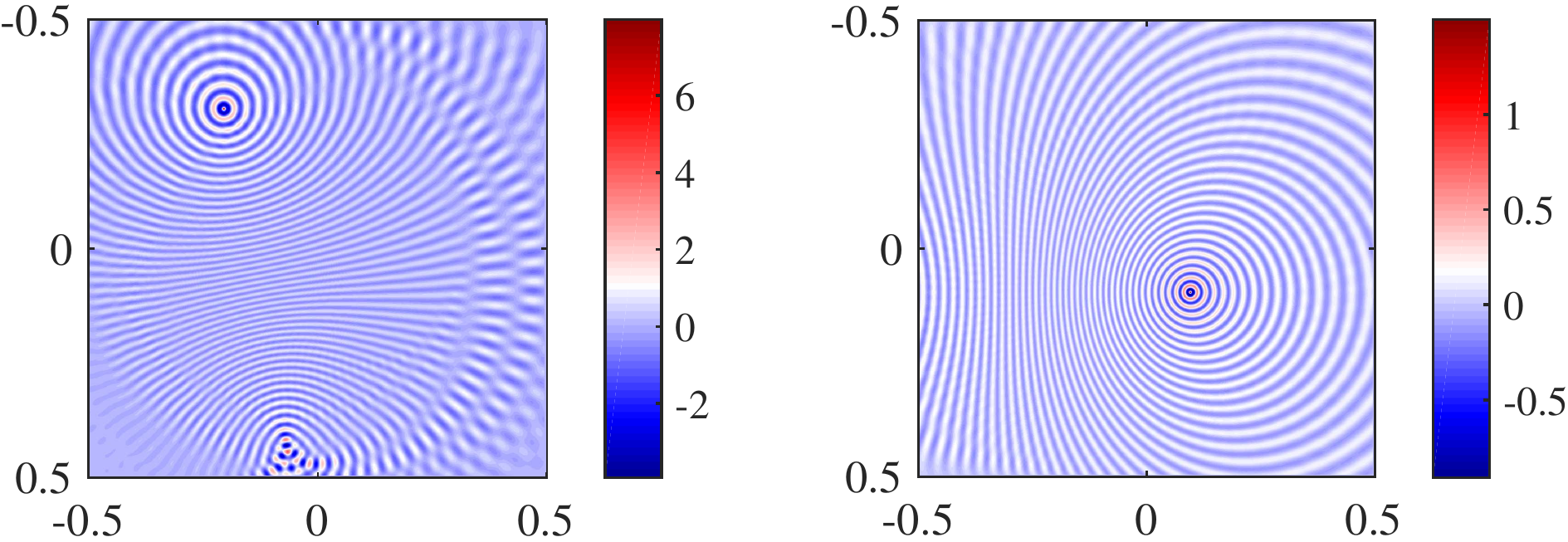} \\		
			\includegraphics[scale=0.45]{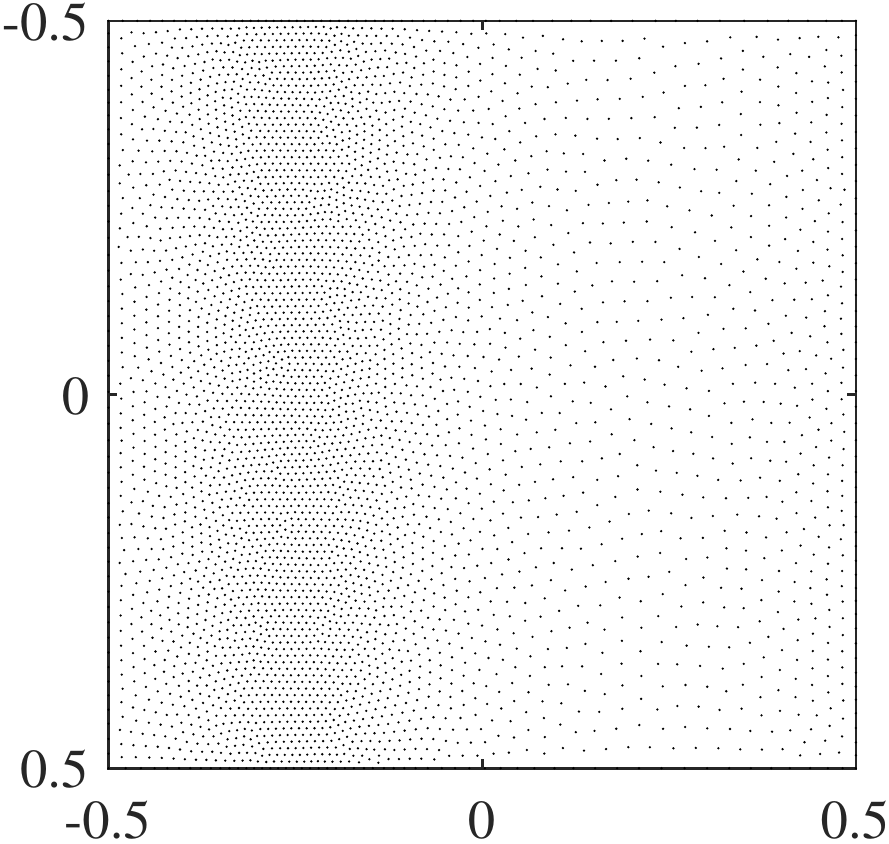} &\includegraphics[scale=0.6]{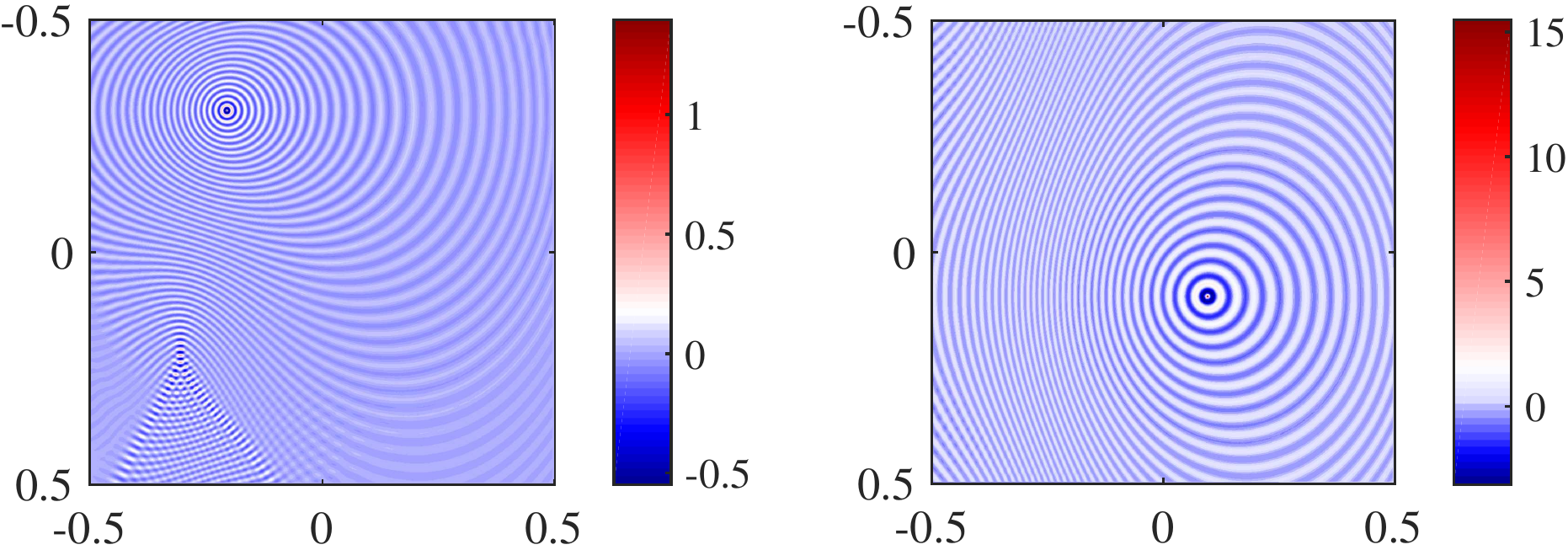} \\
		\end{tabular} \caption{Plot of solutions corresponding to velocity models in \eqref{eq:speed_fang3} and \eqref{eq:speed_fang4}. (Left column) Sketch of their node distribution according to the local wavelength and (right columns) wave fields corresponding to two different single sources.}
		\label{fig:fang34}
	\end{center}
\end{figure}

For this qualitative test, we have calculated approximated solutions of the problem
\begin{equation}\label{eq:helmholtz_ABC2}
\left\{
\begin{array}{rcll}
-\Delta u(\x)-\omega^2c(\x)^{-2}u(\x)&=&\widetilde{\delta}(\x-\x_0), &\mbox{ in } \Omega\\
\frac{\partial}{\partial\n}u(\x)+i\omega c(\x)^{-1}\mathcal{B}u(\x)&=&0, & \mbox{ on } \Gamma=\partial\Omega
\end{array}
\right.
\end{equation}
where $\mathcal{B}$ is  the operator $\mathcal{B}=1+\frac{c(\x)^2}{2\omega^2}\frac{\partial^2}{\partial\bm{\tau}^2}$ corresponding to the ABC of second order and $\Omega=(-0.5,0.5)\times(-0.5,0.5)$. Here we perform two examples with the velocity functions
\begin{equation}
c(x,y)=3 - 2.5e^{  -( (x+0.125)^2 + (y-0.1)^2 )/0.8^2  }
\label{eq:speed_fang3}
\end{equation} 
and
\begin{equation}
c(x,y)=1+0.5\sin(2\pi x).
\label{eq:speed_fang4}
\end{equation}
For these velocity models, nodes distributions are sketched on the left column of Fig. \ref{fig:fang34}. on center and right columns it can be seen the real part of the approximated solution for two different single sources. In table \ref{table:fang34} it shows results for required times to assembly sparse matrices $\mathbf{H}=\mathbf{H}_{\Omega}+\mathbf{H}_{\Gamma}$ and for solution of the system $-\mathbf{H}U=\mathbf{F}$ by LU factorization.

\begin{table}[ht]
	\begin{center}
		\begin{tabular}{ccc|ccc}
			\rowcolor{gray!20} $\omega/2\pi$  & Nodes ($N$)  &  $\kappa(\mathbf{H})$  &Time (s) for $\mathbf{H}_{\Omega}$ & Time (s) for $\mathbf{H}_{\partial\Omega}$ & Time (s) for LU  \\
			\hline\hline 
			2.5 & 1072 & 1.09e+04 & 0.37 & 0.24 & 0.02 \\
			5 & 4404 & 1.46e+05 & 1.18 & 0.45 & 0.10  \\
			10 & 17563 & 6.25e+05 & 4.72 & 1.01 & 0.45  \\
			20 & 70585 & 1.93e+06 & 18.98 & 2.58 & 3.00 \\
			40 & 283458 & 1.27e+07 & 75.10 & 9.05 & 17.14   \\
			\hline
		\end{tabular}
		\\
	\end{center} 
	\caption{Results in computing solutions corresponding to the smooth velocity model in \eqref{eq:speed_fang3}.} \label{table:fang34}		
\end{table}

\subsubsection{Non smooth medium (Test in the 2004 BP velocity model)}

\begin{figure}[ht]
	\begin{center}
		\begin{tabular}{c}
			\includegraphics[scale=0.48]{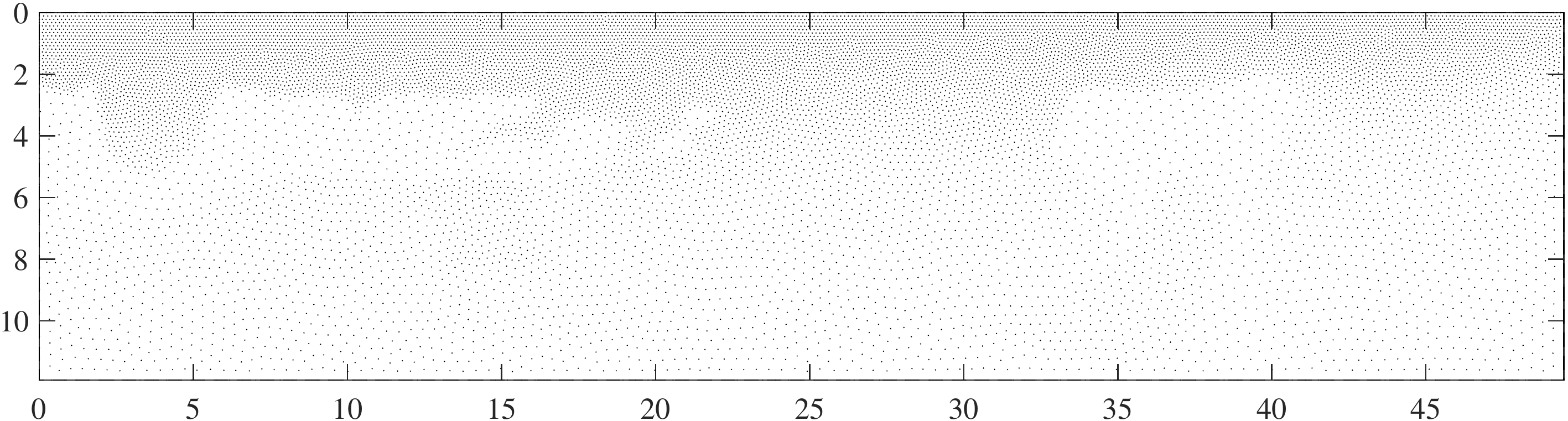}\\
			\includegraphics[scale=0.90]{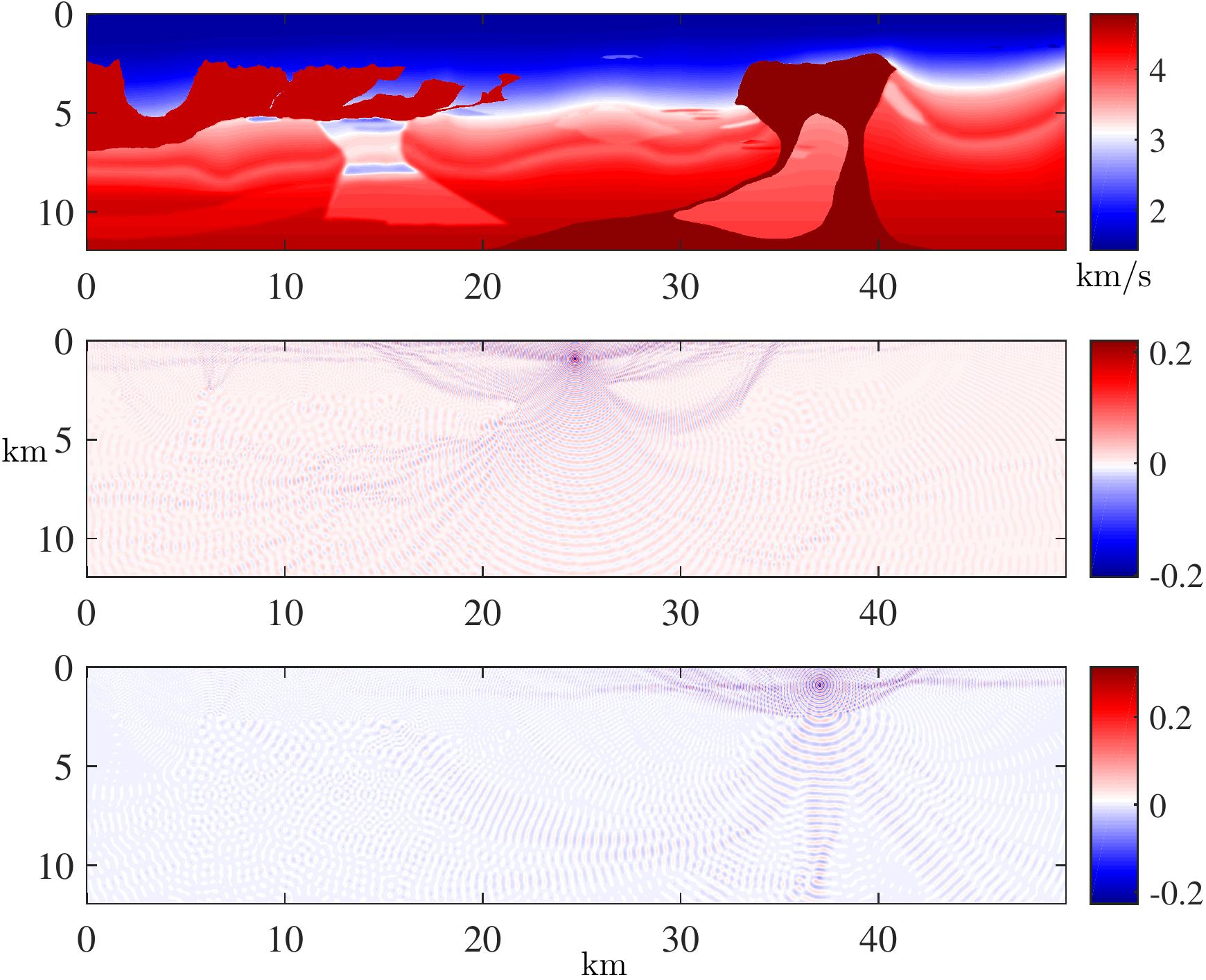}
		\end{tabular}
	\end{center}
	\caption{(Top) Sketch of node distribution for 2004BP model used for a frequency $\omega/2\pi=2$Hz with NPW=10 of the (Middle-top) velocity model.  (Middle-Bottom and bottom) plots of real part of the wavefield with to two single sources located at different positions for a frequency  $\omega/2\pi=10$Hz. In according to the local wavelength, the domain is discretized with $N=345393$, at inner nodes the local interpolation is performed with $19$-stencils and at boundary nodes with $25$-stencils.}\label{fig:2004bp_model}. 
\end{figure}

We consider the 2004 BP velocity benchmark, which is a popular model in research for velocity estimation methods in seismic imaging, which is presented as a challenge due to its complexity and large scale \cite{2004bp}. The  velocity function $c$ can be seen in the density plot in the middle-top in Fig. \ref{fig:2004bp_model}. Roughness of the velocity model such as hard interfaces and
sharp transitions generates strong reflections that hinder the efficiency of known iterative methods due to the ill-conditioning of the matrix gets worse \cite{ZEPEDA_NUNEZ2016347}, increasing the number of iterations. Moreover, for large $\omega$ the interaction of high frequency waves with short wavelength structures such as discontinuities, increases the reflections, further deteriorating the convergence rate. In BRBF-FD scheme the local interpolation matrices $\JJ_k$ becomes dramatically ill-conditioned, however by keeping the condition number of $\widetilde{\JJ}_k$ to the value $\kappa_0\approx 10^{(1+\sqrt{n_s})}$, we have empirically observed that the condition number of $\mathbf{H}$ lies in the range: $10^4\leq\kappa(\mathbf{H})\leq 10^9$, so it is feasible to perform LU factorization. Moreover, qualitatively, the wavelenghts of the wavefield has the expected behavior according to velocity model. To model wave propagation we have solved \eqref{eq:helmholtz_ABC2} with ABC of second order.

\begin{table}[ht]
	\begin{center}
		\begin{tabular}{ccc|ccc}
			\rowcolor{gray!20} $\omega/2\pi$  & Nodes ($N$)  &  $\kappa(\mathbf{H})$  &Time (s) for $\mathbf{H}_{\Omega}$ & Time (s) for $\mathbf{H}_{\partial\Omega}$ & Time (s) for LU  \\
			\hline\hline 
		2 & 38488  &  9.01e+04 & 20.59 & 4.72 & 0.99 \\
		4 & 153223 & 3.43e+05 & 74.27 & 10.52 & 5.42\\	
		6 & 345389 & 9.93e+05 & 173.78 & 19.75 & 16.04\\
		8 & 613529 & 7.22e+06 & 192.71 & 33.50 & 35.76\\	
			\hline
		\end{tabular}
		\\
	\end{center} 
	\caption{Results in computing solutions corresponding to the 2004 BP velocity model.} \label{table:2004bp}		
\end{table}

\section{Conclusions}
We perform local interpolation with a free shape parameter oscillatory RBF based on Bessel function of first kind to obtain a high order RBF-FD scheme for solving Helmholtz equation. In this approach the shape parameter is substituted for the local wavenumber. However due to the local interpolation matrices are extremely ill-conditioned even for small stencils, we overcome this issue with regularization by small perturbation of the diagonal. In some tests we have achieved convergence rates of third and sixth order. 

\section*{Acknowledgments}
This work is supported by Colombian Oil Company ECOPETROL and COLCIENCIAS as a part of the research project grant No. 0266-2013.
    
    \section*{References}
	\bibliographystyle{apalike}
	\bibliography{references.bib}

	
	
	
	
	

\end{document}

\documentclass[12pt,preprint]{elsarticle}

\usepackage{amssymb}



\usepackage[utf8]{inputenc}
\usepackage{graphicx}
\usepackage{dcolumn}
\usepackage{dsfont}
\usepackage{amsthm}
\usepackage{epstopdf}
\usepackage{textcomp}
\usepackage{colortbl}
\usepackage{float}
\usepackage{bm}
\usepackage[format=plain,font={small},justification=centerlast]{caption}
\usepackage{subcaption}
\usepackage{color}
\usepackage{tabularx,ragged2e,booktabs}
\usepackage{MnSymbol}
\usepackage{wasysym}
\usepackage{anysize}
\usepackage{epstopdf}
\epstopdfsetup{update}
\usepackage{verbatim}
\usepackage{tikz}
\usepackage{pgfplots}
\usepackage{mathtools}
\usepackage{natbib}

\DeclareGraphicsExtensions{.pdf,.png,.jpg,.pdf}

\makeatletter
\newenvironment{figurehere}
{\def\@captype{figure}}
{}
\makeatother

\biboptions{comma,square}


\def\im{\mathop{\rm \od{\iota}}\nolimits}
\newcommand{\BR}{{\mathbb{R}}}
\newcommand{\p}{\partial}

\newcommand{\ts}[1]{\textstyle #1}
\newcommand{\bn}[1]{\mbox{\boldmath $#1$}}
\newcommand{\bc}{\begin{center}}
	\newcommand{\ec}{\end{center}}
\newcommand{\be}{\begin{equation}}
\newcommand{\ee}{\end{equation}}
\newcommand{\bea}{\begin{eqnarray}}
\newcommand{\eea}{\end{eqnarray}}
\newcommand{\ba}{\begin{array}}
	\newcommand{\ea}{\end{array}}

\newcommand{\JJ}{\mathbf{J}}
\newcommand{\II}{\mathbf{I}}

\newcommand{\N}{\mathds{N}}
\newcommand{\x}{\mathbf{x}}
\newcommand{\n}{\mathbf{n}}
\newcommand{\e}{\mathbf{e}}
\newcommand{\y}{\mathbf{y}}
\newcommand{\ds}{\displaystyle}
\renewcommand{\a}{\mathbf{a}}
\renewcommand{\r}{\mathbf{r}}
\newcommand{\h}{\mathbf{h}}
\renewcommand{\v}{\mathbf{v}}
\renewcommand{\u}{\mathbf{u}}
\newcommand{\bO}{\mathcal{O}}
\newcommand{\Pcal}{\mathcal{P}}

\newcommand{\Z}{\mathds{Z}}
\newcommand{\Q}{\mathds{Q}}
\newcommand{\R}{\mathds{R}}
\newcommand{\C}{\mathds{C}}
\renewcommand{\P}{\mathds{P}}
\newcommand{\phiep}{\varphi_{\epsilon}}
\newcommand{\mt}{\mathcal{T}}
\newcommand{\mb}{\mathcal{B}}
\newcommand{\ml}{\mathcal{L}}
\newcommand{\ve}{\varepsilon}
\newcommand{\om}{\Omega}
\newcommand{\td}{\tilde}
\newcommand{\al}{\alpha}
\newcommand{\vp}{\varphi}

\newtheorem{thm}{Theorem}
\newtheorem{lemma}{Lemma}
\newtheorem{definition}{Definition}
\newtheorem{rmrk}{Remark}
\newtheorem{pro}{Proposition}

\journal{---}

\begin{document}
	
	\begin{frontmatter}
		
		
		\title{Radial basis function-generated finite differences with Bessel weights for the 2$D$ Helmholtz equation}
		

		
		\author{Mauricio A. Londoño}
		\ead{alejandro.londono@udea.edu.co}
		\author{Hebert Montegranario}
		\ead{hebert.montegranario@udea.edu.co}
		
		\address{Instituto de Matemáticas \\
			Universidad de Antioquia\\
			Calle 67 53-108, Medell\'in, Colombia}
		
		\begin{abstract}
			In this paper we obtain approximated numerical solutions for the 2D Helmholtz equation using a radial basis function-generated finite difference scheme, where  weights are calculated by taking the oscillatory radial basis function given in terms of Bessel functions of first kind. The problem of local interpolation to obtain weights is ill-conditioned and we overcome this difficulty by means of regularization of the interpolation matrix by perturbing the diagonal, where the condition number of the perturbed matrix is controlled according a prescript value. We perform different tests from which observe convergence, algorithm complexity, besides we see that the pollution-effects are mitigated. 
		\end{abstract}
		
		\begin{keyword}
			RBF-FD\sep Helmholtz equation \sep shape parameter \sep pollution effect \sep Oscillatory RBF.
			
			
		\end{keyword}
		
	\end{frontmatter}
	
	
	\section{Introduction}
    
    In this work we consider the 2$D$ Helmholtz equation
    
    \begin{equation}\label{eq:helmholtz1}
    \left\{
    \begin{array}{rcll}
    -\Delta u(\x)-\omega^2c(\x)^{-2}u(\x)&=&f(\x), &\mbox{ in } \Omega\\
    b\frac{\partial}{\partial\n}u(\x)+i\omega c(\x)^{-1}\mathcal{B}u(\x)&=&g(\x), & \mbox{ on } \Gamma=\partial\Omega
    \end{array}
    \right.
    \end{equation}
where $\omega$ is the angular frequency, $c(\x)>0$ is the sound speed of the continuous media, $f(\x)$ is the source term, $\n$ is unitary normal vector to the boundary $\Gamma$, $b$ takes values zero or one, $\mathcal{B}$ is a certain linear operator and $g(\x)$ is certain exact data on $\Gamma$.

    \section{Oscillatory RBF}

Trefftz methods for the Helmholtz problems are schemes of type finite elements where test and trial functions are local solutions of the differential equation to solve. In \cite{trefftz_paper} it can be seen a wide survey about. Inspired by Trefftz methods, in this section we work with a family of oscillatory RBF which is solution of the homogeneous Helmholtz equation. Besides,  given the oscillatory behavior of solutions of Helmholtz equation, it makes sense to consider a such family, whose members are given in terms of Bessel functions of the first kind.

Now, we are going to describe briefly the class of oscillatory radial basis functions (RBF), which are given by
\begin{equation}\label{eq:bessel_RBF}
\varphi_k^{(d)}(r)=\frac{J_{d/2-1}(kr)}{(kr)^{d/2-1}}, \ \ \ d=1,2,\ldots,
\end{equation}
 whose detailed study is presented in \cite{FORNBERG2006_oscillatory_RBF}. Here $J_{\alpha}(r)$ is denoting the Bessel function of the first kind and order $\alpha$. Two remarkable properties of these oscillatory RBF are:

\begin{itemize}

	\item the non-singularity of the interpolation matrix for arbitrarily scattered data in up to $d$ dimensions, when $d>1$,

	\item and that the Laplace eigenvalue problem $\Delta \varphi + k^2\varphi =0$ has as bounded solutions, at the origin, the functions given in \eqref{eq:bessel_RBF}, thus any interpolant of the form

	\begin{equation}\label{eq:oscilla_interpolant}
	s(\x)=\sum_{j=1}^n\alpha_j\varphi^{(d)}_k(||\x-\x_j||
	\end{equation}

	will satisfy too $\Delta s+k^2 s=0$.  

\end{itemize}

In the case $d=2$ the oscillatory RBF \eqref{eq:bessel_RBF} have been little bit studied due to the function \eqref{eq:oscilla_interpolant} does not have local maximum at points where this is negative, which put a restriction to be used for general 2D interpolation. But in this work a such feature becomes a strength, since that we are just interpolate solutions of Helmholtz problems, which locally can be seen as plane waves satisfying the homogeneous Helmholtz equation. In early works, as in \cite{LIN2012_oscillatory_radial_Basis_function_helmholtz}, oscillatory RBF based on Bessel functions have been employed to solve the 2D Helmholtz equation with constant wavenumber within the approach of global collocation method.
They presented the RBF  
\begin{equation}\label{eq:RBFlin2012}
\phi_{C,k}(r)=J_0(k\sqrt{r^2+C^2}),
\end{equation}
which has two shape parameters but fixed $k$ to the wavenumber. The ill-conditioning of the interpolation matrix arises from \eqref{eq:RBFlin2012} is overcome by way of a regularized singular value decomposition method. 

For our interest, the Helmholtz problem 2D with large wavenumber, we take the special case $d=2$. So we work with the oscillatory Bessel-RBF
\begin{equation}\label{eq:bessel_RBF2D}
\phi_{k}(r)=J_0(kr),
\end{equation}
where in the approach RBF-FD the shape parameter $k$ will be evaluated at the wave number $k(\x)=\omega/c(\x)$ corresponding to the center of the stencil. As is well known, to assemble the sparse matrix, which discretizes the Helmholtz problem is necessary to solve a small linear equation system at each node.  As we will see, interpolation matrices are ill-conditioned and we will overcome this issue with the method of diagonal increments (MDI) \cite{method_of_diagonal_increments}, \cite{SARRA2014_regularized_positive_definite} adding to the diagonal entries a small regularization parameter $\beta>0$, thus we solve, instead of the linear system  $\mathbf{b}=\mathbf{A}\mathbf{y}$,  the equation
\begin{equation}
\mathbf{b}=(\mathbf{A}+\beta\mathbf{I})\widetilde{\mathbf{y}},
\end{equation}
where $\II$ is the identity matrix. We will give a brief explanation about MDI in the section \ref{sec:MDI}, where it shows that the matrix $\widetilde{\mathbf{A}}=\mathbf{A}+\beta\II$ is better conditioned than $\mathbf{A}$ and  $\mathbf{y}\approx\widetilde{\mathbf{y}}$. Before we describe the goodness of discretizing Helmholtz problems with Bessel-RFB.


\section{Bessel-RBF-FD}

Suppose that $u$ is a solution of the Helmholtz equation 
$\Delta u(\x)+k(\x)^2u(\x)=0$, for $\x\in\Omega$ and for $\x\in\partial\Omega$ it satisfies certain boundary condition. If $X=\{\x_i\}_{i=1}^N\subset\Omega\cup\partial\Omega$ is a set of nodes, for $\x_i\in X\cap\Omega$ we take a stencil $S_i=\{\x_j^i\}_{j=1}^{n_i}\subset X$ based on $\x_i$, with $\x_1^i=\x_i$. For $\x\in\Omega_i=\mbox{ConvexHull}(S_i)$ we define, with $k_i=k(\x_i)$, the interpolant 
\begin{equation}\label{eq:interpolant_bessel}
\widetilde{u}(\x)=\sum_{j=1}^{n_i}\alpha_j^i J_0(k_i\|\x-\x_j^i\|).
\end{equation}
With the local interpolation matrix, $\mathbf{J}_{k_i}=(J_0(k_i\|\x_l^i-\x_j^i\|))_{1\leq l,j\leq n_i}$, which is positive definite, and forcing the condition $\widetilde {u}(\x_l^i)=u(\x_l^i)$, then from \eqref{eq:interpolant_bessel} we have the linear equation
 \begin{equation}\label{eq:system_interpolation_ill}
 U_i=\mathbf{J}_{k_i}\bm{\alpha}_i,
 \end{equation}
where $U_i=\left(\begin{array}{cccc}
u(\x_1^i)&u(\x_2^i)&\cdots&u(\x_{n_i}^i)
\end{array}\right)^t$ and $\bm{\alpha}_i=\left(\begin{array}{cccc}
\alpha_1^i&\alpha_2^i&\cdots&\alpha_{n_i}^i
\end{array}\right)^t$.
In view that $\phi_k$, defined in \eqref{eq:bessel_RBF2D}, satisfies the homogeneous Helmholtz equation, then
\begin{eqnarray}
\Delta\widetilde{u}(\x)|_{\x=\x_i}&=&\sum_{j=1}^{n_i}\alpha_j^i \Delta J_0(k_i\|\x-\x_j^i\|)|_{\x=\x_i}\\
&=&-k_i^2\sum_{j=1}^n\alpha_j^i J_0(k_i\|\x_i-\x_j^i\|)\\
&=&-k_i^2\widetilde{u}(\x_i).
\end{eqnarray}
Hence the interpolant \eqref{eq:interpolant_bessel} satisfies the homogeneous Helmholtz equation. On the other hand, applying  $\Delta_{S_i,k_i}$ to the solution $u$, we have
\begin{eqnarray}
\Delta_{S_i,k_i}u(\x_i)&=&\Delta \mathbf{J}_{k_i}\mathbf{J}_{k_i}^{-1}U_i\\
&=&-k^2_i\mathbf{e}_1U_i\\
&=&-k_i^2u(\x_i)
\end{eqnarray}
Note that $\Delta_{S_i,k_i}u(\x_i)-\Delta_{S_i,k_i}\widetilde{u}(\x_i)=-k_i^2(u(\x_i)-\widetilde{u}(\x_i))$, thus, for solutions of the homogeneous Helmholtz problem the local truncation error for the Laplace operator has a theoretical error depending of wavenumber at $\x_i$ and of the error of the local interpolant. The error of the approximated  solutions is produced by the interpolant \eqref{eq:interpolant_bessel} and by the ill-conditioning of the matrix $\JJ_{k_i}$, in solving the linear system \eqref{eq:system_interpolation_ill}. Next we will deal with solutions fo these systems.

\section{Method of diagonal increments (MDI)}\label{sec:MDI}

The interpolation matrix $\JJ_k$ is ill-conditioned, especially for certain node distributions. In literature there are several methods to overcome the ill-conditioning when the shape parameter is small \cite{stable_general_bessel_rbf}, ever in our case we are taking the shape parameter as the wavenumber $k$, which can be large. So we have chosen the MDI. For our case $\JJ_k$ will be considered  ill-conditioned when the condition number, $\kappa(\JJ_k)$, satisfies  $\kappa(\JJ_k)>10^{15}$, which hinders that the solution $\bm{\alpha}$ of $U_i=\mathbf{J}_k\bm{\alpha}$ be accurately calculated, in double precision, 	through Cholesky facorization. When $\JJ_k$ is ill-conditioned we solve the better conditioned problem $U_i=(\mathbf{J}_k+\beta\mathbf{I})\widetilde{\bm{\alpha}}$ instead, where $\II$ is  the identity matrix and $\beta$  a small positive real number.  Next, we will give some important aspects about the spectrum of $\JJ_k+\beta\II$ in order to give an estimate of the error $\|\bm{\alpha}-\widetilde{\bm{\alpha}}\|_2$.

\begin{rmrk}\label{rmrk:convergence_neumann_series_jk}

	$\JJ_k$ is positive definite, thus its spectrum is real and positive. If $\{\lambda_m\}_{m=1}^n$ is its spectrum,  with $\lambda_1\geq\lambda_2\geq\cdots\geq\lambda_n$, then $\{\lambda_m+\beta\}_{m=1}^n$ is the specturm of $\widetilde{\JJ}_k=\JJ_k+\beta\II$ and $\{\frac{\beta}{\lambda_m+\beta}\}_{m=1}^n$ is the spectrum of $\beta\widetilde{\JJ}_k^{-1}$, hence its spectral norm is $\|\beta\widetilde{\JJ}_k^{-1}\|_2=\frac{\beta}{\lambda_n+\beta}$. The  above implies that the Neumann series  $\sum_{m=0}^{\infty}(\beta\widetilde{\JJ}_k^{-1})^m$ converges and the equality 
	\begin{equation}
	(\II-\beta\widetilde{\JJ}_k^{-1})^{-1}=\sum_{m=0}^{\infty}(\beta\widetilde{\JJ}_k^{-1})^m
	\end{equation} holds.

\end{rmrk}

\begin{rmrk}\label{rmrk:better_coditioning_wideJk}

	If $\{\lambda_m\}_{m=1}^n$ is the spectrum  of $\JJ_k$,  then the condition number of $\JJ_k$ is given by $\kappa(\JJ_k)=\frac{\lambda_1}{\lambda_n}$ and $\kappa(\widetilde{\JJ}_k)=\frac{\lambda_1+\beta}{\lambda_n+\beta}$, which implies that 
	\begin{equation}
	\kappa(\widetilde{\JJ}_k)<\kappa({\JJ}_k)
	\end{equation}
	whit this, the matrix $\widetilde{\JJ}_k$ is better conditioned than $\JJ_k$.

\end{rmrk}

\begin{rmrk}\label{rmrk:aprrox_jk_by_neumann_series}

Given that $\widetilde{\JJ}_k=\JJ_k+\beta\II$, then $\JJ_k^{-1}=\widetilde{\JJ}_k^{-1}(\II-\beta\widetilde{\JJ}_k^{-1})^{-1}$ and by using the Neumann series we have
\begin{eqnarray}
\JJ_k^{-1}&=&\widetilde{\JJ}_k^{-1}\sum_{m=0}^{\infty}(\beta\widetilde{\JJ}_k^{-1})^m\\
&=&\frac{1}{\beta}\sum_{m=1}^{\infty}(\beta\widetilde{\JJ}_k^{-1})^m
\end{eqnarray}

\end{rmrk}

\begin{rmrk}

	If $\bm{\alpha}$ is the true solution of the equation $\JJ_k\bm{\alpha}=U$ and $\widetilde{\bm{\alpha}}$ is the solution for the equation $\widetilde{\JJ}_k\widetilde{\bm{\alpha}}=U$ then, from the remark \ref{rmrk:aprrox_jk_by_neumann_series}, \begin{equation}\label{eq:alpha_series}
	\bm{\alpha}=\frac{1}{\beta}\sum_{m=1}^{\infty}(\beta\widetilde{\JJ}_k^{-1})^mU
	\end{equation} and
	\begin{equation}
	\bm{\alpha}=\sum_{m=0}^{\infty}(\beta\widetilde{\JJ}_k^{-1})^m\widetilde{\bm{\alpha}}.
	\end{equation}
\end{rmrk}

If we truncate the series in \eqref{eq:alpha_series} up to order $M$ we obtain an approximation of the true   solution $\bm{\alpha}$, we denote it by  
\begin{equation}\label{eq:approx_alpha_barM}
\widetilde{\bm{\alpha}}_M=\frac{1}{\beta}\sum_{m=1}^{M}(\beta\widetilde{\JJ}_k^{-1})^mU.
\end{equation}
With this, the error of the approximation $\widetilde{\bm{\alpha}}_M$ can be bounded by using the formula
\begin{eqnarray}
\bm{\alpha}-\widetilde{\bm{\alpha}}_M&=&\frac{1}{\beta}\sum_{m=M+1}^{\infty}(\beta\widetilde{\JJ}_k^{-1})^mU\\
&=&\frac{1}{\beta}(\II-\beta\widetilde{\JJ}_k^{-1})^{-1}(\beta\widetilde{\JJ}_k^{-1})^{M+1}U.
\end{eqnarray}
Taking the Euclidean norm, then we have in terms of the spectral norm,
\begin{equation}
\|\bm{\alpha}-\widetilde{\bm{\alpha}}_M\|_2\leq\| (\beta\widetilde{\JJ}_k^{-1})^{M+1}\|_2\|(\II-\beta\widetilde{\JJ}_k^{-1})^{-1} \|_2\|U\|_2.
\end{equation}
From the remark \ref{rmrk:convergence_neumann_series_jk} we can conclude that
\begin{equation}\label{eq:error_bound_MDI}
\|\bm{\alpha}-\widetilde{\bm{\alpha}}_M\|_2\leq\frac{\beta}{\lambda_n}\left(\frac{\beta}{\lambda_n+\beta}\right)^{M}\|U\|_2.
\end{equation}
In particular,
\begin{equation}\label{eq:error_bound_R}
	\|\bm{\alpha}-\widetilde{\bm{\alpha}}\|_2\leq\frac{\beta}{\lambda_n}\|U\|_2.
\end{equation}

An iterative proceeder to compute \eqref{eq:approx_alpha_barM}, with the better conditioned matrix $\beta\widetilde{\JJ}_k$, can be obtained just by noting that, with $\widetilde{\bm{\alpha}}=\widetilde{\JJ}_k^{-1}U$,
  \begin{eqnarray}\label{eq:approx_alpha_barM_iter}
  \widetilde{\bm{\alpha}}_M&=&\frac{1}{\beta}\sum_{m=1}^{M}(\beta\widetilde{\JJ}_k^{-1})^mU\\
  &=&\sum_{m=1}^{M}(\beta\widetilde{\JJ}_k^{-1})^{m-1}\widetilde{\bm{\alpha}}\\
  &=&\widetilde{\bm{\alpha}}+\beta\widetilde{\JJ}_k^{-1}\left(\widetilde{\bm{\alpha}}+\beta\widetilde{\JJ}_k^{-1}\left(\widetilde{\bm{\alpha}}+\cdots\right)\right).
  \end{eqnarray}
  Hence we can compute $\widetilde{\bm{\alpha}}_M$  as: 
  \begin{eqnarray}
  	\label{eq:iter_first_step}\widetilde{\bm{\alpha}}_0&=&\widetilde{\JJ}_k^{-1}U\\
  	 \label{eq:iter_m_step}\widetilde{\bm{\alpha}}_{m}&=&\widetilde{\bm{\alpha}}_0+\beta\widetilde{\JJ}^{-1}_k\widetilde{\bm{\alpha}}_{m-1}, \ \ \ \mbox{\ \  for } m=1,2,\ldots,M
  \end{eqnarray}
Since $0<\frac{\beta}{\lambda_n+\beta}<1$ the convergence order is controlled theoretically, however when $\lambda_n$ is near to the machine epsilon (e.g., in double precision the machine epsilon is approximately 2.22e-16) the ratio $\frac{\beta}{\lambda_n+\beta}$ is very close to 1 and in this case the convergence may be too slow, but with a small $M$ (up to $M=5$ it works) is enough to improve the error in \eqref{eq:error_bound_R}.

\section{Numerical tests}
For this section we have tried to reproduce examples and tests presented in \cite{Fang2017} due to the solver presented by authors is recent and very accurate, also is they show a detailed local description of solutions of Helmholtz equation in terms of its phases and ray direction local.
\subsection{Local truncation error}

For small stencils we take $10^7 \leq\kappa_0\leq10^{14}$, with $\kappa_0=10^{7+\sqrt{n}}$ where $n$ is the size of the stencil,  and we take $\beta=\frac{\lambda_1-\kappa_0  \lambda_n}{\kappa_0 -1}$ ensuring, from remark \ref{rmrk:better_coditioning_wideJk}, that $\kappa(\widetilde{\JJ}_k)\approx\kappa_0$, which is an adequate condition number to work in double precision. 

\begin{figure}
	\begin{center}	
		\begin{tabular}{|c|c|c|c|}
			\hline
			\hline
			\includegraphics[scale=0.5]{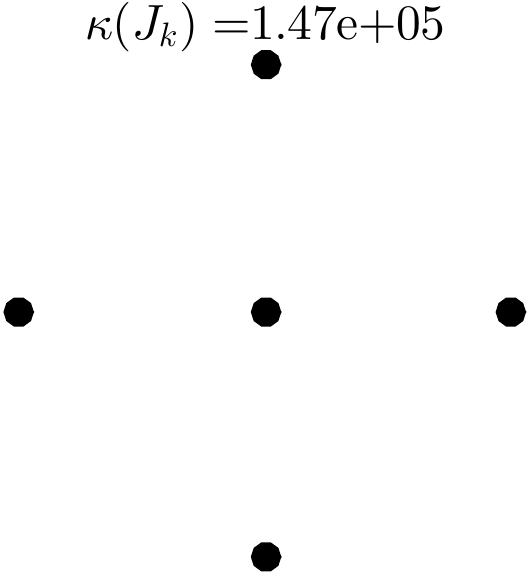}&\includegraphics[scale=0.5]{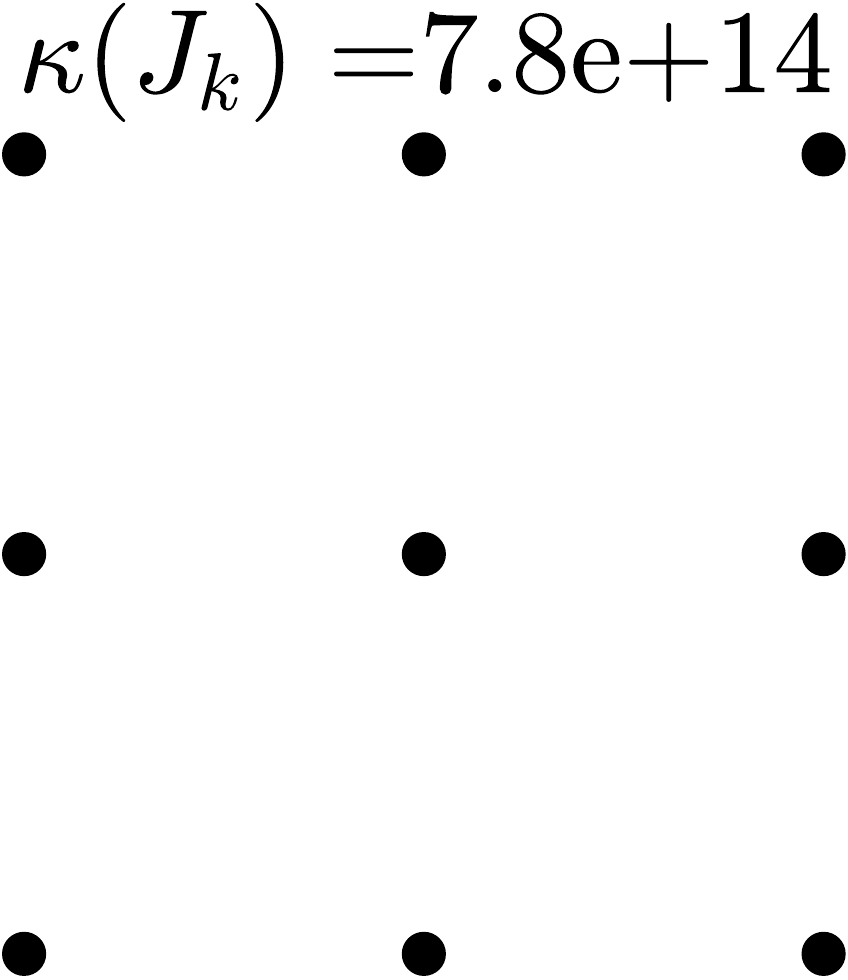} &\includegraphics[scale=0.5]{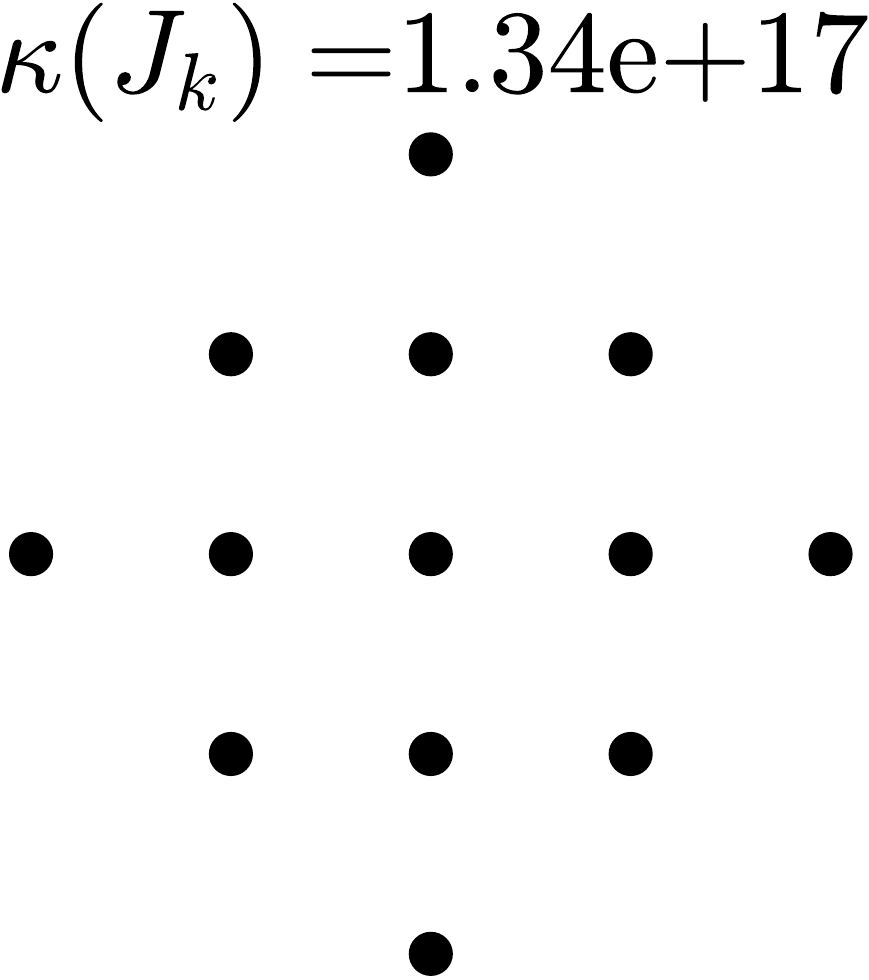} &\includegraphics[scale=0.5]{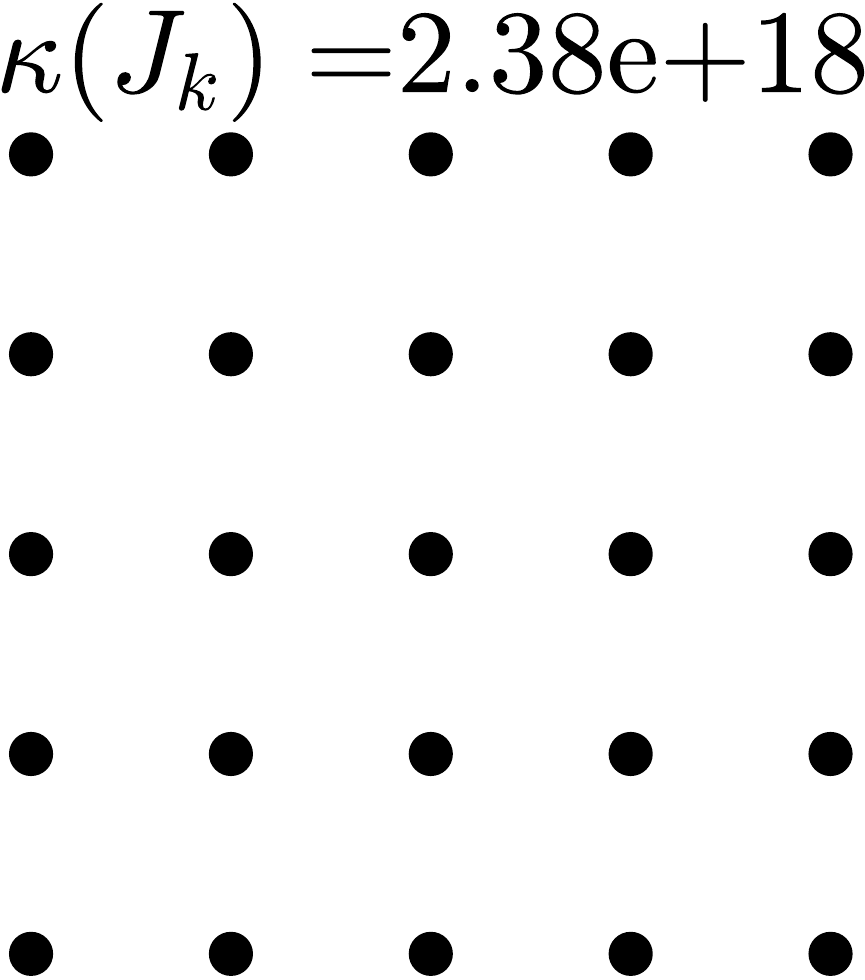}\\
			\hline
			\hline
			\includegraphics[scale=0.5]{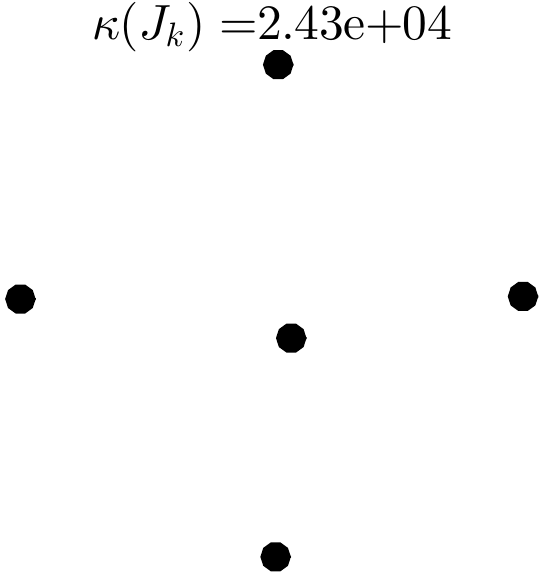}&\includegraphics[scale=0.5]{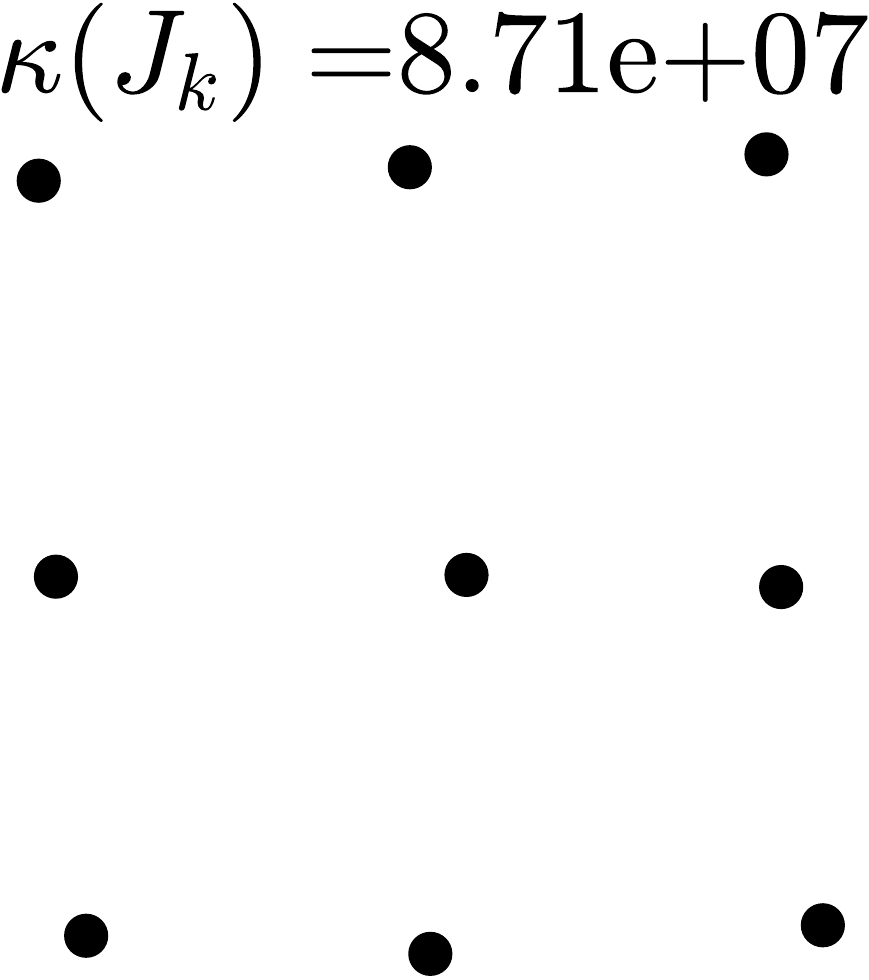} &\includegraphics[scale=0.5]{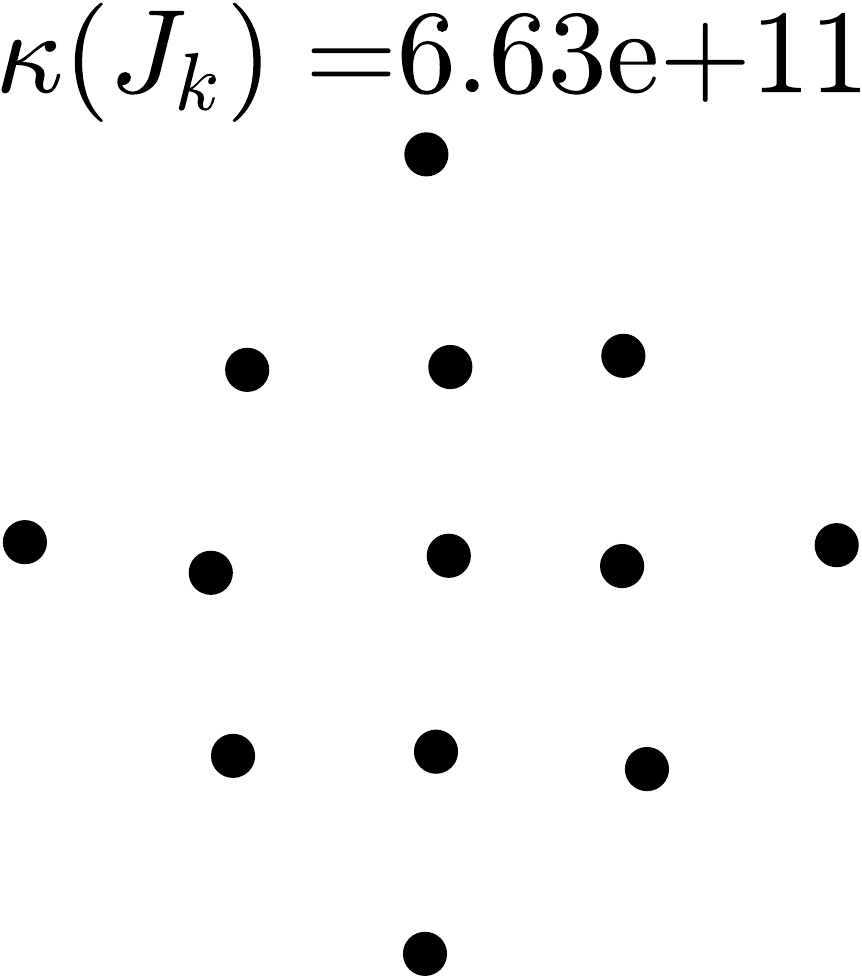} &\includegraphics[scale=0.5]{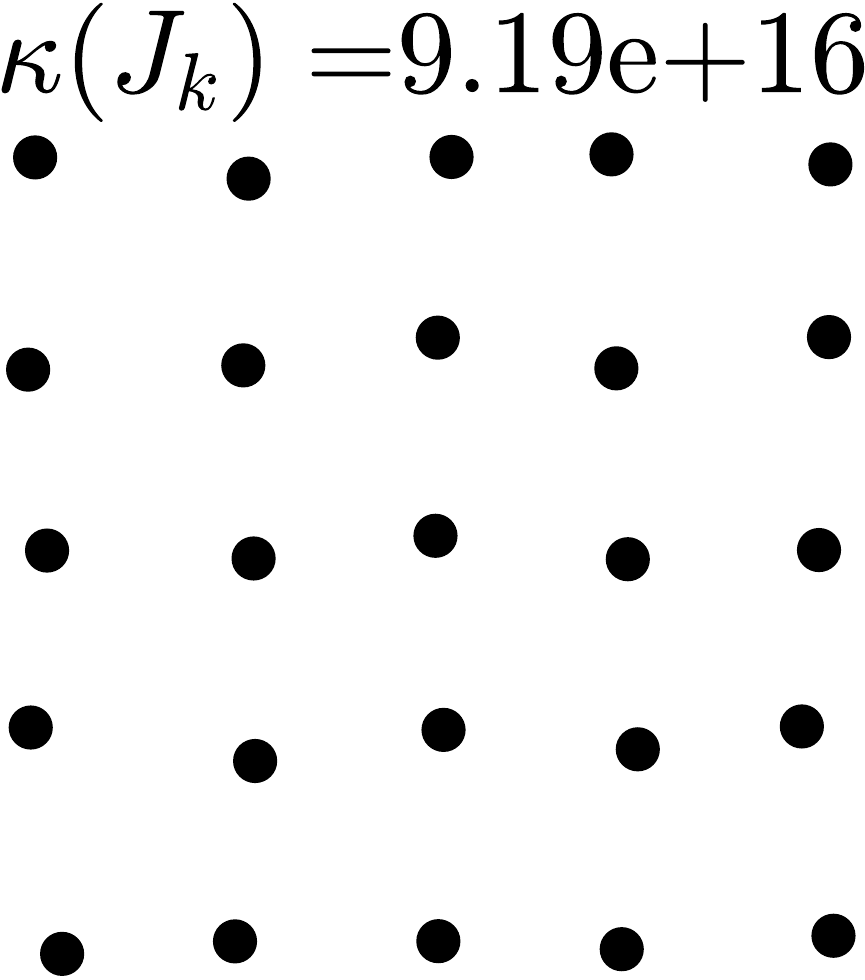}\\
			\hline
		\end{tabular}
	\end{center}
	\caption{Plots of some small stencils with the respective approximated condition number of $\JJ_k$, with $k=100$ and $h=\frac{2\pi}{6k}$. Top row: Stencils are taken from a regular square grid. Bottom row: perturbed position from stencils of top row.} 	\label{fig:stencils_square}
\end{figure}

\begin{figure}
	\begin{center}	
		\begin{tabular}{|c|c|c|c|c|}
			\hline
			\hline
			\includegraphics[scale=0.45]{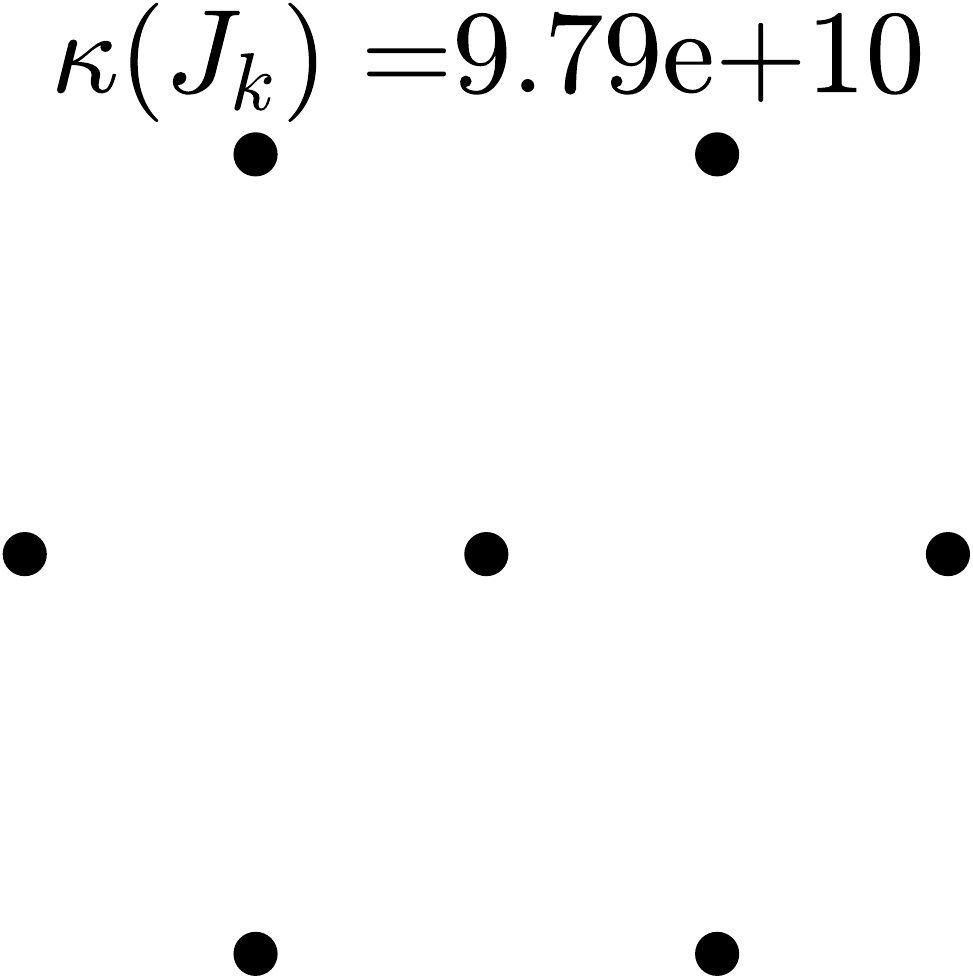}&\includegraphics[scale=0.45]{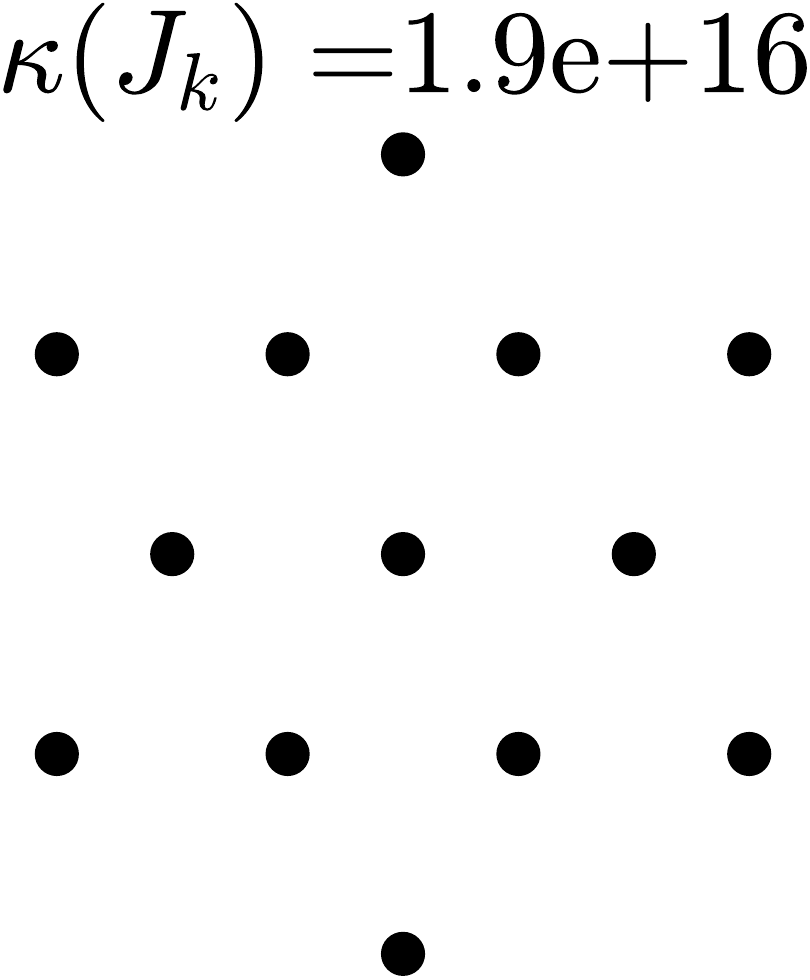} &\includegraphics[scale=0.45]{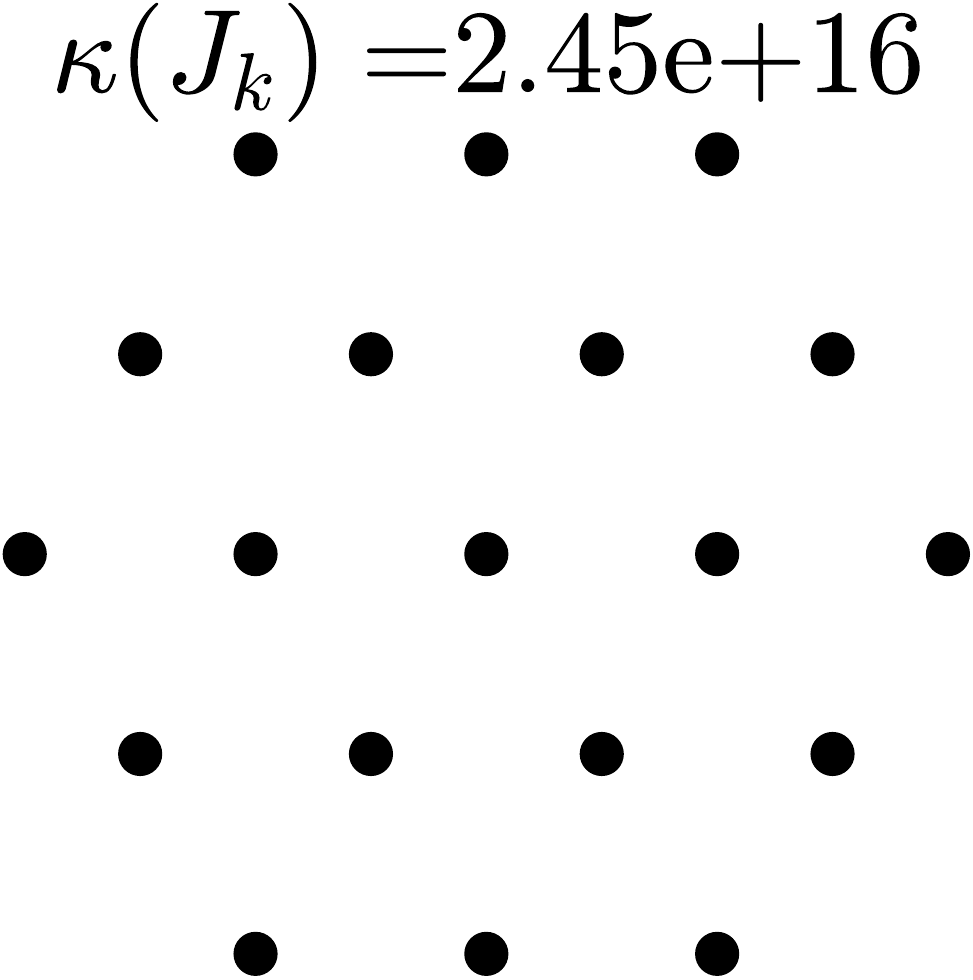} &\includegraphics[scale=0.45]{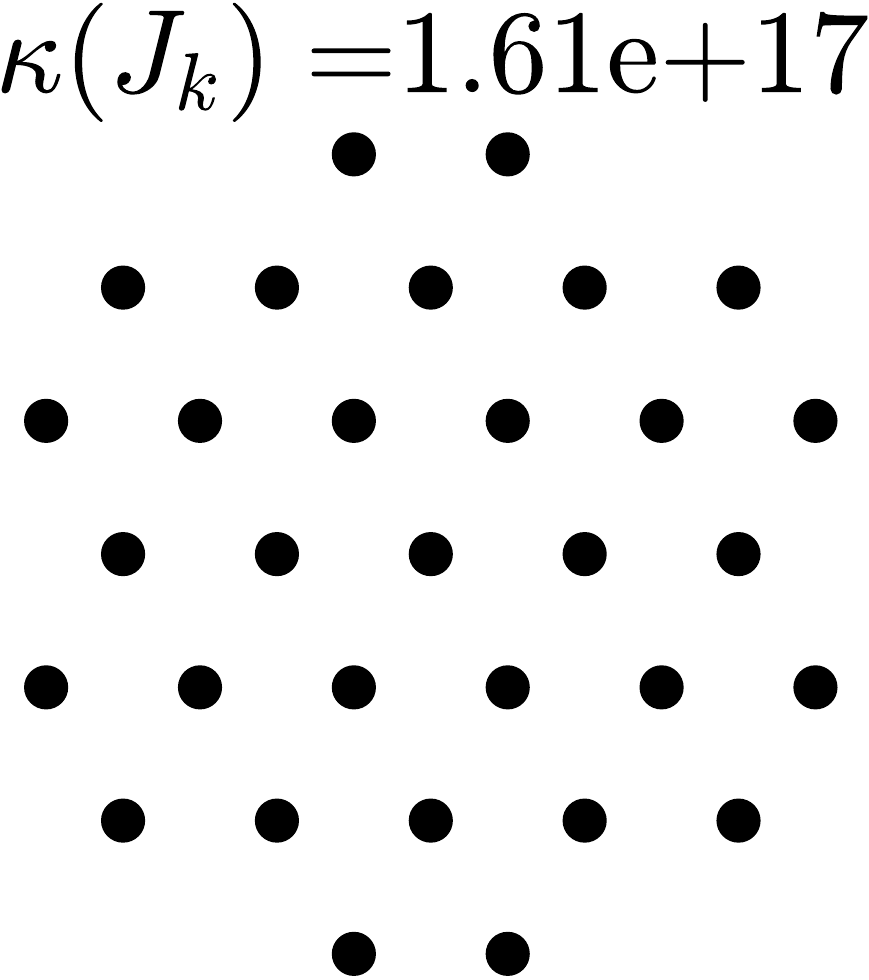}
			&\includegraphics[scale=0.45]{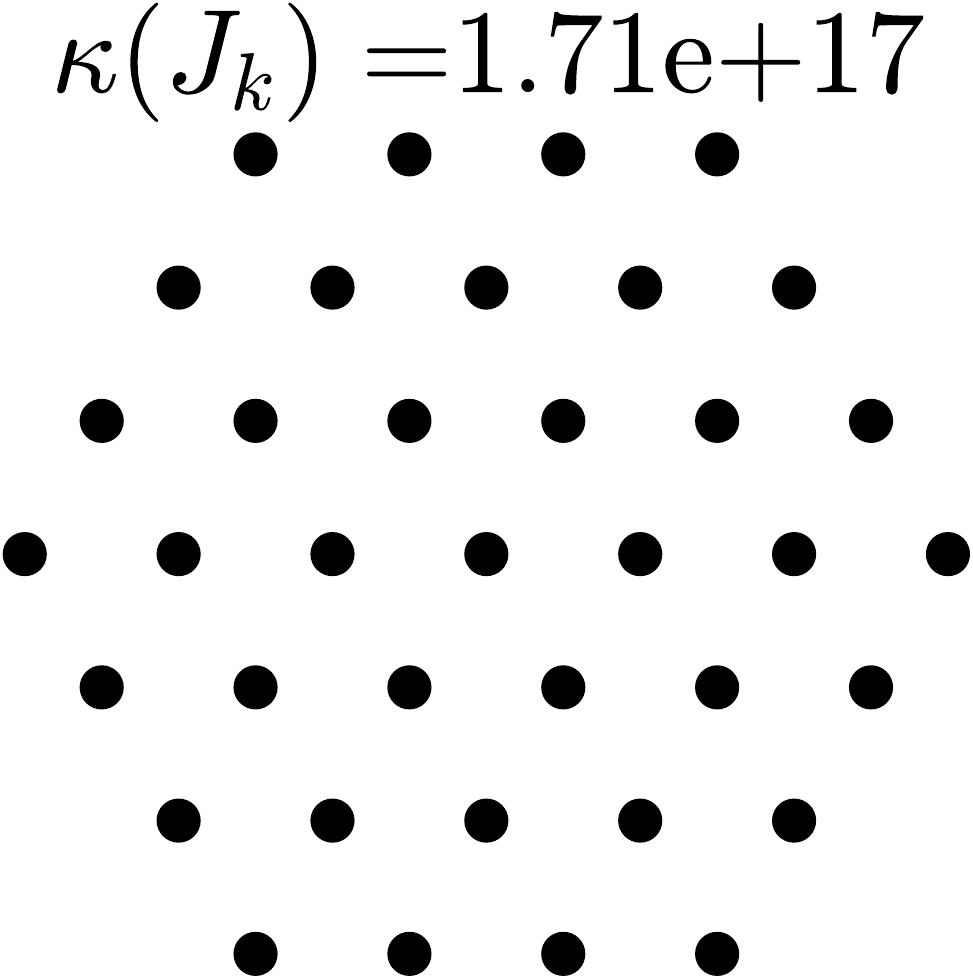}\\
			\hline
			\hline
			\includegraphics[scale=0.45]{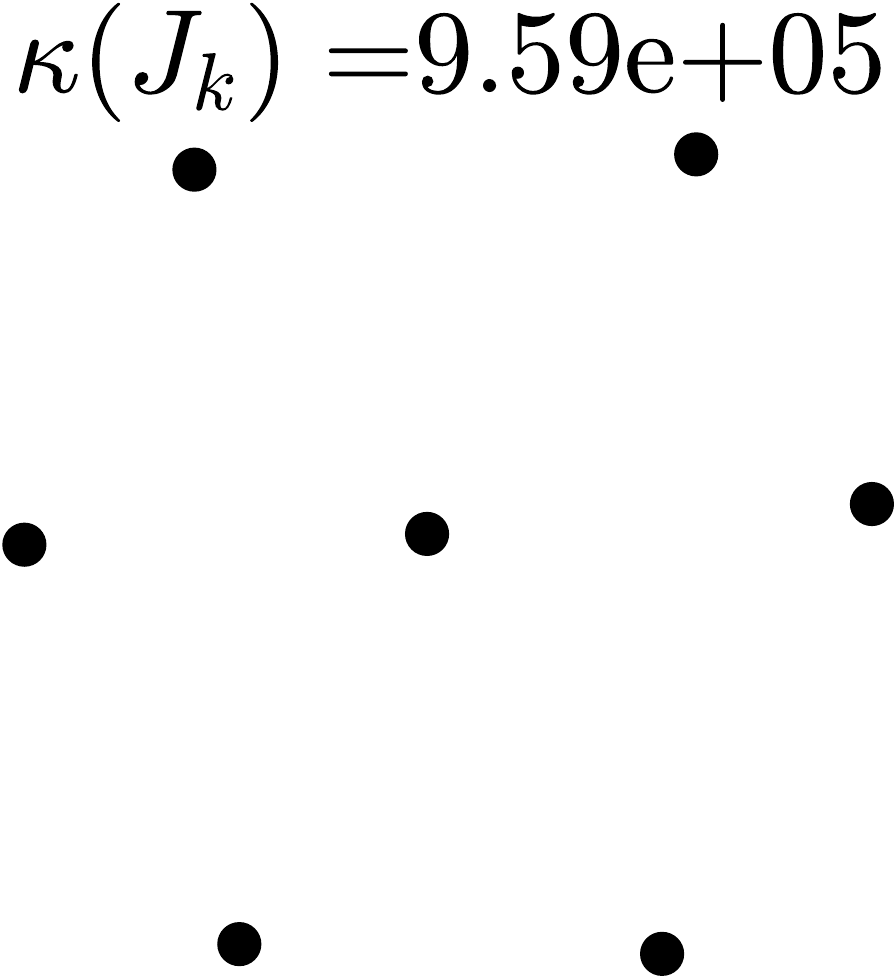}&\includegraphics[scale=0.45]{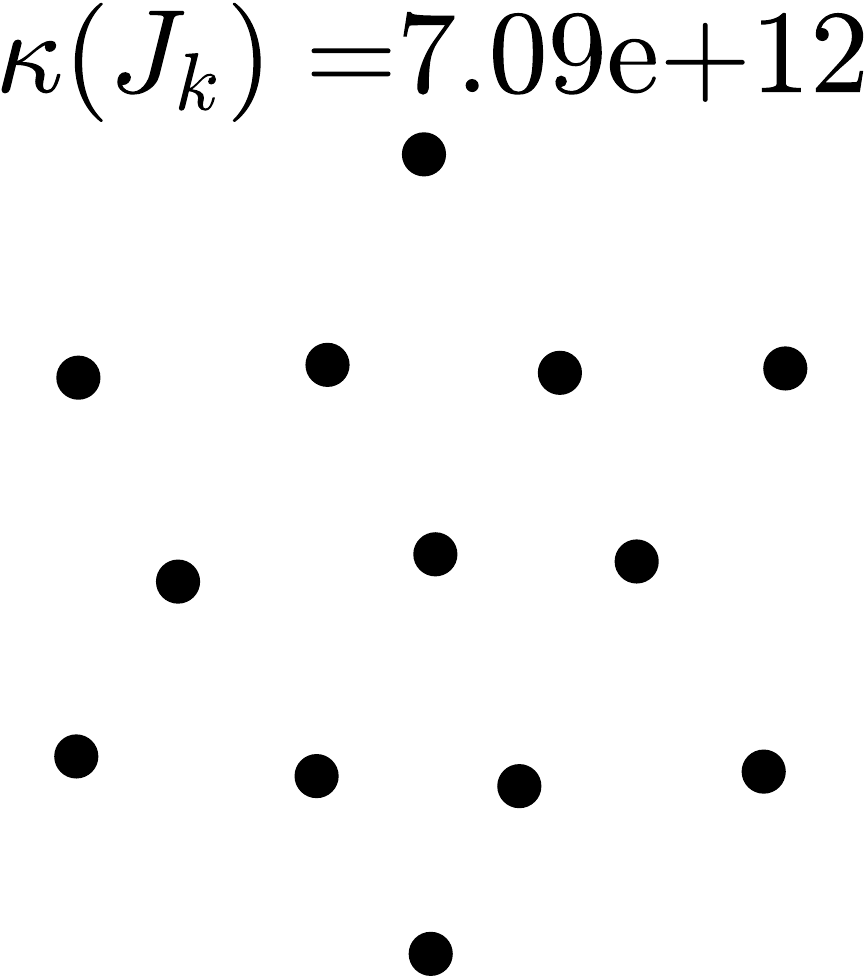} &\includegraphics[scale=0.45]{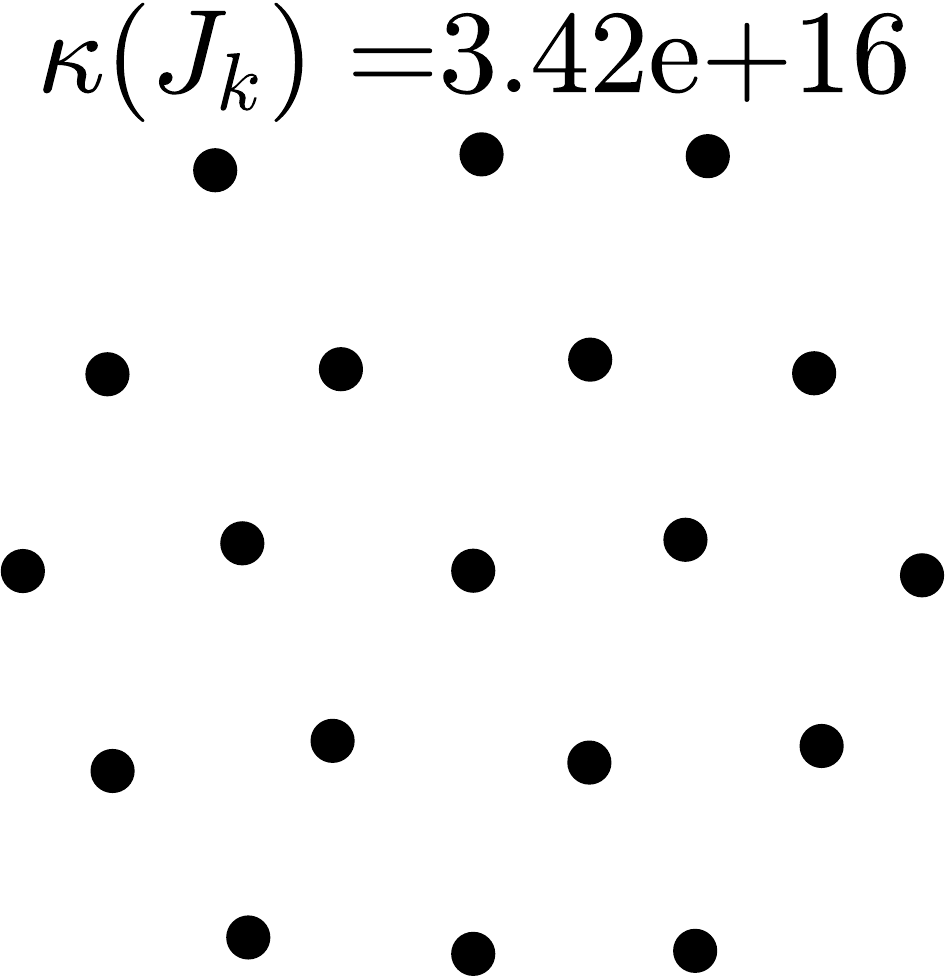} &\includegraphics[scale=0.45]{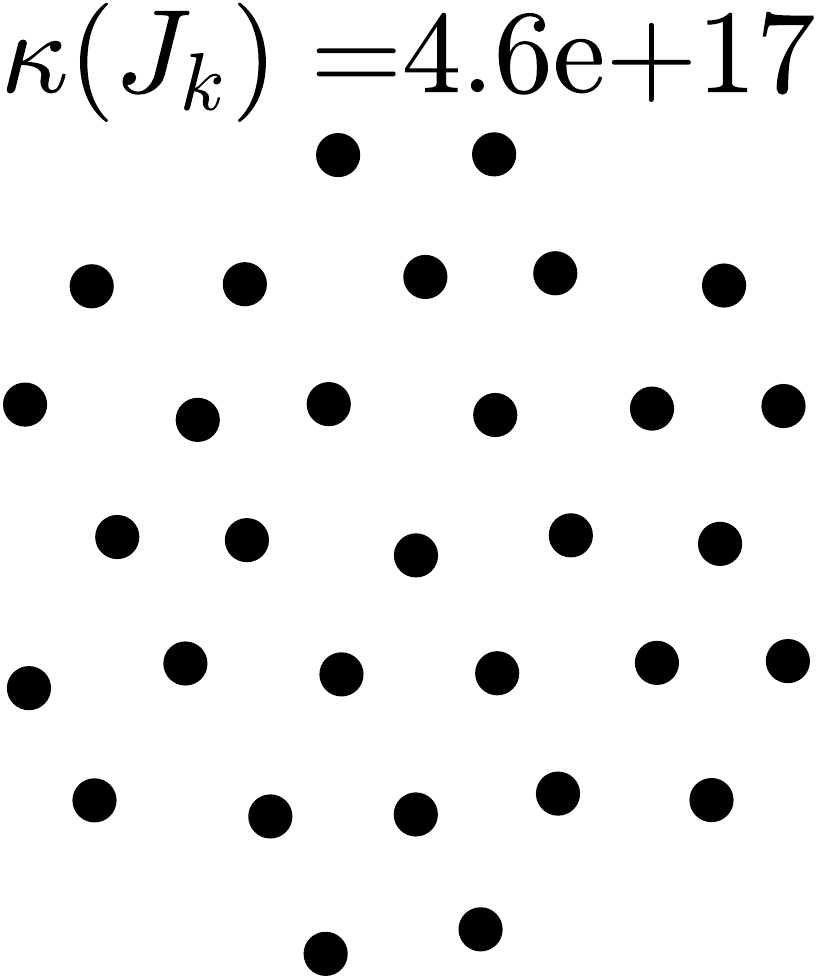}
			&\includegraphics[scale=0.45]{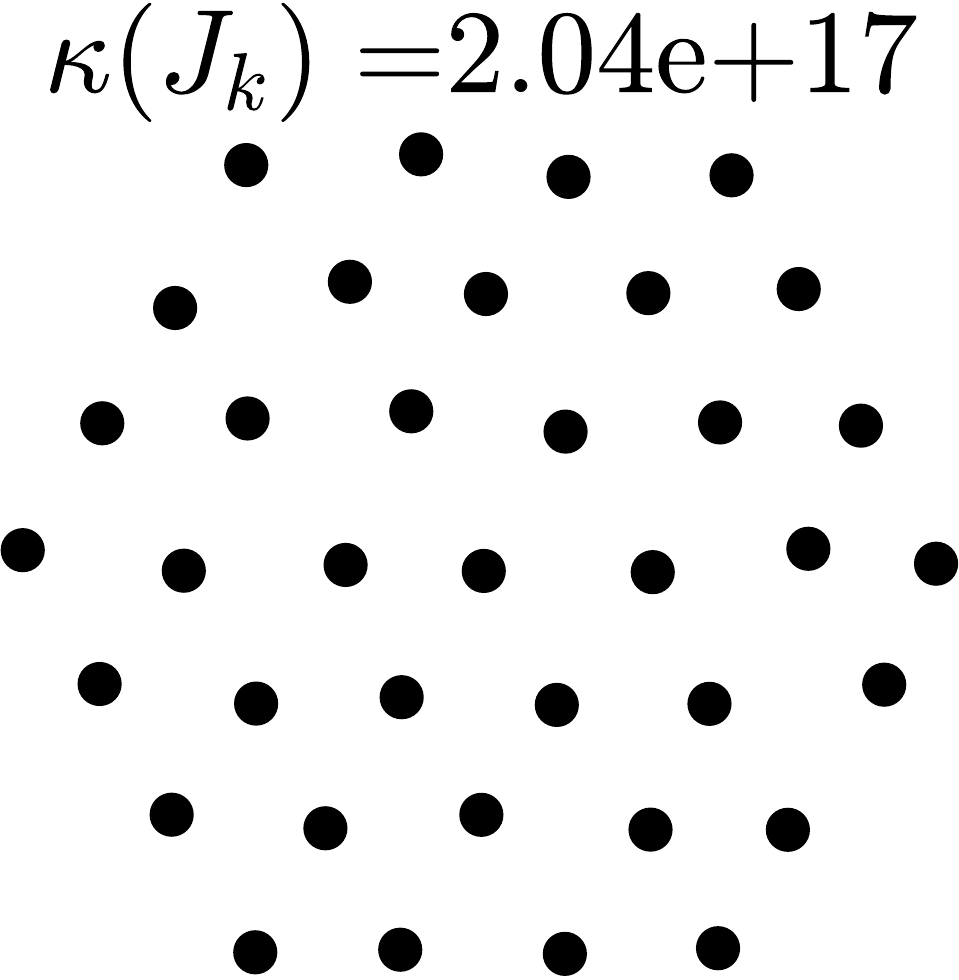}\\
			\hline
		\end{tabular}
	\end{center}
	\caption{Plots of some small stencils with the respective approximated condition number of $\JJ_k$, with $k=100$ and $h=\frac{2\pi}{6k}$. Top row: Stencils are taken from a regular hexagonal grid. Bottom row: perturbed position from stencils of top row.} \label{fig:stencils_hexagonal}
\end{figure}

We have noted empirically that the matrix $\JJ_k$ is worse conditioned for stencil with nodes put symmetrically on a regular grids, e.g. with square and hexagonal grids. See figures \ref{fig:stencils_square} and \ref{fig:stencils_hexagonal}, where we can observe that first size of severe ill-conditioning occur with symmetric stencils of $13$ nodes. However, in this case, with a small perturbation in the position of the nodes, its associated interpolation matrix $\JJ_k$ has a better condition number.

For numerical tests we consider the Helmholtz solutions

\begin{equation}\label{eq:sol_imped1}
u_1(x,y)=\sqrt{k}H^{(1)}_0(k\sqrt{(x-2)^2+(y-2)^2}),
\end{equation}
and
\begin{equation}\label{eq:sol_imped2}
\begin{array}{rcl}
u_2(x,y)&=&\sqrt{k}H^{(1)}_0(k\sqrt{(x+20)^2+(y+20)^2})+2\sqrt{k}H^{(1)}_0(k\sqrt{(x-20)^2+(y-20)^2})\\
&+&0.5\sqrt{k}H^{(1)}_0(k\sqrt{(x+20)^2+(y-20)^2})-\sqrt{k}H^{(1)}_0(k\sqrt{(x-20)^2+(y+20)^2}),
\end{array}
\end{equation}
where $H_0^{(1)}$ is the Hankel function of the first kind and $k$ is the constant wavenumber. $u_1$ corresponds to solution of the single source problem located at $\x_s=(2,2)$ whereas $u_2$ corresponds to solution of the problem with four single sources located at $\x_{s_1}=(-20,-20)$, $\x_{s_2}=(20,20)$, $\x_{s_3}=(-20,20)$ and $\x_{s_1}=(20,-20)$.
In the first test, we consider the exact solution to the Helmholtz equation \eqref{eq:sol_imped2} to get the local truncation error of the approximation $\Delta_{S,k}$, also we validate the ill-conditioning of the matrix $\JJ_k$ and the better conditioning of $\widetilde{\JJ}_k$. We compare 
to interpolate data produced with the function  $u(x,y)=\dfrac{\partial u_2}{\partial x}(x,y)$, at stencils nodes centered at $(x_c,y_c)=(0.5,0.5)$. 
With stencils $S$ as in Fig. \ref{fig:stencils_square} and Fig. \ref{fig:stencils_hexagonal} we solve the system $\JJ_k\bm{\alpha}=U$ with $U=u|_S$. In tables \ref{table:Jk_bar_square} and \ref{table:Jk_bar_hexagonal} it can be seen how is the behavior of the condition number $\kappa(\JJ_k)$ and $\kappa(\widetilde{\JJ}_k)$ for symmetric nodes distributions respectively. We have compared the relative errors of calculate approximations of $\alpha$: with direct method, the MDI and MDIIT, whit error: Error D, Error RD and Error IT as in tables \ref{table:Jk_bar_square} and \ref{table:Jk_bar_hexagonal}.
\begin{table}[ht]
	\begin{center}
		\begin{tabular}{c|cc|ccc}
			\hline
			\rowcolor{gray!20} Stencil size ($n$)  &  $\kappa(\JJ_k)$ & $\kappa(\widetilde{\JJ}_k)$ & Error D & Error RD  & Error IT \\
			\hline\hline 
	        5 & 1.47e+05 & 1.72e+09 & 3.22e-13 & 8.78e+04 & 3.22e-13\\
	        9 & 4.68e+12 & 1.00e+10 & 1.18e-08 & 1.60e-03 & 1.57e-03\\
	        13 & 6.75e+15 & 4.03e+10 & 5.67e-08 & 1.31e-03 & 1.48e-07\\
	        25 & 9.42e+16 & 1.00e+12 & 3.44e-06 & 3.34e-06 & 2.79e-06\\
	        37 & 9.59e+16 & 1.21e+13 & 1.80e-05 & 1.43e-04 & 7.26e-08\\ \hline
		\end{tabular}
		\\
	\end{center} 
	\caption{This table shows values of the condition number of the matrices $\JJ_k$ and  $\widetilde{\JJ}_k$, and local truncation errors of the approximation $\Delta_{S,k}u(\x)\approx\Delta u(\x)$. We have used stencils as in Fig. \ref{fig:stencils_square}. We compare errors $|\Delta_{S,k}u(\x)-\Delta u(\x)|$ produced by using Cholesky (Error D), by \eqref{eq:iter_first_step} (Error RD) and   \eqref{eq:iter_m_step} (Error IT) in using $5$ iterations applied to solution. The used function is $u=u_2$ as in \eqref{eq:sol_imped2}.}\label{table:Jk_bar_square}	
\end{table}

\begin{table}[ht]
	\begin{center}
		\begin{tabular}{c|cc|ccc}
			\rowcolor{gray!20} Stencil size ($n$)  &  $\kappa(\JJ_k)$ & $\kappa(\widetilde{\JJ}_k)$ & Error D & Error RD  & Error IT \\
			\hline\hline 
	   	7 & 2.06e+09 & 4.42e+09 & 3.48e-10 & 2.19e-02 & 5.06e-02\\
		13 & 6.09e+15 & 4.03e+10 & 6.26e-09 & 4.15e-06 & 9.01e-08\\
		19 & 3.57e+16 & 2.29e+11 & 1.12e-06 & 2.75e-06 & 2.77e-06\\
		31 & 2.03e+17 & 3.70e+12 & 5.47e-06 & 8.83e-06 & 7.71e-06\\
		37 & 1.18e+17 & 1.21e+13 & 3.99e-07 & 4.33e-05 & 3.40e-05\\ \hline
		\end{tabular}
		\\
	\end{center} 
		\caption{This table shows values of the condition number of the matrices $\JJ_k$ and  $\widetilde{\JJ}_k$, and local truncation errors of the approximation $\Delta_{S,k}u(\x)\approx\Delta u(\x)$. We have used stencils as in Fig. \ref{fig:stencils_hexagonal}. We compare errors $|\Delta_{S,k}u(\x)-\Delta u(\x)|$ produced by using Cholesky (Error D), by \eqref{eq:iter_first_step} (Error RD) and   \eqref{eq:iter_m_step} (Error IT) in using $5$ iterations applied to solution. The used function is $u=u_2$ as in \eqref{eq:sol_imped2}.}\label{table:Jk_bar_hexagonal}		
\end{table}

\subsection{Pollution-effect and convergence}

\subsubsection{Test 1}

In this test we calculate the approximated solution for the problem

    \begin{equation}\label{eq:helmholtz_impedance}
    \left\{
    \begin{array}{rcll}
    -\Delta u(\x)-\omega^2u(\x)&=&0, &\mbox{ in } \Omega\\
    \frac{\partial}{\partial\n}u(\x)+i\omega u(\x)&=&g(\x), & \mbox{ on } \Gamma=\partial\Omega
    \end{array}
    \right.
    \end{equation}
     with the known data $g(\x)$, $\Omega=(-0.5,0.5)\times(-0.5,0.5)$.  Results are verified with solutions $u_1$ and $u_2$ in \eqref{eq:sol_imped1} and \eqref{eq:sol_imped2}. 
     
     We see in tables \ref{table:fang1square_pollution} and \ref{table:fang1hexagon_pollution} how it behaves the error when the resolution is kept constant with $N_g=6$ nodes per wavelength. We have fixed the condition number of the local interpolation matrices, $\widetilde{\JJ}_k$, to $\kappa_0=10^{(7+\sqrt{n})}$. In both cases, for uniform square and hexagonal grids, we see that the errors are $\bO(\omega)$, but with a small positive slope and in some cases, increasing the stencil size, the error fits to $m\omega+b$ with $m$ small and negative, see table \ref{table:fang1hexagon19_pollution}. Although pollution-effects does not eliminated this is kept relatively constant, in accord to order error and does not increase dramatically.

     Results of convergence tests are summarized in tables \ref{table:fang1convergence_square} and \ref{table:fang1convergence_hexa}. For the local interpolation matrices $\widetilde{\JJ}_k$ we have chosen, empirically, $\beta$ such that the condition number  $\kappa_0=10^{(7+(hk)^{-1}+\sqrt{n})}$.

\begin{table}
	\begin{center}
		\begin{tabular}{cccc|cc}
			\rowcolor{gray!20} $\frac{\omega}{2\pi}$  &  $\frac{1}{h}$ & $N$ & $\kappa(\mathbf{H})$ & $\|u_1-\widetilde{u}_1\|_{\infty}$ & $\|u_2-\widetilde{u}_2\|_{\infty}$   \\
			\hline\hline 
	    	10 & 60 & 3721 & 1.36e+04 & 1.76e-04 & 1.72e-04  \\  
	    	20 & 120 & 14641 & 7.46e+04 & 2.54e-04  & 1.71e-04   \\
			40 & 240 & 58081 & 4.19e+05 & 5.15e-04 & 2.82e-04    \\
			80 & 480 & 231361 & 2.34e+06 & 1.04e-03 & 5.90e-04\\
			120 & 720 & 519841 & 5.25e+06 & 1.57e-03 & 8.65e-04  \\
			\hline
		\end{tabular}
		\\
	\end{center} 
	\caption{Results for approximated solutions of \eqref{eq:helmholtz_impedance}. This approximation was calculated from a square uniform grid of $\Omega\cap\partial\Omega$. For inner nodes the stencil size is $n=9$, at boundary nodes $n=15$, the number of nodes per wave length is kept constant with $N_g=6$. Solution for the system $-\mathbf{H}\widetilde{U}=G$ is obtained for LU factorization. Errors scale respect to the frequency as $\|u_1-\widetilde{u}_1\|_{\infty}\sim(1.3e-05)(\omega/2\pi)+1.7e-05$ and $\|u_2-\widetilde{u}_2\|_{\infty}\sim(6.6e-06)(\omega/2\pi)+5.8e-05$ } \label{table:fang1square_pollution}		
\end{table}

\begin{table}
	\begin{center}
		\begin{tabular}{cccc|cc}
			\rowcolor{gray!20} $\frac{\omega}{2\pi}$  &  $\frac{1}{h}$ & $N$ & $\kappa(\mathbf{H})$ & $\|u_1-\widetilde{u}_1\|_{\infty}$ & $\|u_2-\widetilde{u}_2\|_{\infty}$   \\
			\hline\hline 
			10 & 60  & 4237   & 1.96e+04 & 3.95e-05 & 3.30e-05  \\  
	    	20 & 120 & 16752  & 1.09e+05 & 4.49e-05 & 3.35e-05  \\
	     	40 & 240 & 66861  & 5.25e+05 & 9.30e-05 & 3.59e-05   \\
			80 & 480 & 266680 & 3.46e+06 & 6.69e-05 & 5.62e-05\\
			120 & 720 & 599458 & 8.21e+06 & 9.05e-05 & 4.25e-05 \\
		 \hline
		\end{tabular}
		\\
	\end{center} 
	\caption{For inner points the stencil size is $n=13$ with an hexagonal uniform grid, at boundary points the stencil size is $n_b=25$, the number of nodes per wave length is kept constant with $N_g=6$. Solution for the system $-\mathbf{H}\widetilde{U}=G$ is obtained for LU factorization. Error scale respect to the frequency as $\|u_1-\widetilde{u}_1\|_{\infty}\sim(3.9e-07)(\omega/2\pi)+4.7e-05$ and $\|u_2-\widetilde{u}_2\|_{\infty}\sim(1.4e-07)(\omega/2\pi)+3.3e-05$ } \label{table:fang1hexagon_pollution}		
\end{table}

\begin{table}
	\begin{center}
		\begin{tabular}{cccc|cc}
			\rowcolor{gray!20} $\frac{\omega}{2\pi}$  &  $\frac{1}{h}$ & $N$ & $\kappa(\mathbf{H})$ & $\|u_1-\widetilde{u}_1\|_{\infty}$ & $\|u_2-\widetilde{u}_2\|_{\infty}$   \\
			\hline\hline 
		    10 & 60 & 4237     & 8.18e+04 & 1.97e-05 & 1.23e-05   \\  
			20 & 120 & 16752   & 8.58e+05 & 2.22e-05 & 1.77e-05   \\
			40 & 240 & 66861   & 6.12e+05 & 1.94e-05 & 1.47e-05   \\
			80 & 480 & 266680  & 5.22e+06 & 1.78e-05 & 1.28e-05   \\
			120 & 720 & 599458 & 9.53e+06 & 1.85e-05  & 9.86e-06  \\
			\hline
		\end{tabular}
		\\
	\end{center} 
	\caption{Results for approximated solutions of \eqref{eq:helmholtz_impedance}. This approximation was calculated from a square uniform grid of $\Omega\cap\partial\Omega$. For inner nodes the stencil size is $n=19$, at boundary nodes $n=25$, the number of nodes per wave length is kept constant with $N_g=6$. Solution for the system $-\mathbf{H}\widetilde{U}=G$ is obtained for LU factorization. Errors scale respect to the frequency as $\|u_1-\widetilde{u}_1\|_{\infty}\sim( -2.5e-08)(\omega/2\pi)+2.1e-05$ and $\|u_2-\widetilde{u}_2\|_{\infty}\sim(-4.4e-08)(\omega/2\pi)+ 1.6e-05$ } \label{table:fang1hexagon19_pollution}		
\end{table}

\begin{table}
	\begin{center}
		\begin{tabular}{cccc|cc}
			\rowcolor{gray!20} NPW=$N_g$  &  $\frac{1}{h}$ & $N$ & $\kappa(\mathbf{H})$ & $\|u_1-\widetilde{u}_1\|_{\infty}$ & $\|u_2-\widetilde{u}_2\|_{\infty}$   \\
			\hline\hline 
		6.0 & 120.0 & 14641 & 7.80e+04 & 2.54e-04 & 1.70e-04 \\
		8.6 & 171.4 & 29584 & 1.05e+05 & 2.67e-05 & 1.56e-05 \\
		12.2 & 244.9 & 60025 & 1.47e+05 & 3.07e-06  & 1.79e-06 \\
		17.5 & 350.0 & 122500 & 1.99e+05 &3.55e-07  & 2.06e-07 \\
		25.0 & 500.0 & 251001 & 2.99e+05& 4.44e-08  & 2.55e-08 \\
			\hline
		\end{tabular}
		\\
	\end{center} 
	\caption{With $\omega/2\pi=20$, stencil size for inner nodes $n=9$ and stencil size for boundary nodes $n_b=15$ taken in a square uniform grid. Errors scale respect to $h$  as $\|u_1-\widetilde{u}_1\|_{\infty}\sim( 1.74e+09) h^{6.2}$ and  $\|u_2-\widetilde{u}_2\|_{\infty}\sim (3.2e+09) h^{6.4}$.} \label{table:fang1convergence_square}		
\end{table}

\begin{table}
	\begin{center}
		\begin{tabular}{cccc|cc}
			\rowcolor{gray!20} NPW=$N_g$  &  $\frac{1}{h}$ & $N$ & $\kappa(\mathbf{H})$ & $\|u_1-\widetilde{u}_1\|_{\infty}$ & $\|u_2-\widetilde{u}_2\|_{\infty}$   \\
			\hline\hline 
	    	6.0 & 120.0 & 16752 & 3.39e+17 & 4.43e-04 & 3.45e-04   \\
	    	8.6 & 171.4 & 33960 & 2.29e+15 & 2.07e-04  & 1.01e-04  \\
			12.2 & 244.9 & 69338 & 2.10e+07 & 8.76e-05 & 7.58e-05  \\
			17.5 & 350.0 & 141753 & 2.65e+06 & 7.67e-06 & 1.04e-06 \\
	    	25.0 & 500.0 & 289291 & 2.14e+20 & 8.50e-07 & 3.82e-07   \\
			\hline
		\end{tabular}
		\\
	\end{center} 
	\caption{With $\omega/2\pi=20$, stencil size for inner nodes $n=13$ and stencil size for boundary nodes $n_b=19$, taken in an hexagonal uniform grid. The errors scale respect to $h$  as $\|u_1-\widetilde{u}_1\|_{\infty}\sim(1613 ) h^{3.1}$ and  $\|u_2-\widetilde{u}_2\|_{\infty}\sim (5.5e+04) h^{3.9}$.} \label{table:fang1convergence_hexa}		
\end{table}

\subsubsection{Test 2}

\begin{figure}[ht]
	\begin{center}
		\begin{tabular}{cc}			
			\includegraphics[width=7.7cm]{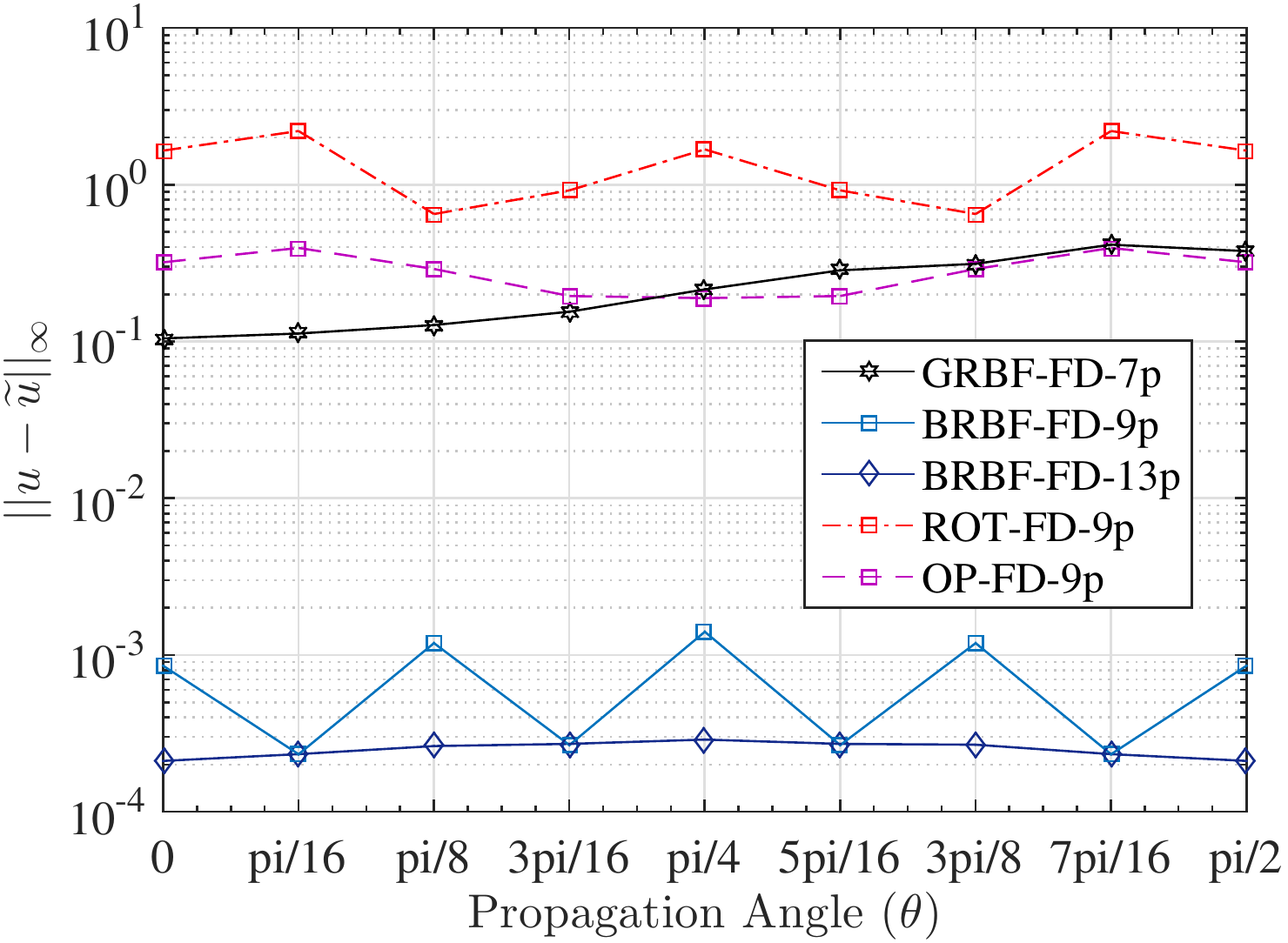} &          \includegraphics[width=7.7cm]{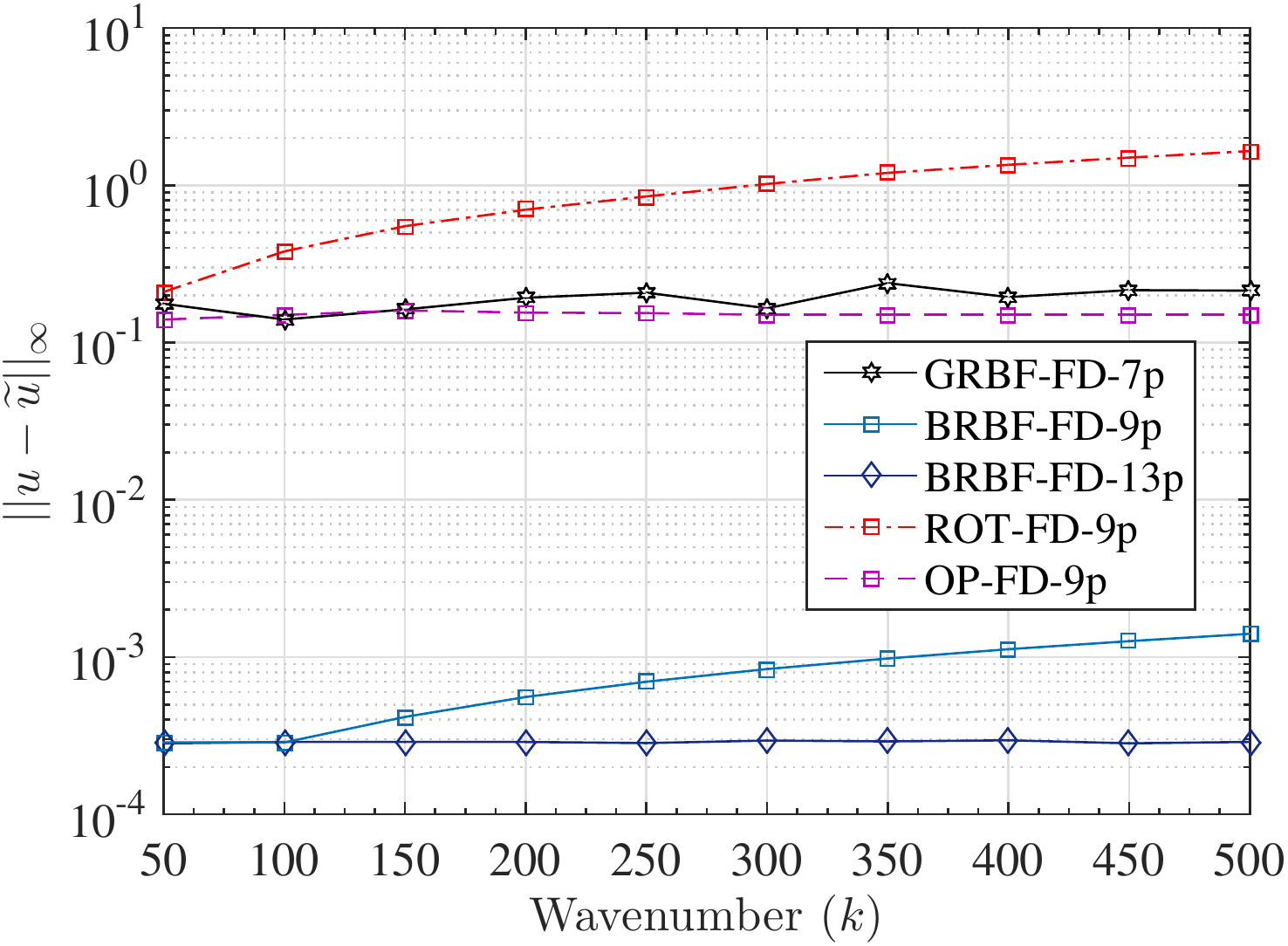}  			
		\end{tabular}       	
		\caption{Comparison of results between BRBF-FD-9p, BRBF-FD-13p,  GRBF-FD-7p and those reported in \cite{chen9}. (Left) Results for $k=500$ and $h=1/500$ varying the propagation angle. (Right) With $k$ varying, $\theta=\pi/4$ and $h=1/k$.}
		\label{fig:testpoll}
	\end{center}
\end{figure}

In this example we consider the problem

\begin{equation}
\begin{cases}\label{eq:abc_problem}
-\Delta u(x,y)-k^2u(x,y)=0& \mbox{ in } \Omega\\
\frac{\partial u}{\partial \mathbf{n}}u(x,y)+iku(x,y)=g(x,y)&  \mbox{ on } \partial\Omega,
\end{cases}
\end{equation}
with $\Omega=(0,1)\times(0,1)$, which has analytic solution given by the plane wave $u(x,y;k,\theta)=e^{ik(x\cos\theta+y\sin\theta)}$  when the data $g$ in the impedance boundary condition is given by
\begin{equation}
g(x,y)=
\begin{cases}
i(k-k_2)e^{ik_1x} & \mbox{ if } x\in (0,1) \mbox{ and } y=0\\
i(k+k_1)e^{i(k_1+k_2y)}& \mbox{ if } x=1 \mbox{ and } y\in(0,1)\\
i(k+k_2)e^{i(k_1x+k_2)}& \mbox{ if } x\in (0,1) \mbox{ and } y=1\\
i(k-k_1)e^{ik_2y}& \mbox{ if }  x=0 \mbox{ and } y\in(0,1)
\end{cases}
\end{equation}
with $k_1=k\cos\theta$ and $k_2=k\sin\theta$. This example is a standard benchmark within to test the numerical dispersion of solvers for Helmholtz equation. The approximated solutions for this problem was calculated using Gaussian RBF-FD on hexagonal grids with 7-stencil (GRBF-FD-7p), Bessel RBF-FD with 9-stencils (BRBF-FD-9p) and 13-stencils (BRBF-FD-13p) on uniform Cartesian grids and we compare with results reported in \cite{chen9}. In all methods it has fixed $\mbox{NPW}=2\pi$ to keep constant resolution. For GRBF-FD-7p  we have used the approximations in \eqref{eq:wlap} with shape parameter $\ve_{op}$  given by \eqref{eq:optimal_shape_p} to approximate the Laplace operator and all partial derivative operators involved in the boundary condition. To solve the local interpolations in BRBF-FD-9p and BRBF-FD-13p we have used the condition number $\kappa_0=10^{7+1/hk+\sqrt{n_s}}$.  Our results can be seen in Fig. \ref{fig:testpoll}. On left plot we can see that there is less anisotropy in the error of BRBF-FD-13p. On right plot we can see that the dispersion and pollution effects are mitigated. Besides, we can see that the behavior of the error in the case GRBF-FD-7p is similar to reported in \cite{chen9}; in cases BRBF-FD-9p and BRBF-FD-13p the error it improvements at least with two accuracy orders.  We can stand out results for Bessel-RBF-13p, since even in increasing the wavenumber by fixing the resolution $kh=1$, the error remains constant. We can say about GRBF-FD-7p that although its accuracy is rather modest compared to BRBF-FD, its low complexity and easy implementation it makes a feasible option.  

The aim of this chapter is to show approximated solutions for some Helmholtz problems to test the performance of  Gaussian RBF-FD and Bessel RBF-FD schemes. We have compared with the fundamental solution of Helmholtz equation and also, we have computed solutions with smooth and non-continuous coefficients to have a qualitative sight. At last section we show results of inverted data with FWI,  applied in three well-known benchmarks: Marmousi, 2D Oversthrust and 2004 BP velocity model. 

All codes were typed in Matlab R2016a and run in a laptop \textcolor{black}{ with Core i7 processor at 2.8 Ghz with 12 Gbytes of RAM.}
\section{Numerical results for some Helmholtz problems}\label{sec:numerical_tests}

In this section we test the BRBF-FD scheme to computing solutions with second and third order ABC.

\subsection{Approximated fundamental solutions}

It is known that the problem $-\Delta u-k^2u=\delta(\x)$ in the free-space $\R^2$ has a unique solution when it imposes the the Sommerfeld radiation condition
\begin{equation}
\lim_{\|\x\|\rightarrow\infty}\|\x\|^\frac{1}{2}\left(\frac{\partial}{\partial r}-ik\right)u(\x)=0.
\end{equation}
Particularly, the associated Green's function, which is solution of $-\Delta u(\x)-k^2u(\x)=\delta(\x-\x_0)$,  \cite{Ciraolo2009_radiation_condition_for_the_2D_helmholtz_equation}  is given by $u(\x)=G(\x,\x_0)$, with
\begin{equation}
G(\x,\x_0)=\frac{i}{4}H_0^{(1)}(k\|\x-\x_0\|)
	\end{equation}

\begin{figure}[ht]
	\includegraphics[scale=0.8]{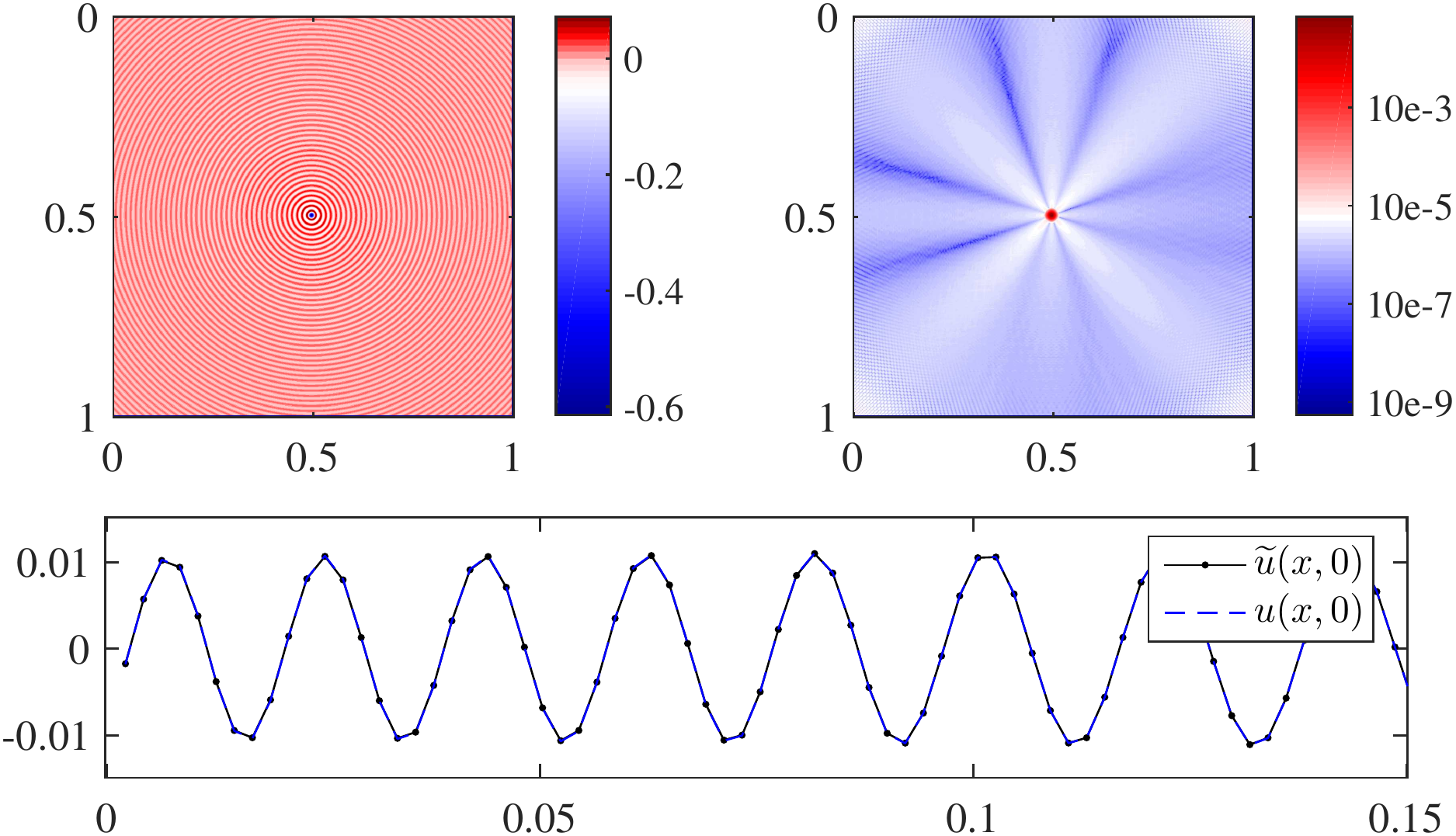}
	\caption{(Top-left), approximated solution $u(\x)$ (top-right) plot of  $|\widetilde{u}(\x)-u(\x)|$, (bottom) comparison on boundary values. Source at $\x=(0.5,0.5)$}\label{fig:hankel_c_500}
\end{figure}

\begin{figure}[ht]
		\includegraphics[scale=0.8]{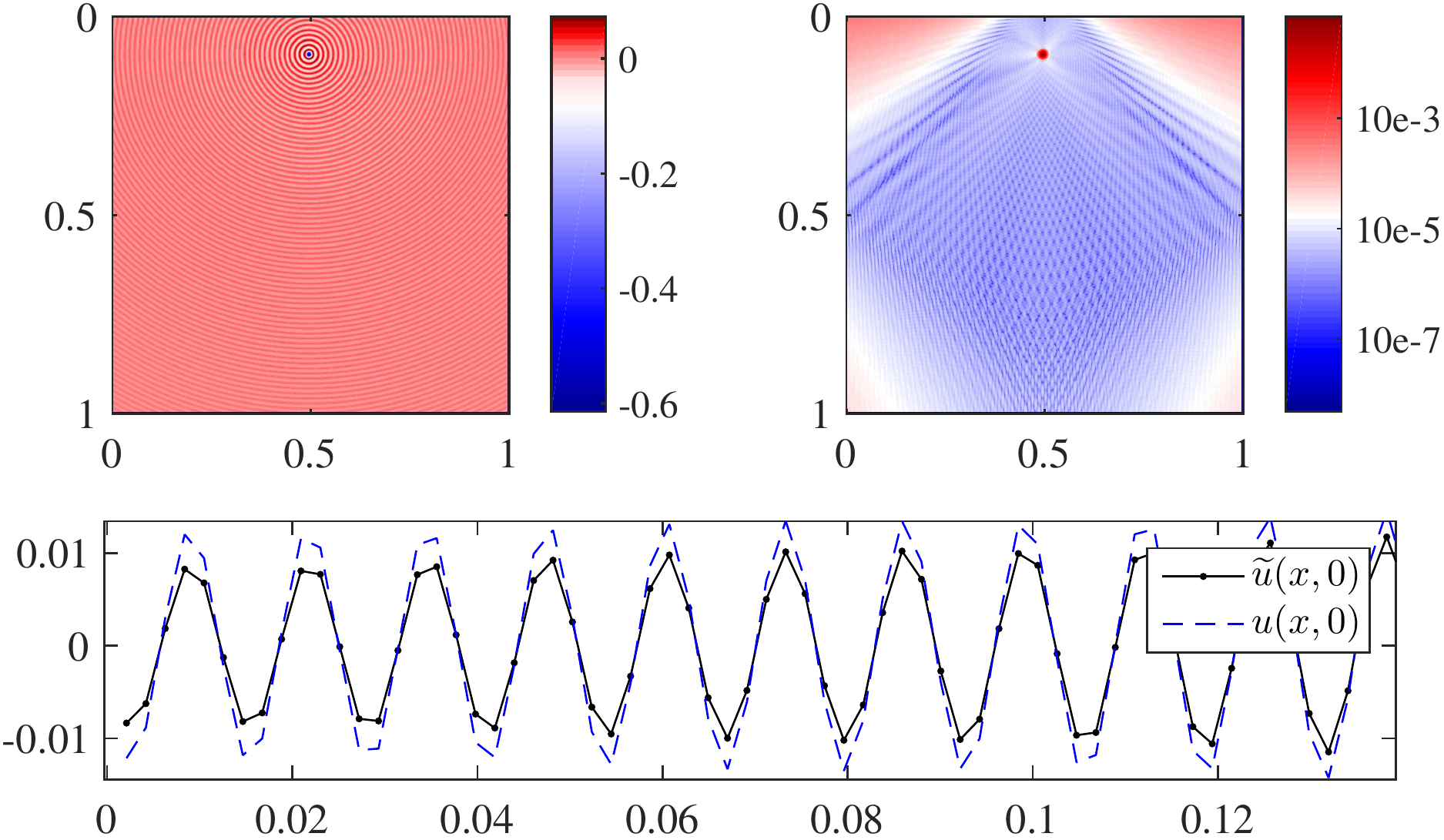}
		\caption{(Top-left), approximated solution $u(\x)$ (top-right) plot of  $|\widetilde{u}(\x)-u(\x)|$, (bottom) comparison on boundary values. Source at $\x=(0.5,0.1)$}\label{fig:hankel_m_500}
\end{figure}

We compute approximated Green's functions in the free space truncated to a bounded domain $\Omega\subset\R^2$ for $\x\neq\x_0$, where $\x_0\in\Omega$,  through the boundary problem
 \begin{equation}\label{eq:helmholtz_impedance}
 \left\{
 \begin{array}{rcll}
 -\Delta u(\x)-k^2u(\x)&=&\widetilde{\delta}(\x-\x_0), &\mbox{ in } \Omega\\
 \frac{\partial}{\partial\n}u(\x)+ik \left(1+\frac{3}{4k^2}\frac{\partial^2}{\partial\bm{\tau}^2}-\frac{i}{4k^3}\frac{\partial^3 }{\partial \mathbf{n}\partial\bm{\tau}^2 }\right) u(\x)&=&0, & \mbox{ on } \Gamma=\partial\Omega,
 \end{array}
 \right.
 \end{equation}
 where the boundary condition corresponds to the ABC in the Padé approximation. 
  The single source is given by the Gaussian function
\begin{equation}\label{eq:single_source}
\widetilde{\delta}(\x-\x_0)=\frac{1}{2\pi\sigma^2}e^{-\frac{\|\x-\x_0\|}{2\sigma^2}}
\end{equation}
for a value of $\sigma$ such that $\int_{\Omega}\widetilde{\delta}(\x-\x_0)d\x\approx1$.
We have used the BRFB-FD scheme with a square regular grid in $\Omega=(0,1)\times(0,1)$ with square $9$-stencils at inner nodes and $19$-stencils at boundary nodes. The results
shown in Fig. \ref{fig:hankel_c_500} and Fig. \ref{fig:hankel_m_500} where calculated with wavenumber $k=500$ with NPW$=6$, i.e., $h=2\pi/6k$. We point out that results show a good accuracy at $\x\neq\x_0$ and on the boundary the wavelength of the numerical solution matches very good with the exact one when the source is put at the center of the square domain, this is a good indication that dispersion errors are not significant, however in Fig. \ref{fig:hankel_m_500} (bottom) we see that amplitude has a considerable discrepancy respect to the exact one, this is due to the approximated ABC.


\subsection{Heterogeneous medium}

\subsubsection{Smooth medium}
\begin{figure}[ht]
	\begin{center}
		\begin{tabular}{ll}		
			\includegraphics[scale=0.45]{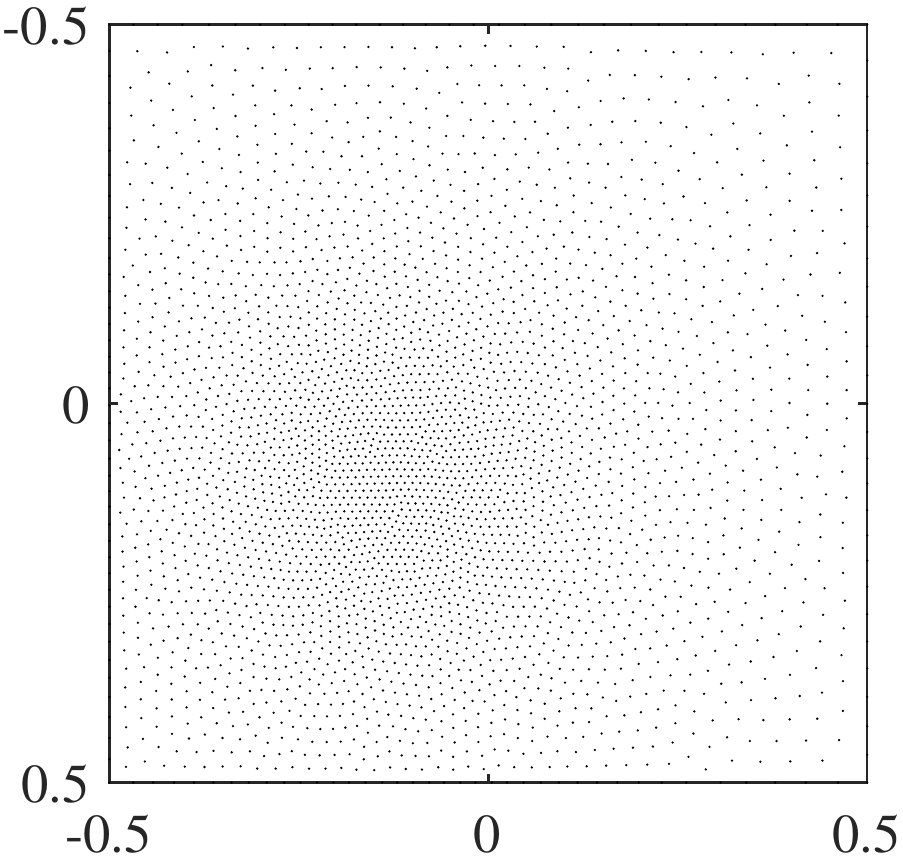} &\includegraphics[scale=0.6]{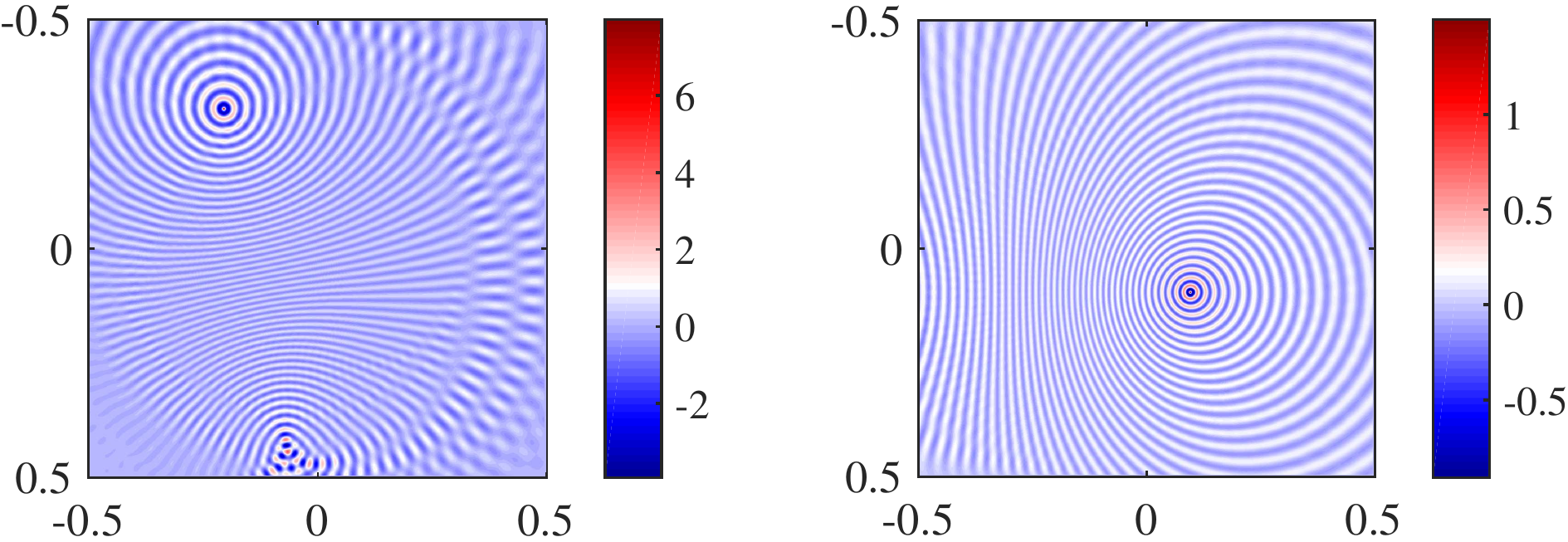} \\		
			\includegraphics[scale=0.45]{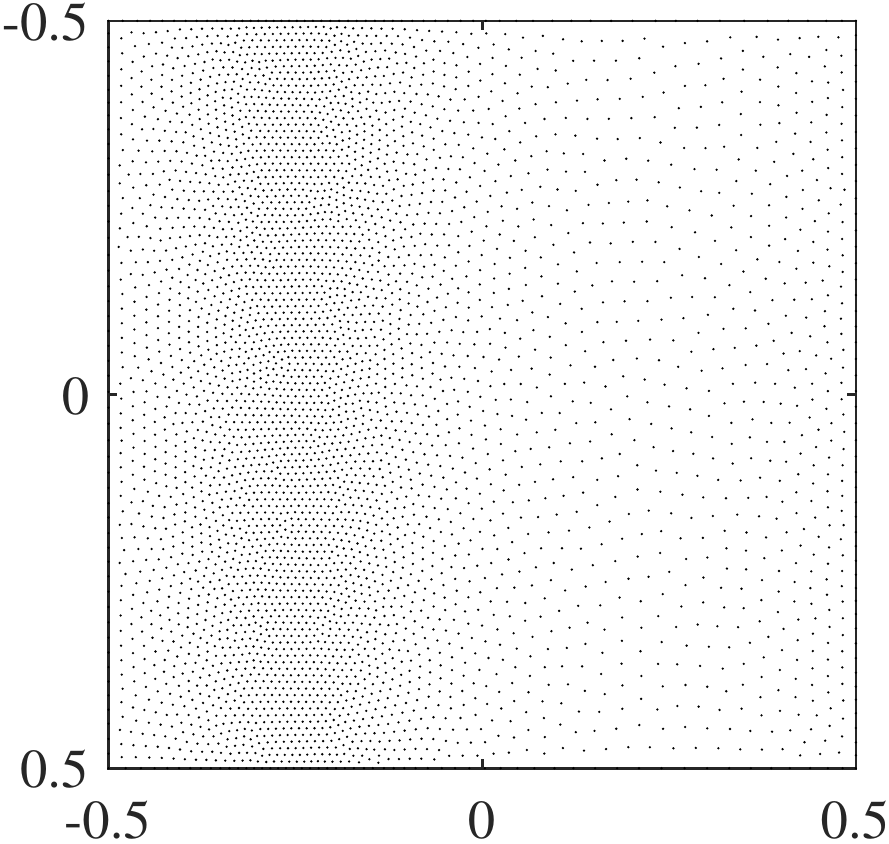} &\includegraphics[scale=0.6]{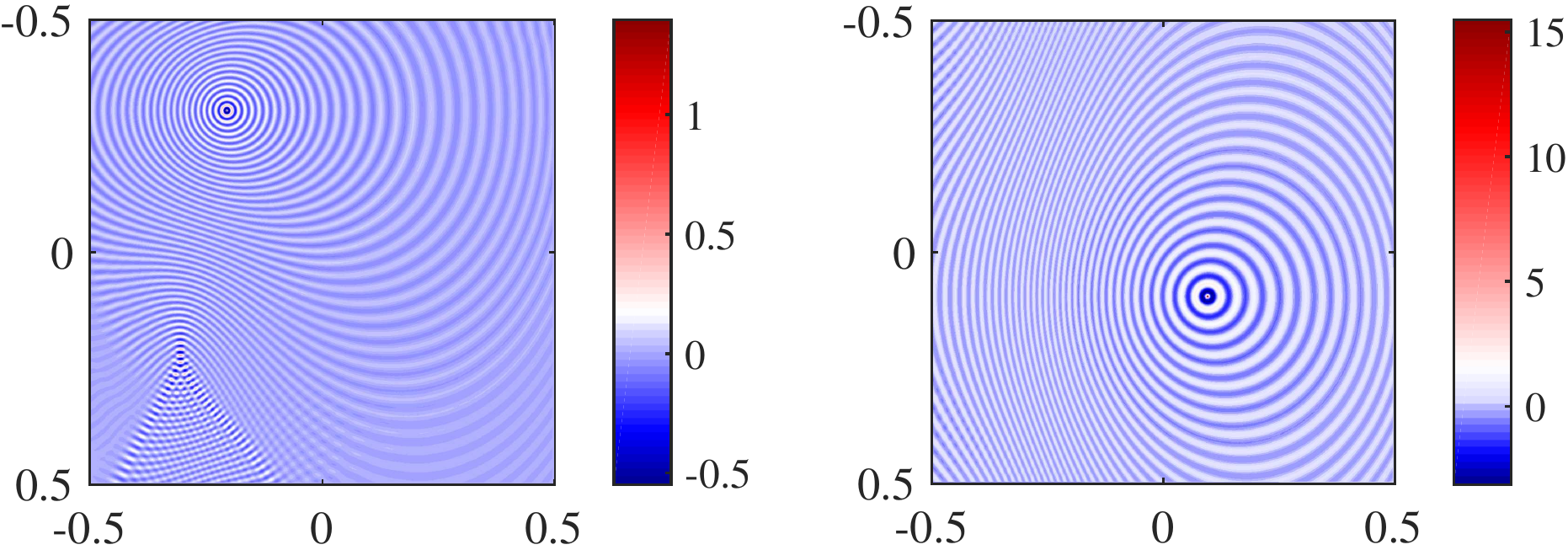} \\
		\end{tabular} \caption{Plot of solutions corresponding to velocity models in \eqref{eq:speed_fang3} and \eqref{eq:speed_fang4}. (Left column) Sketch of their node distribution according to the local wavelength and (right columns) wave fields corresponding to two different single sources.}
		\label{fig:fang34}
	\end{center}
\end{figure}

For this qualitative test, we have calculated approximated solutions of the problem
\begin{equation}\label{eq:helmholtz_ABC2}
\left\{
\begin{array}{rcll}
-\Delta u(\x)-\omega^2c(\x)^{-2}u(\x)&=&\widetilde{\delta}(\x-\x_0), &\mbox{ in } \Omega\\
\frac{\partial}{\partial\n}u(\x)+i\omega c(\x)^{-1}\mathcal{B}u(\x)&=&0, & \mbox{ on } \Gamma=\partial\Omega
\end{array}
\right.
\end{equation}
where $\mathcal{B}$ is  the operator $\mathcal{B}=1+\frac{c(\x)^2}{2\omega^2}\frac{\partial^2}{\partial\bm{\tau}^2}$ corresponding to the ABC of second order and $\Omega=(-0.5,0.5)\times(-0.5,0.5)$. Here we perform two examples with the velocity functions
\begin{equation}
c(x,y)=3 - 2.5e^{  -( (x+0.125)^2 + (y-0.1)^2 )/0.8^2  }
\label{eq:speed_fang3}
\end{equation} 
and
\begin{equation}
c(x,y)=1+0.5\sin(2\pi x).
\label{eq:speed_fang4}
\end{equation}
For these velocity models, nodes distributions are sketched on the left column of Fig. \ref{fig:fang34}. on center and right columns it can be seen the real part of the approximated solution for two different single sources. In table \ref{table:fang34} it shows results for required times to assembly sparse matrices $\mathbf{H}=\mathbf{H}_{\Omega}+\mathbf{H}_{\Gamma}$ and for solution of the system $-\mathbf{H}U=\mathbf{F}$ by LU factorization.

\begin{table}[ht]
	\begin{center}
		\begin{tabular}{ccc|ccc}
			\rowcolor{gray!20} $\omega/2\pi$  & Nodes ($N$)  &  $\kappa(\mathbf{H})$  &Time (s) for $\mathbf{H}_{\Omega}$ & Time (s) for $\mathbf{H}_{\partial\Omega}$ & Time (s) for LU  \\
			\hline\hline 
			2.5 & 1072 & 1.09e+04 & 0.37 & 0.24 & 0.02 \\
			5 & 4404 & 1.46e+05 & 1.18 & 0.45 & 0.10  \\
			10 & 17563 & 6.25e+05 & 4.72 & 1.01 & 0.45  \\
			20 & 70585 & 1.93e+06 & 18.98 & 2.58 & 3.00 \\
			40 & 283458 & 1.27e+07 & 75.10 & 9.05 & 17.14   \\
			\hline
		\end{tabular}
		\\
	\end{center} 
	\caption{Results in computing solutions corresponding to the smooth velocity model in \eqref{eq:speed_fang3}.} \label{table:fang34}		
\end{table}

\subsubsection{Non smooth medium (Test in the 2004 BP velocity model)}

\begin{figure}[ht]
	\begin{center}
		\begin{tabular}{c}
			\includegraphics[scale=0.48]{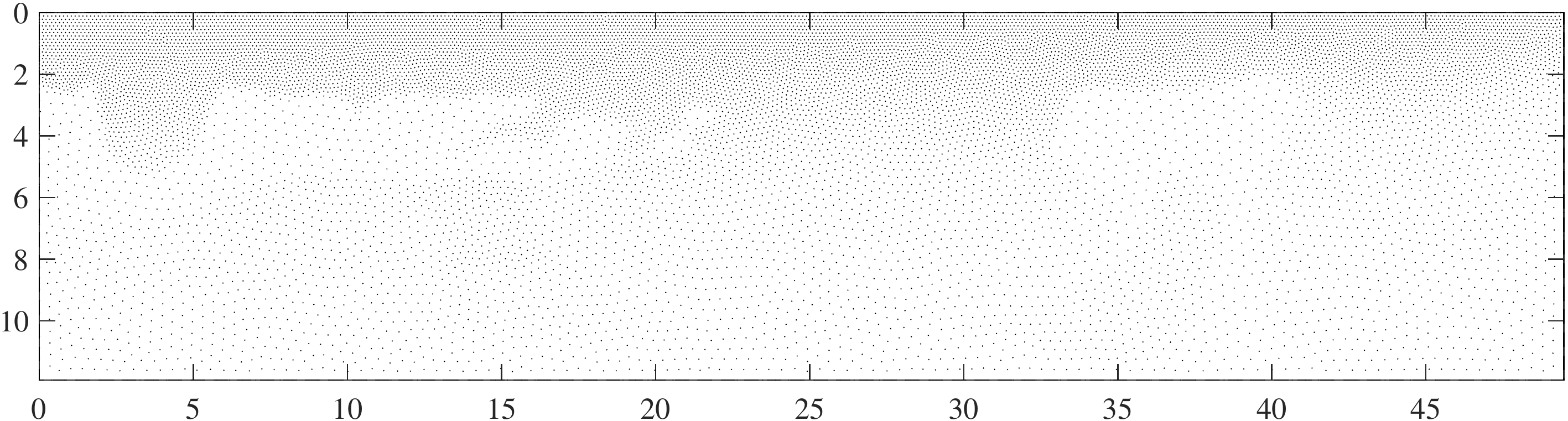}\\
			\includegraphics[scale=0.90]{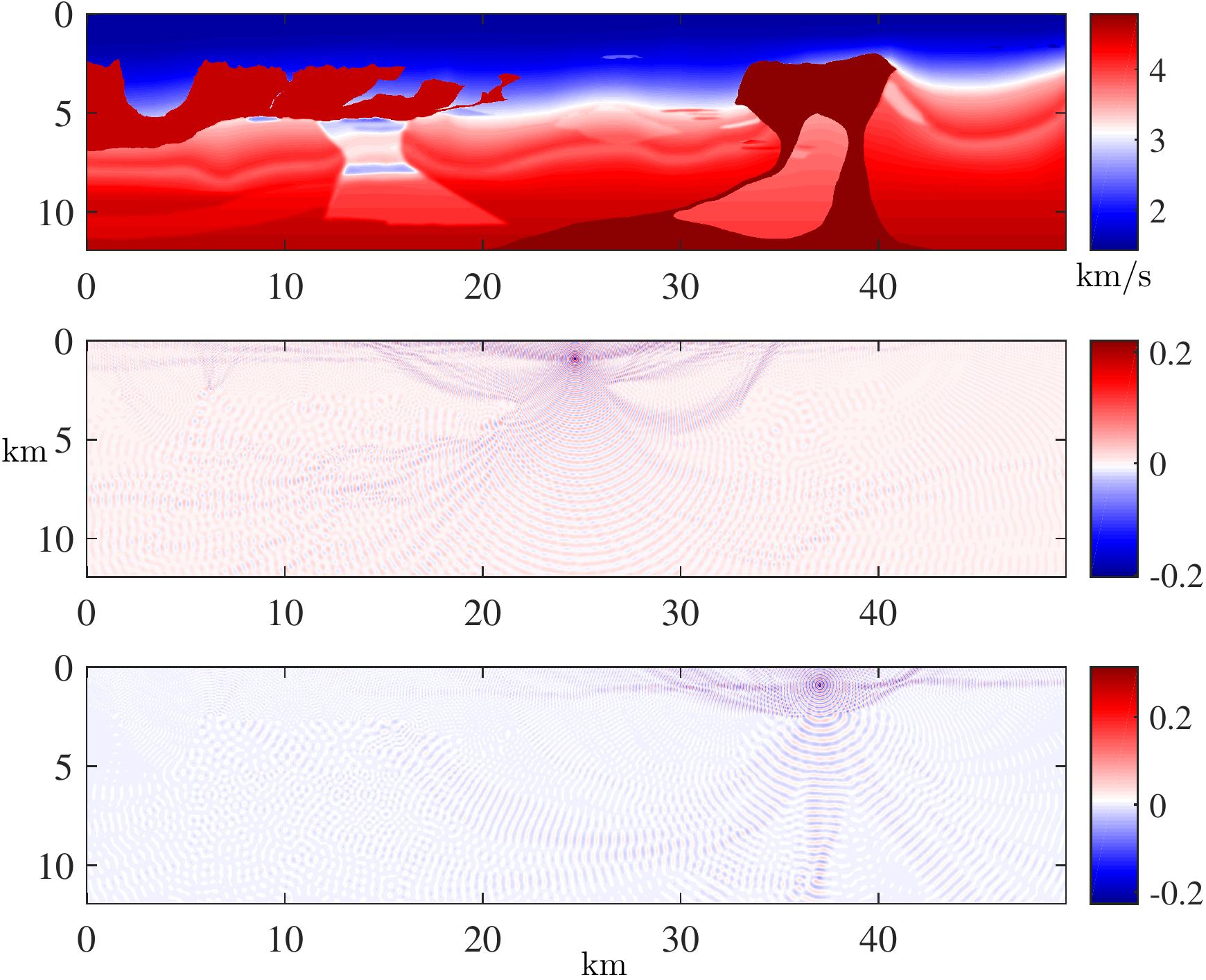}
		\end{tabular}
	\end{center}
	\caption{(Top) Sketch of node distribution for 2004BP model used for a frequency $\omega/2\pi=2$Hz with NPW=10 of the (Middle-top) velocity model.  (Middle-Bottom and bottom) plots of real part of the wavefield with to two single sources located at different positions for a frequency  $\omega/2\pi=10$Hz. In according to the local wavelength, the domain is discretized with $N=345393$, at inner nodes the local interpolation is performed with $19$-stencils and at boundary nodes with $25$-stencils.}\label{fig:2004bp_model}. 
\end{figure}

We consider the 2004 BP velocity benchmark, which is a popular model in research for velocity estimation methods in seismic imaging, which is presented as a challenge due to its complexity and large scale \cite{2004bp}. The  velocity function $c$ can be seen in the density plot in the middle-top in Fig. \ref{fig:2004bp_model}. Roughness of the velocity model such as hard interfaces and
sharp transitions generates strong reflections that hinder the efficiency of known iterative methods due to the ill-conditioning of the matrix gets worse \cite{ZEPEDA_NUNEZ2016347}, increasing the number of iterations. Moreover, for large $\omega$ the interaction of high frequency waves with short wavelength structures such as discontinuities, increases the reflections, further deteriorating the convergence rate. In BRBF-FD scheme the local interpolation matrices $\JJ_k$ becomes dramatically ill-conditioned, however by keeping the condition number of $\widetilde{\JJ}_k$ to the value $\kappa_0\approx 10^{(1+\sqrt{n_s})}$, we have empirically observed that the condition number of $\mathbf{H}$ lies in the range: $10^4\leq\kappa(\mathbf{H})\leq 10^9$, so it is feasible to perform LU factorization. Moreover, qualitatively, the wavelenghts of the wavefield has the expected behavior according to velocity model. To model wave propagation we have solved \eqref{eq:helmholtz_ABC2} with ABC of second order.

\begin{table}[ht]
	\begin{center}
		\begin{tabular}{ccc|ccc}
			\rowcolor{gray!20} $\omega/2\pi$  & Nodes ($N$)  &  $\kappa(\mathbf{H})$  &Time (s) for $\mathbf{H}_{\Omega}$ & Time (s) for $\mathbf{H}_{\partial\Omega}$ & Time (s) for LU  \\
			\hline\hline 
		2 & 38488  &  9.01e+04 & 20.59 & 4.72 & 0.99 \\
		4 & 153223 & 3.43e+05 & 74.27 & 10.52 & 5.42\\	
		6 & 345389 & 9.93e+05 & 173.78 & 19.75 & 16.04\\
		8 & 613529 & 7.22e+06 & 192.71 & 33.50 & 35.76\\	
			\hline
		\end{tabular}
		\\
	\end{center} 
	\caption{Results in computing solutions corresponding to the 2004 BP velocity model.} \label{table:2004bp}		
\end{table}

\section{Conclusions}
We perform local interpolation with a free shape parameter oscillatory RBF based on Bessel function of first kind to obtain a high order RBF-FD scheme for solving Helmholtz equation. In this approach the shape parameter is substituted for the local wavenumber. However due to the local interpolation matrices are extremely ill-conditioned even for small stencils, we overcome this issue with regularization by small perturbation of the diagonal. In some tests we have achieved convergence rates of third and sixth order. 

\section*{Acknowledgments}
This work is supported by Colombian Oil Company ECOPETROL and COLCIENCIAS as a part of the research project grant No. 0266-2013.
    
    \section*{References}
	\bibliographystyle{apalike}
	\bibliography{references.bib}

	
	
	
	
	

\end{document}

